\documentclass{article}
\usepackage{amsmath}
\usepackage{amssymb}
\usepackage{amsthm}
\begin{document}
\def\e#1\e{\begin{equation}#1\end{equation}}
\def\eq#1{{\rm(\ref{#1})}}
\theoremstyle{plain}
\newtheorem{thm}{Theorem}[section]
\newtheorem{lem}[thm]{Lemma}
\newtheorem{prop}[thm]{Proposition}
\newtheorem{cor}[thm]{Corollary}
\theoremstyle{definition}
\newtheorem{dfn}[thm]{Definition}
\newtheorem{ex}[thm]{Example}
\def\dim{\mathop{\rm dim}}
\def\Re{\mathop{\rm Re}}
\def\Im{\mathop{\rm Im}}
\def\supp{\mathop{\rm supp}}
\def\sn{{\ts\mathop{\rm sn}}}
\def\vol{\mathop{\rm vol}}
\def\area{\mathop{\rm area}}
\def\U{\mathbin{\rm U}}
\def\SU{\mathop{\rm SU}}
\def\SO{\mathop{\rm SO}}
\def\ge{\geqslant} 
\def\le{\leqslant} 
\def\R{\mathbin{\mathbb R}}
\def\Z{\mathbin{\mathbb Z}}
\def\C{\mathbin{\mathbb C}}
\def\al{\alpha}
\def\be{\beta}
\def\la{\lambda}
\def\ga{\gamma}
\def\de{\delta}
\def\ep{\epsilon}
\def\om{\omega}
\def\ka{\kappa}
\def\th{\theta}
\def\vp{\varphi}
\def\si{\sigma}
\def\ze{\zeta}
\def\De{\Delta}
\def\La{\Lambda}
\def\Om{\Omega}
\def\Si{\Sigma}
\def\Th{\Theta}
\def\Up{\Upsilon}
\def\d{{\rm d}}
\def\pd{\partial}
\def\db{{\bar\partial}}
\def\ts{\textstyle}
\def\sst{\scriptscriptstyle}
\def\bs{\boldsymbol}
\def\w{\wedge}
\def\lt{\ltimes}
\def\sm{\setminus}
\def\ot{\otimes}
\def\bigot{\bigotimes}
\def\iy{\infty}
\def\ra{\rightarrow}
\def\t{\times}
\def\ha{{\textstyle\frac{1}{2}}}
\def\ti{\tilde}
\def\ov{\overline}
\def\ovB{\,\overline{\!B}}
\def\ms#1{\vert#1\vert^2}
\def\md#1{\vert #1 \vert}
\def\bmd#1{\big\vert #1 \big\vert}
\def\cnm#1#2{\Vert #1 \Vert_{C^{#2}}} 
\def\lnm#1#2{\Vert #1 \Vert_{L^{#2}}} 
\title{$\U(1)$-invariant special Lagrangian 3-folds. III.\\
Properties of singular solutions}
\author{Dominic Joyce \\ Lincoln College, Oxford}
\date{}
\maketitle

\section{Introduction}
\label{up1}

Special Lagrangian submanifolds (SL $m$-folds) are a distinguished
class of real $m$-dimensional minimal submanifolds in $\C^m$, which are
calibrated with respect to the $m$-form $\Re(\d z_1\w\cdots\w\d z_m)$.
They can also be defined in Calabi--Yau manifolds, are important in
String Theory, and are expected to play a r\^ole in the eventual
explanation of Mirror Symmetry between Calabi--Yau 3-folds.

This is the third in a suite of three papers \cite{Joyc4,Joyc5}
studying special Lagrangian 3-folds $N$ in $\C^3$ invariant under
the $\U(1)$-action
\e
{\rm e}^{i\th}:(z_1,z_2,z_3)\mapsto
({\rm e}^{i\th}z_1,{\rm e}^{-i\th}z_2,z_3)
\quad\text{for ${\rm e}^{i\th}\in\U(1)$.}
\label{up1eq1}
\e
These three papers and \cite{Joyc6} are reviewed in \cite{Joyc7}.
Locally we can write $N$ as
\e
\begin{split}
N=\Bigl\{(z_1&,z_2,z_3)\in\C^3:
\Im(z_3)=u\bigl(\Re(z_3),\Im(z_1z_2)\bigr),\\
&\Re(z_1z_2)=v\bigl(\Re(z_3),\Im(z_1z_2)\bigr),
\quad\ms{z_1}-\ms{z_2}=2a\Bigr\},
\end{split}
\label{up1eq2}
\e
where $a\in\R$ and $u,v:\R^2\ra\R$ are continuous functions. It was
shown in \cite{Joyc4} that when $a\ne 0$, $N$ is an SL 3-fold
in $\C^3$ if and only if $u,v$ satisfy
\e
\frac{\pd u}{\pd x}=\frac{\pd v}{\pd y}\quad\text{and}\quad
\frac{\pd v}{\pd x}=-2\bigl(v^2+y^2+a^2\bigr)^{1/2}\frac{\pd u}{\pd y},
\label{up1eq3}
\e
and then $u,v$ are smooth and $N$ is nonsingular.

Here is what we mean by saying $N$ is locally of the form
\eq{up1eq2}. {\it Every} connected $\U(1)$-invariant SL 3-fold
$N$ may be written in the form \eq{up1eq2} for $u,v$ continuous,
{\it multi-valued\/} `functions' (multifunctions). Since
\eq{up1eq3} is a {\it nonlinear Cauchy-Riemann equation},
$u+iv$ is a bit like a holomorphic function of $x+iy$, and so
the multifunctions $(u,v)$ behave like holomorphic multifunctions
in complex analysis, such as~$\sqrt{z}$.

Thus we expect them to have isolated {\it branch points} in $\R^2$.
Over simply connected open sets $U\subset\R^2$ not containing any
branch points, the multifunctions $(u,v)$ decompose into {\it
sheets}, each of which is a continuous, single-valued function.
This is like choosing a branch of $\sqrt{z}$ on a simply-connected
open subset $U\subset\C\sm\{0\}$. In terms of $N$, we expect there
to be a discrete set of `branch point' $\U(1)$-orbits, and small
$\U(1)$-invariant open sets in $N$ away from these can be written
in the form \eq{up1eq2} for single-valued~$(u,v)$.

As $\frac{\pd u}{\pd x}=\frac{\pd v}{\pd y}$ there exists
$f:\R^2\ra\R$ with $\frac{\pd f}{\pd y}=u$ and
$\frac{\pd f}{\pd x}=v$, satisfying
\e
\Bigl(\Bigl(\frac{\pd f}{\pd x}\Bigr)^2+y^2+a^2
\Bigr)^{-1/2}\frac{\pd^2f}{\pd x^2}+2\,\frac{\pd^2f}{\pd y^2}=0.
\label{up1eq4}
\e
In \cite[Th.~7.6]{Joyc4} we proved existence and uniqueness for
the Dirichlet problem for \eq{up1eq4} in strictly convex domains
$S$ in $\R^2$ when $a\ne 0$. This gives existence and uniqueness
of a large class of nonsingular $\U(1)$-invariant SL 3-folds in
$\C^3$ satisfying certain boundary conditions.

When $a=0$, if $v=y=0$ the factor $-2(v^2+y^2+a^2)^{1/2}$ in
\eq{up1eq3} becomes zero, and then \eq{up1eq3} is no longer
elliptic. Because of this, when $a=0$ the appropriate thing to
do is to consider {\it weak solutions} of \eq{up1eq3}, which
may have {\it singular points} $(x,0)$ with $v(x,0)=0$. At
such a point $u,v$ may not be differentiable, and
$\bigl(0,0,x+iu(x,0)\bigr)$ is a singular point of the SL
3-fold $N$ in~$\C^3$. 

In \cite[Th.~7.1]{Joyc5} we proved existence and uniqueness for
a suitable class of weak solutions to the Dirichlet problem for
\eq{up1eq4} in strictly convex domains $S$ in $\R^2$ when $a=0$.
This gives existence and uniqueness of a large class of {\it singular}
$\U(1)$-invariant SL 3-folds in $\C^3$ satisfying certain boundary
conditions.

The goal of this paper is to study these singular SL 3-folds
in more detail. It will be shown that under mild conditions the
singularities are {\it isolated}, and can be roughly classified
by a {\it multiplicity} $n\ge 1$, and within each multiplicity
by one of two {\it types}. For singularities of multiplicity
$n$, the `germ' of the singularity is described to leading
order by real parameters~$\ga_1,\ldots,\ga_n$.

Examples of singular SL 3-folds will be constructed with every
multiplicity $n\ge 1$ and type, and realizing all possible values
of $\ga_1,\ldots,\ga_n$. They provide an infinite family of local
models for singularities of compact SL 3-folds in (almost)
Calabi--Yau 3-folds. We will also use our results to construct
large families of {\it special Lagrangian fibrations} of open
subsets of $\C^3$, including singular fibres, which may be of
every multiplicity and type. These will be studied further
in \cite{Joyc6}, in connection with the SYZ Conjecture. For
a brief summary of \cite{Joyc4,Joyc5,Joyc6} and this paper,
see~\cite{Joyc7}.

Section \ref{up2} introduces special Lagrangian geometry, and
\S\ref{up3} gives some background in analysis and Geometric
Measure Theory. Section \ref{up4} then reviews the preceding
papers \cite{Joyc4,Joyc5}. The new material begins in \S\ref{up5},
where we classify the possible {\it tangent cones}, in the sense
of Geometric Measure Theory, at a singular point of an SL 3-fold
of the form~\eq{up1eq2}.

Let $(u,v)$ and $(\hat u,\hat v)$ satisfy \eq{up1eq3} in a domain
$S$ in $\R^2$, for $a\ne 0$. In \cite{Joyc4} we defined the
{\it multiplicity} of a zero of $(u,v)-(\hat u,\hat v)$, and
related the number of zeroes of $(u,v)-(\hat u,\hat v)$ in
$S^\circ$, counted with multiplicity, to boundary data on $\pd S$.
Section \ref{up6} extends these to the singular case $a=0$,
showing that zeroes of $(u,v)-(\hat u,\hat v)$ are {\it isolated},
with multiplicity a {\it positive integer}.

As part of the proof, in \S\ref{up62} we study singular solutions
$(u,v)$ of \eq{up1eq3} with $a=0$ and $v(x,0)\equiv 0$, so that
$u,v$ is singular all along the $x$-axis in $S$. We show that
such solutions have the symmetry $u(x,-y)=u(x,y)$, $v(x,-y)=-v(x,y)$,
and the corresponding SL 3-folds $N$ are actually the union of two
{\it nonsingular} $\U(1)$-invariant SL 3-folds $N_\pm$ intersecting
in a real analytic real curve $\ga$, which is the singular set of
$N$ and the fixed point set of the $\U(1)$-action~\eq{up1eq1}.

Section \ref{up7} applies the results of \S\ref{up6} to construct
{\it special Lagrangian fibrations} on open subsets of $\C^3$, which
will be the central tool in \cite{Joyc6}. Let $(u,v)$ be a singular
solution of \eq{up1eq3} in $S$ with $a=0$. In \S\ref{up8} we show
that either $u(x,-y)\equiv u(x,y)$ and $v(x,-y)\equiv -v(x,y)$, or
the singularities of $u,v$ in $S^\circ$ are {\it isolated}.

For isolated singularities in $S^\circ$ we define the {\it
multiplicity} $n\ge 1$ and {\it type}, discuss the tangent cones
of the corresponding SL 3-fold singularities, and give counting
formulae for singular points with multiplicity. Section \ref{up9}
shows that singularities with every multiplicity $n\ge 1$ and
type exist, and occur in codimension $n$ in the family of all
$\U(1)$-invariant SL 3-folds.
\medskip

\noindent{\it Acknowledgements.} I would like to thank Robert
Bryant for kindly proving Theorem \ref{up6thm2} below when I could
not, and for taking the trouble to write up the proof \cite{Brya}.
I would also like to thank Mark Gross for helpful conversations,
and the referee for useful criticism. I was supported by an EPSRC
Advanced Fellowship whilst writing this paper.

\section{Special Lagrangian submanifolds in $\C^m$}
\label{up2}

For introductions to special Lagrangian geometry, see the author
\cite{Joyc3} and Harvey and Lawson \cite[\S III]{HaLa}. We begin
by defining {\it calibrations} and {\it calibrated submanifolds},
following Harvey and Lawson~\cite{HaLa}.

\begin{dfn} Let $(M,g)$ be a Riemannian manifold. An {\it oriented
tangent $k$-plane} $V$ on $M$ is a vector subspace $V$ of
some tangent space $T_xM$ to $M$ with $\dim V=k$, equipped
with an orientation. If $V$ is an oriented tangent $k$-plane
on $M$ then combining $g\vert_V$ with the orientation on $V$
gives a natural {\it volume form} $\vol_V$ on $V$, which is a 
$k$-form on~$V$.

Now let $\vp$ be a closed $k$-form on $M$. We call $\vp$ a
{\it calibration} on $M$ if for every oriented
$k$-plane $V$ on $M$ we have $\vp\vert_V\le \vol_V$. Here
$\vp\vert_V=\al\cdot\vol_V$ for some $\al\in\R$, and 
$\vp\vert_V\le\vol_V$ if $\al\le 1$. Let $N$ be an 
oriented submanifold of $M$ with dimension $k$. Then 
each tangent space $T_xN$ for $x\in N$ is an oriented
tangent $k$-plane. We call $N$ a {\it calibrated 
submanifold} if $\vp\vert_{T_xN}=\vol_{T_xN}$ for all~$x\in N$.
\label{up2def1}
\end{dfn}

It is easy to show that calibrated submanifolds are automatically
{\it minimal submanifolds} \cite[Th.~II.4.2]{HaLa}. Here is the 
definition of special Lagrangian submanifolds in $\C^m$, taken
from~\cite[\S III]{HaLa}.

\begin{dfn} Let $\C^m$ have complex coordinates $(z_1,\dots,z_m)$, 
and define a metric $g$, a real 2-form $\om$ and a complex $m$-form 
$\Om$ on $\C^m$ by
\e
\begin{split}
g=\ms{\d z_1}+\cdots+\ms{\d z_m},\quad
\om&=\frac{i}{2}(\d z_1\w\d\bar z_1+\cdots+\d z_m\w\d\bar z_m),\\
\text{and}\quad\Om&=\d z_1\w\cdots\w\d z_m.
\end{split}
\label{up2eq1}
\e
Then $\Re\Om$ and $\Im\Om$ are closed real $m$-forms on $\C^m$. Let
$L$ be an oriented real submanifold of $\C^m$ of real dimension 
$m$. We say that $L$ is a {\it special Lagrangian submanifold} 
of $\C^m,$ or {\it SL\/ $m$-fold}\/ for short, if $L$ is calibrated 
with respect to $\Re\Om$, in the sense of Definition~\ref{up2def1}. 
\label{up2def2}
\end{dfn}

Harvey and Lawson \cite[Cor.~III.1.11]{HaLa} give the following
alternative characterization of special Lagrangian submanifolds.

\begin{prop} Let\/ $L$ be a real\/ $m$-dimensional submanifold 
of\/ $\C^m$. Then $L$ admits an orientation making it into a
special Lagrangian submanifold of\/ $\C^m$ if and only if\/
$\om\vert_L\equiv 0$ and\/~$\Im\Om\vert_L\equiv 0$.
\label{up2prop1}
\end{prop}

An $m$-dimensional submanifold $L$ in $\C^m$ is called {\it Lagrangian} 
if $\om\vert_L\equiv 0$. Thus special Lagrangian submanifolds are 
Lagrangian submanifolds satisfying the extra condition that 
$\Im\Om\vert_L\equiv 0$, which is how they get their name.

\section{Background material from analysis}
\label{up3}

Here are some definitions we will need to make sense of
analytic results from \cite{Joyc4,Joyc5}. A good reference for
\S\ref{up31}--\S\ref{up32} is Gilbarg and Trudinger~\cite{GiTr}.

\subsection{Banach spaces of functions on subsets of $\R^n$}
\label{up31}

We shall work in a special class of subsets of $\R^n$ called
{\it domains}.

\begin{dfn} A closed, bounded, contractible subset $S$ in $\R^n$ will
be called a {\it domain} if it is a disjoint union $S=S^\circ\cup\pd S$,
where the {\it interior} $S^\circ$ of $S$ is a connected open set in
$\R^n$ with $S=\overline{S^\circ}$, and the {\it boundary}
$\pd S=S\sm S^\circ$ is a compact embedded hypersurface in~$\R^n$.

A domain $S$ in $\R^2$ is called {\it strictly convex} if $S$ is 
convex and the curvature of $\pd S$ is nonzero at every point. So, 
for example, $x^2+y^2\le 1$ is strictly convex but $x^4+y^4\le 1$ is 
not, as its boundary has zero curvature at $(\pm 1,0)$ and~$(0,\pm 1)$.
\label{up3def1}
\end{dfn}

We will use a number of different spaces of real functions on~$S$.

\begin{dfn} Let $S$ be a domain in $\R^n$. For $k\ge 0$, define
$C^k(S)$ to be the space of continuous functions $f:S\ra\R$ with
$k$ continuous derivatives, and norm $\cnm{f}k=\sum_{j=0}^k\sup_S
\bmd{\pd^jf}$. Define $C^\iy(S)=\bigcap_{k=0}^\iy C^k(S)$. For
$k\ge 0$ and $\al\in(0,1]$, define the {\it H\"older space}
$C^{k,\al}(S)$ to be the subset of $f\in C^k(S)$ for which
\e
[\pd^k f]_\al=\sup_{x\ne y\in S}
\frac{\bmd{\pd^kf(x)-\pd^kf(y)}}{\md{x-y}^\al}
\label{up3eq1}
\e
is finite. The {\it H\"older norm} on $C^{k,\al}(S)$ is
$\cnm{f}{k,\al}=\cnm{f}k+[\pd^kf]_\al$. For $q\ge 1$, define the
{\it Lebesgue space} $L^p(S)$ to be the set of locally integrable
functions $f$ on $S$ for which the norm
\begin{equation*}
\lnm{f}p=\left(\int_S\md{f}^p\d{\bf x}\right)^{1/p}
\end{equation*}
is finite. Then $C^k(S),C^{k,\al}(S)$ and $L^p(S)$ are Banach spaces.
Here $\pd$ is the vector operator $(\frac{\pd}{\pd x_1},\ldots,
\frac{\pd}{\pd x_n})$, and lengths $\bmd{\pd^jf}$ are computed
using the standard Euclidean metric on~$\R^n$.
\label{up3def2}
\end{dfn}

\subsection{Quasilinear elliptic operators and weak solutions}
\label{up32}

Quasilinear elliptic operators are a class of nonlinear partial
differential operators, which are linear and elliptic in their
highest-order derivatives.

\begin{dfn} Let $S$ be a domain in $\R^n$. A {\it second-order 
quasilinear operator} $Q:C^2(S)\ra C^0(S)$ is an operator of the form
\e
\bigl(Qu\bigr)(x)=
\sum_{i,j=1}^na^{ij}(x,u,\pd u)\frac{\pd^2u}{\pd x_i\pd x_j}(x)
+b(x,u,\pd u),
\label{up3eq2}
\e
where $a^{ij}$ and $b$ are continuous maps $S\t\R\t(\R^n)^*\ra\R$,
and $a^{ij}=a^{ji}$ for all $i,j=1,\ldots,n$. We call the functions
$a^{ij}$ and $b$ the {\it coefficients} of $Q$. We call $Q$ 
{\it elliptic} if the symmetric $n\t n$ matrix $(a^{ij})$ is 
positive definite at every point of~$S\t\R\t(\R^n)^*$. 
\label{up3def3}
\end{dfn}

A second-order quasilinear operator $Q:C^2(S)\ra C^0(S)$ is in
{\it divergence form} if it is written as
\e
\bigl(Qu\bigr)(x)=
\sum_{j=1}^n\frac{\pd}{\pd x_j}\bigl(a^j(x,u,\pd u)\bigr)
+b(x,u,\pd u)
\label{up3eq3}
\e
for functions $a^j\in C^1\bigr(S\t\R\t(\R^n)^*\bigr)$ for $j=1,\ldots,n$ 
and $b\in C^0\bigr(S\t\R\t(\R^n)^*\bigr)$.

\begin{dfn} Let $Q$ be a second-order quasilinear operator
in divergence form on a domain $S$ in $\R^n$, given by
equation \eq{up3eq3}. Then we say that $u\in L^1(S)$ is a
{\it weak solution} of the equation $Qu=f$ for $f\in L^1(S)$ 
if $u$ is weakly differentiable with weak derivative $\pd u$, 
and $a^j(x,u,\pd u),b(x,u,\pd u)\in L^1(S)$ with
\e
-\sum_{j=1}^n\int_S\frac{\pd\psi}{\pd x_j}\cdot a^j(x,u,\pd u)
\d{\bf x}+\int_S\psi\cdot b(x,u,\pd u)\d{\bf x}
=\int_S\psi\cdot f\,\d{\bf x}
\label{up3eq4}
\e
for all $\psi\in C^1(S)$ with~$\psi\vert_{\pd S}\equiv 0$.
\label{up3def4}
\end{dfn}

If $Q$ is a second-order quasilinear operator, we shall interpret
the equation $Qu=f$ in three different senses:
\begin{itemize}
\setlength{\itemsep}{0pt}
\setlength{\parsep}{0pt}
\item We say that $Qu=f$, or $Qu=f$ {\it holds classically}, if $u\in
C^2(S)$ and $f\in C^0(S)$ and $Qu=f$ in $C^0(S)$ in the usual way.
\item We say that $Qu=f$ {\it holds with weak derivatives} if $u$
is twice weakly differentiable and $Qu=f$ holds in $L^p(S)$ for
some $p\ge 1$, defining $Qu$ using the weak derivatives of $u$.
\item We say that $Qu=f$ {\it holds weakly} if $u$ is a weak
solution of $Qu=f$, as in Definition \ref{up3def4}. Note that
this requires only that $u$ be {\it once} weakly differentiable,
and the second derivatives of $u$ need not exist even weakly.
\end{itemize}

Clearly the first sense implies the second, which implies the third.
If $Q$ is {\it elliptic} and $a^j,b,f$ are suitably regular, one can
usually show that a weak solution to $Qu=f$ is a classical solution,
so that the three senses are equivalent. But we will be dealing with
singular equations that are not elliptic at every point, and then the
three senses will be distinct.

\subsection{Geometric Measure Theory and tangent cones}
\label{up33}

Much of this paper concerns the singularities of SL 3-folds,
which are examples of {\it singular minimal submanifolds}. Now
there is an elegant theory of singular submanifolds called
{\it Geometric Measure Theory}, which is especially powerful
in dealing with singular minimal submanifolds. An introduction
to the subject is provided by Morgan \cite{Morg} and an in-depth
treatment by Federer \cite{Fede}, and Harvey and Lawson \cite[\S
II]{HaLa} relate Geometric Measure Theory to calibrated geometry.

Let $m\le n$ be nonnegative integers. One defines a class of
$m$-dimensional {\it rectifiable currents} in $\R^n$, which are
measure-theoretic generalizations of compact, oriented $m$-submanifolds
$N$ with boundary $\pd N$ in $\R^n$, with integer {\it multiplicities}.
Here $N$ with multiplicity $k$ is like $k$ copies of $N$ superimposed,
and changing the orientation of $N$ changes the sign of the multiplicity.
This enables us to add and subtract submanifolds.

If $T$ is an $m$-dimensional rectifiable current, one can define
the {\it volume} $\vol(T)$ of $T$, by Hausdorff $m$-measure. If
$\vp$ is a compactly-supported $m$-form on $\R^n$ then one can
define $\int_T\vp$. Thus we can consider $T$ as an element
$\vp\mapsto\int_T\vp$ of the dual space $({\cal D}^m)^*$ of
the vector space ${\cal D}^m$ of compactly-supported $m$-forms
on $\R^n$. This induces a topology on the space of rectifiable
currents in $\R^n$. The {\it interior} $T^\circ$ of $T$ is
$T\sm\pd T$ (that is,~$\supp T\sm\supp\pd T$).

An $m$-dimensional rectifiable current $T$ is called an {\it integral
current} if $\pd T$ is an $(m\!-\!1)$-dimensional rectifiable current.
By \cite[5.5]{Morg}, \cite[4.2.17]{Fede}, integral currents have strong
{\it compactness properties}.

\begin{thm} Let\/ $m\le n$, $C>0$, and\/ $K$ be a closed ball in $\R^n$.
Then the set of\/ $m$-dimensional integral currents $T$ in $K$ with\/
$\vol(T)\le C$ and\/ $\vol(\pd T)\le C$ is compact in an appropriate
weak topology.
\label{up3thm1}
\end{thm}

This is significant as it enables us to easily construct
{\it minimal\/} (volume-minimizing) integral currents satisfying
certain conditions, by choosing a suitable minimizing sequence
and extracting a convergent subsequence. The important question
then becomes, how close are such minimal integral currents to
being manifolds, and how bad are their singularitities? One can
(partially) answer this using regularity results. Here, for example,
is a major result of Almgren \cite[8.3]{Morg}, slightly rewritten.

\begin{thm} Let\/ $T$ be an $m$-dimensional minimal rectifiable
current in $\R^n$. Then $T^\circ$ is a smooth, embedded,
oriented, minimal\/ $m$-submanifold with locally constant
positive integer multiplicity, except on a singular set of
Hausdorff dimension at most\/~$m-2$.
\label{up3thm2}
\end{thm}

Next we discuss {\it tangent cones} of minimal rectifiable
currents, a generalization of tangent spaces of submanifolds,
as in~\cite[9.7]{Morg}.

\begin{dfn} A locally rectifiable current $C$ in $\R^n$ is
called a {\it cone} if $C=tC$ for all $t>0$, where $t:\R^n\ra\R^n$
acts by dilations in the obvious way. Let $T$ be a rectifiable
current in $\R^n$, and let ${\bf x}\in T^\circ$. We say
that $C$ is a {\it tangent cone} to $T$ at $\bf x$ if there
exists a decreasing sequence $r_1>r_2>\cdots$ such that
$r_j\ra 0$ and $r_j^{-1}(T-{\bf x})\ra C$ as a locally
rectifiable current as~$j\ra\iy$.
\label{up3def5}
\end{dfn}

The next result follows from Morgan \cite[p.~94-95]{Morg}, Federer
\cite[5.4.3]{Fede} and Harvey and Lawson~\cite[Th.~II.5.15]{HaLa}.

\begin{thm} Let\/ $T$ be a minimal rectifiable current in $\R^n$.
Then for all\/ ${\bf x}\in T^\circ$, there exists a tangent
cone $C$ to $T$ at\/ $\bf x$. Moreover $C$ is itself a minimal
locally rectifiable current with\/ $\pd C=\emptyset$, and if\/
$T$ is calibrated with respect to a constant calibration $\vp$
on $\R^n$, then $C$ is also calibrated with respect to~$\vp$.
\label{up3thm3}
\end{thm}

Note that the theorem does {\it not\/} claim that the tangent cone
$C$ is unique, and in fact it is an important open question whether
a minimal rectifiable current has a unique tangent cone at each
point of $T^\circ$, \cite[p.~93]{Morg}. However, using the idea
of {\it density} we can constrain the choice of tangent cones.

\begin{dfn} Let $T$ be an $m$-dimensional minimal locally
rectifiable current in $\R^n$. For each ${\bf x}\in T^\circ$
define the {\it density} $\Th(T,{\bf x})$ of $T$ at $\bf x$ by
\e
\Th(T,{\bf x})=\lim_{r\ra 0_+}
\frac{\vol\bigl(T\cap\ovB_r({\bf x})\bigr)}{\om_mr^m},
\label{up3eq5}
\e
where $\ovB_r({\bf x})$ is the closed ball of radius $r$ about
${\bf x}\in\R^n$, and $\vol(\ldots)$ the volume of $m$-dimensional
rectifiable currents, and $\om_m$ the volume of the unit ball
in $\R^m$. By \cite[9.4]{Morg} the limit in \eq{up3eq5} exists for
all ${\bf x}\in T^\circ$. By \cite[5.4.5(1)]{Fede} the density is
an {\it upper-semicontinuous} function on~$T^\circ$.
\label{up3def6}
\end{dfn}

Note that if $C$ is an $m$-dimensional locally rectifiable cone in
$\R^n$, and ${\cal S}^{n-1}$ is the unit sphere in $\R^n$, then
$C\cap{\cal S}^{n-1}$ is an $(m-1)$-dimensional rectifiable current
in ${\cal S}^{n-1}$, and $\Th(C,0)=\vol(C\cap{\cal S}^{n-1})/
\vol({\cal S}^{m-1})$. By \cite[9.9]{Morg}, the density at
${\bf x}\in T^\circ$ of a minimal rectifiable current $T$
agrees with the density at $0$ of any tangent cone $C$ to $T$
at~$\bf x$.

\begin{thm} Let\/ $T$ be a minimal rectifiable current in $\R^n$,
and\/ ${\bf x}\in T^\circ$. Then each tangent cone $C$ to $T$
at\/ $\bf x$ has~$\Th(C,0)=\Th(T,{\bf x})$.
\label{up3thm4}
\end{thm}

Therefore all tangent cones $C$ to $T$ at $\bf x$ must have the same
density at 0. A multiplicity 1 minimal $m$-submanifold $T$ in $\R^n$
has unique tangent cone $T_{\bf x}T$ at $\bf x$ and density 1 at every
point. Here is a kind of converse to this,~\cite[5.4.7]{Fede}.

\begin{thm} Let\/ $T$ be an $m$-dimensional minimal rectifiable
current in $\R^n$. Then $\Th(T,{\bf x})\ge 1$ for all\/
$x\in T^\circ$. There exists $\Up>1$ depending only on $m,n$
such that if\/ ${\bf x}\in T^\circ$ and\/ $\Th(T,{\bf x})<\Up$
then $\Th(T,{\bf x})=1$, and\/ $T$ is a smooth, embedded,
multiplicity $1$ minimal\/ $m$-submanifold near~$\bf x$.
\label{up3thm5}
\end{thm}

The way this works is that the only tangent cones $C$ with density
1 at 0 are multiplicity 1 subspaces $\R^m$ in $\R^n$, and if a
minimal rectifiable current has this as a tangent cone at $\bf x$,
then it is a submanifold near~$\bf x$.

\section{Review of material from \cite{Joyc4} and \cite{Joyc5}}
\label{up4}

We now recapitulate those results from \cite{Joyc4} and \cite{Joyc5}
that we will need later. Readers are referred to \cite{Joyc4,Joyc5}
for proofs, discussion and motivation. The material of \S\ref{up46}
is new.

\subsection{Finding the equations}
\label{up41}

The following result \cite[Prop.~4.1]{Joyc4} is the starting
point for everything in \cite{Joyc4,Joyc5} and this paper.

\begin{prop} Let\/ $S$ be a domain in $\R^2$ or $S=\R^2$, let\/
$u,v:S\ra\R$ be continuous, and\/ $a\in\R$. Define
\e
\begin{split}
N=\bigl\{(z_1,z_2,z_3)\in\C^3:\,& z_1z_2=v(x,y)+iy,\quad z_3=x+iu(x,y),\\
&\ms{z_1}-\ms{z_2}=2a,\quad (x,y)\in S\bigr\}.
\end{split}
\label{up4eq1}
\e
Then 
\begin{itemize}
\setlength{\itemsep}{0pt}
\setlength{\parsep}{0pt}
\item[{\rm(a)}] If\/ $a=0$, then $N$ is a (possibly singular)
special Lagrangian $3$-fold in $\C^3$, with boundary over
$\pd S$, if\/ $u,v$ are differentiable and satisfy
\e
\frac{\pd u}{\pd x}=\frac{\pd v}{\pd y}
\quad\text{and}\quad
\frac{\pd v}{\pd x}=-2\bigl(v^2+y^2\bigr)^{1/2}\frac{\pd u}{\pd y},
\label{up4eq2}
\e
except at points $(x,0)$ in $S$ with\/ $v(x,0)=0$, where $u,v$ 
need not be differentiable. The singular points of\/ $N$ are those
of the form $(0,0,z_3)$, where $z_3=x+iu(x,0)$ for $x\in\R$ 
with\/~$v(x,0)=0$.
\item[{\rm(b)}] If\/ $a\ne 0$, then $N$ is a nonsingular SL\/
$3$-fold in $\C^3$, with boundary over $\pd S$, if and only if\/
$u,v$ are differentiable on all of\/ $S$ and satisfy
\e
\frac{\pd u}{\pd x}=\frac{\pd v}{\pd y}\quad\text{and}\quad
\frac{\pd v}{\pd x}=-2\bigl(v^2+y^2+a^2\bigr)^{1/2}\frac{\pd u}{\pd y}.
\label{up4eq3}
\e
\end{itemize}
\label{up4prop1}
\end{prop}

By showing that $u,v$ each satisfy second-order elliptic equations
and using the maximum principle, we prove~\cite[Cor.~4.4]{Joyc4}:

\begin{cor} Let\/ $S$ be a domain in $\R^2$, let\/ $a\ne 0$,
and suppose $u,v\in C^1(S)$ satisfy \eq{up4eq3}. Then the
maxima and minima of\/ $u$ and\/ $v$ are achieved on~$\pd S$.
\label{up4cor}
\end{cor}

\subsection{Examples}
\label{up42}

Here are some examples of SL 3-folds $N$ in the form \eq{up4eq1},
and the corresponding functions $u,v$. Let $a\ge 0$, and define
\e
\begin{split}
N_a=\Bigl\{(z_1&,z_2,z_3)\in\C^3:\ms{z_1}-2a=\ms{z_2}=\ms{z_3},\\
&\Im\bigl(z_1z_2z_3\bigr)=0,\quad \Re\bigl(z_1z_2z_3\bigr)\ge 0\Bigr\}.
\end{split}
\label{up4eq4}
\e
Then $N_a$ is a nonsingular SL 3-fold diffeomorphic to
${\cal S}^1\t\R^2$ when $a>0$, and $N_0$ is an SL $T^2$-cone
with one singular point at $(0,0,0)$. The $N_a$ are invariant
under the $\U(1)^2$-action
\e
({\rm e}^{i\th_1},{\rm e}^{i\th_2}):(z_1,z_2,z_3)\longmapsto
({\rm e}^{i\th_1}z_1,{\rm e}^{i\th_2}z_2,{\rm e}^{-i\th_1-i\th_2}z_3),
\label{up4eq5}
\e
and are part of a family of explicit $\U(1)^2$-invariant SL 3-folds
written down by Harvey and Lawson \cite[\S III.3.A]{HaLa}. By
\cite[Th.~5.1]{Joyc4}, these SL 3-folds can be written in the
form~\eq{up4eq1}.

\begin{thm} Let\/ $a\ge 0$. Then there exist unique $u_a,v_a:\R^2\ra\R$
such that the SL $3$-fold\/ $N_a$ of \eq{up4eq4} agrees with\/ $N$ in
\eq{up4eq1} with\/ $u_a,v_a$ in place of $u,v$ and\/ $S=\R^2$. Furthermore:
\begin{itemize}
\setlength{\itemsep}{0pt}
\setlength{\parsep}{0pt}
\item[{\rm(a)}] $u_a,v_a$ are smooth on $\R^2$ and satisfy \eq{up4eq3}, 
except at\/ $(0,0)$ when $a=0$, where they are only continuous.
\item[{\rm(b)}] $u_a(x,y)<0$ when $y>0$ for all\/ $x$, and\/ $u_a(x,0)=0$
for all\/ $x$, and\/ $u_a(x,y)>0$ when $y<0$ for all\/~$x$.
\item[{\rm(c)}] $v_a(x,y)>0$ when $x>0$ for all\/ $y$, and\/ $v_a(0,y)=0$
for all\/ $y$, and\/ $v_a(x,y)<0$ when $x<0$ for all\/~$y$.
\item[{\rm(d)}] $u_a(0,y)=-y\bigl(\md{a}+\sqrt{y^2+a^2}\,\,\bigr)^{-1/2}$ 
for all\/~$y$.
\item[{\rm(e)}] $v_a(x,0)=x\bigl(x^2+2\md{a}\bigr)^{1/2}$ for all\/~$x$.
\end{itemize}
\label{up4thm1}
\end{thm}

In \cite[Ex.~5.2 \& Ex.~5.4]{Joyc4} we give two further examples:

\begin{ex} Let $\al,\be,\ga\in\R$ and define $u(x,y)=\al x+\be$ 
and $v(x,y)=\al y+\ga$. Then $u,v$ satisfy \eq{up4eq2} for any 
value of~$a$. 
\label{up4ex1}
\end{ex}

\begin{ex} Let $S=\R^2$, $u(x,y)=\md{y}-\ha\cosh 2x$ and
$v(x,y)=-y\sinh 2x$. Then $u,v$ satisfy \eq{up4eq2}, except that
$\frac{\pd u}{\pd y}$ is not well-defined on the $x$-axis. So
equation \eq{up4eq1} with $a=0$ defines an explicit special
Lagrangian 3-fold $N$ in $\C^3$. It turns out that $N$ is the
union of two nonsingular SL 3-folds intersecting in a real
curve, which are constructed in \cite[Ex.~7.4]{Joyc2} by
evolving paraboloids in~$\C^3$.
\label{up4ex2}
\end{ex}

\subsection{Results motivated by complex analysis}
\label{up43}

Section 6 of \cite{Joyc4} proves analogues for solutions of
\eq{up4eq3} of results on the zeroes of holomorphic functions,
using the ideas of {\it winding number} and {\it multiplicity},
defined in~\cite[Def.~6.1 \& Def.~6.3]{Joyc4}.

\begin{dfn} Let $C$ be a compact oriented 1-manifold, and
$\ga:C\ra\R^2\sm\{0\}$ a differentiable map. Then the
{\it winding number of\/ $\ga$ about\/ $0$ along} $C$ is
$\frac{1}{2\pi}\int_C\ga^*(\d\th)$, where $\d\th$ is the
closed 1-form $(x\,\d y-y\,\d x)/(x^2+y^2)$ on~$\R^2\sm\{0\}$.

In fact the winding number is simply the {\it topological
degree} of $\ga$. Thus it is actually well-defined for
$\ga$ only {\it continuous}, and is invariant under
{\it continuous deformations} of $\ga$, which will be
important in \S\ref{up6}. We gave the $\d\th$ definition for
differentiable $\ga$ first only for the sake of explicitness.
\label{up4def1}
\end{dfn}

\begin{dfn} Let $S$ be a domain in $\R^2$, let $a\ne 0$, and 
suppose $(u_1,v_1)$ and $(u_2,v_2)$ are solutions of \eq{up4eq3} in
$C^1(S)$. Let $k\ge 1$ be an integer and $(b,c)\in S^\circ$. We say 
that $(u_1,v_1)-(u_2,v_2)$ {\it has a zero of multiplicity $k$ at\/} 
$(b,c)$ if $\pd^ju_1(b,c)=\pd^ju_2(b,c)$ and $\pd^jv_1(b,c)=\pd^j
v_2(b,c)$ for $j=0,\ldots,k-1$, but $\pd^ku_1(b,c)\ne\pd^ku_2(b,c)$ 
and $\pd^kv_1(b,c)\ne\pd^kv_2(b,c)$. If $(u_1,v_1)\not\equiv(u_2,v_2)$
then every zero of $(u_1,v_1)-(u_2,v_2)$ has a unique multiplicity.
\label{up4def2}
\end{dfn}

In \cite[Prop.~6.5]{Joyc4} we show that zeroes of $(u_1,v_1)-(u_2,v_2)$
resemble zeroes of holomorphic functions to leading order.

\begin{prop} Let\/ $S$ be a domain in $\R^2$, let\/ $a\ne 0$, and 
let\/ $(u_1,v_1)$ and\/ $(u_2,v_2)$ be solutions of\/ \eq{up4eq3} in
$C^1(S)$. Suppose $(u_1,v_1)-(u_2,v_2)$ has a zero of multiplicity 
$k\ge 1$ at\/ $(b,c)$ in $S^\circ$. Then there exists a nonzero 
complex number $C$ such that
\e
\begin{split}
\la u_1(x,y)+iv_1(x,y)=\la u_2(x,y)&+iv_2(x,y)
+C\bigl(\la(x-b)+i(y-c)\bigr)^k\\
&+O\bigl(\md{x-b}^{k+1}+\md{y-c}^{k+1}\bigr),
\end{split}
\label{up4eq6}
\e
where~$\la=\sqrt{2}\bigl(v_1(b,c)^2+c^2+a^2\bigr)^{1/4}$.
\label{up4prop2}
\end{prop}

In \cite[Th.~6.7]{Joyc4} we give a formula for the number of
zeroes of~$(u_1,v_1)-(u_2,v_2)$.

\begin{thm} Let\/ $S$ be a domain in $\R^2$ and\/ $(u_1,v_1)$,
$(u_2,v_2)$ solutions of\/ \eq{up4eq3} in $C^1(S)$ for some
$a\ne 0$, with\/ $(u_1,v_1)\neq(u_2,v_2)$ at every point of\/
$\pd S$. Then $(u_1,v_1)-(u_2,v_2)$ has finitely many zeroes in
$S$. Let there be $n$ zeroes, with multiplicities $k_1,\ldots,k_n$.
Then the winding number of\/ $(u_1,v_1)-(u_2,v_2)$ about\/ $0$
along $\pd S$ is~$\sum_{i=1}^nk_i$.
\label{up4thm2}
\end{thm}

\subsection{Generating $u,v$ from a potential $f$}
\label{up44}

In \cite[Prop.~7.1]{Joyc4} we show that solutions $u,v\in C^1(S)$ 
of \eq{up4eq3} come from a potential $f\in C^2(S)$ satisfying a
second-order quasilinear elliptic equation.

\begin{prop} Let\/ $S$ be a domain in $\R^2$ and\/ $u,v\in C^1(S)$
satisfy \eq{up4eq3} for $a\ne 0$. Then there exists $f\in C^2(S)$
with\/ $\frac{\pd f}{\pd y}=u$, $\frac{\pd f}{\pd x}=v$ and
\e
P(f)=\Bigl(\Bigl(\frac{\pd f}{\pd x}\Bigr)^2+y^2+a^2
\Bigr)^{-1/2}\frac{\pd^2f}{\pd x^2}+2\,\frac{\pd^2f}{\pd y^2}=0.
\label{up4eq7}
\e
This $f$ is unique up to addition of a constant, $f\mapsto f+c$.
Conversely, all solutions of\/ \eq{up4eq7} yield solutions 
of\/~\eq{up4eq3}. 
\label{up4prop3}
\end{prop}

Equation \eq{up4eq7} may also be written in {\it divergence form} as
\e
P(f)=\frac{\pd}{\pd x}\Bigl[A\Bigl(a,y,\frac{\pd f}{\pd x}\Bigr)\Bigr]
+2\,\frac{\pd^2f}{\pd y^2}=0,
\label{up4eq8}
\e
where $A(a,y,v)=\int_0^v\bigl(w^2+y^2+a^2\bigr)^{-1/2}\,\d w$, so that
$\frac{\pd A}{\pd v}=\bigl(v^2+y^2+a^2\bigr)^{-1/2}$. Note that $A$
is undefined when $a=y=0$.

In \cite[Th.~7.6]{Joyc4} we prove existence and uniqueness of 
solutions to the Dirichlet problem for \eq{up4eq7} in {\it strictly 
convex domains} in $\R^2$, as in Definition~\ref{up3def1}.

\begin{thm} Let\/ $S$ be a strictly convex domain in $\R^2$, 
and let\/ $a\ne 0$, $k\ge 0$ and\/ $\al\in(0,1)$. Then for 
each\/ $\phi\in C^{k+2,\al}(\pd S)$ there exists a unique 
solution $f$ of\/ \eq{up4eq7} in $C^{k+2,\al}(S)$ with\/ 
$f\vert_{\pd S}=\phi$. This $f$ is real analytic in $S^\circ$, 
and satisfies $\cnm{f}{1}\le C\cnm{\phi}{2}$, for some $C>0$
depending only on~$S$.
\label{up4thm3}
\end{thm}

After considerable work, this was extended to the case $a=0$
in~\cite[Th.~7.1]{Joyc5}.

\begin{thm} Let\/ $S$ be a strictly convex domain in $\R^2$ invariant
under the involution $(x,y)\mapsto(x,-y)$, let\/ $k\ge 0$ and\/ $\al\in
(0,1)$. Then for each\/ $\phi\in C^{k+3,\al}(\pd S)$ there exists a
unique weak solution $f$ of\/ \eq{up4eq8} in $C^1(S)$ with\/ $a=0$ and\/
$f\vert_{\pd S}=\phi$. Furthermore $f$ is twice weakly differentiable
and satisfies \eq{up4eq7} with weak derivatives.

Let\/ $u=\frac{\pd f}{\pd y}$ and\/ $v=\frac{\pd f}{\pd x}$. Then
$u,v\in C^0(S)$ are weakly differentiable and satisfy \eq{up4eq2}
with weak derivatives, and\/ $v$ satisfies \eq{up4eq9} weakly with\/
$a=0$. The weak derivatives $\frac{\pd u}{\pd x},\frac{\pd u}{\pd y},
\frac{\pd v}{\pd x},\frac{\pd v}{\pd y}$ satisfy $\frac{\pd u}{\pd x}=
\frac{\pd v}{\pd y}\in L^p(S)$ for $p\in[1,2]$, and\/ $\frac{\pd u}{\pd
y}\in L^q(S)$ for $q\in[1,2)$, and\/ $\frac{\pd v}{\pd x}$ is bounded
on $S$. Also $u,v$ are $C^{k+2,\al}$ in $S$ and real analytic in
$S^\circ$ except at singular points $(x,0)$ with\/~$v(x,0)=0$.
\label{up4thm4}
\end{thm}

Combined with Propositions \ref{up4prop1} and \ref{up4prop3},
these two theorems give an existence and uniqueness result for
$\U(1)$-invariant SL 3-folds in $\C^3$ satisfying certain boundary
conditions. In \cite[Th.~7.2]{Joyc5} we show that in the last two
theorems $f$ depends continuously on~$\phi,a$.

\begin{thm} Let\/ $S$ be a strictly convex domain in $\R^2$ invariant
under the involution $(x,y)\mapsto(x,-y)$, let\/ $k\ge 0$ and\/
$\al\in(0,1)$. Then the map $C^{k+3,\al}(\pd S)\t\R\ra C^1(S)$ taking
$(\phi,a)\mapsto f$ is continuous, where $f$ is the unique solution
of\/ \eq{up4eq7} (with weak derivatives) with\/ $f\vert_{\pd S}=\phi$
constructed in Theorem \ref{up4thm3} when $a\ne 0$, and in Theorem
\ref{up4thm4} when $a=0$. This map is also continuous in stronger
topologies on $f$ than the $C^1$ topology.
\label{up4thm5}
\end{thm} 

In \cite[Th.~7.11]{Joyc4} we prove an analogue of Theorem
\ref{up4thm2} where we count zeroes of $(u_1,v_1)-(u_2,v_2)$ not
in terms of a winding number, but in terms of the stationary
points of the difference of potentials $\phi_1-\phi_2$ on~$\pd S$.

\begin{thm} Let\/ $S$ be a strictly convex domain in $\R^2$,
let\/ $a\ne 0$, $\al\in(0,1)$, and\/ $f_1,f_2\in C^{2,\al}(S)$
satisfy \eq{up4eq7} with\/ $f_j\vert_{\pd S}=\phi_j$. Set\/
$u_j=\frac{\pd f_j}{\pd y}$ and\/ $v_j=\frac{\pd f_j}{\pd x}$,
so that\/ $u_j,v_j\in C^{1,\al}(S)$ satisfy \eq{up4eq3}.
Suppose $\phi_1-\phi_2$ has exactly $l$ local maxima and\/
$l$ local minima on $\pd S$. Then $(u_1,v_1)-(u_2,v_2)$ has
$n$ zeroes in $S^\circ$ with multiplicities $k_1,\ldots,k_n$,
where~$\sum_{i=1}^nk_i\le l-1$.
\label{up4thm6}
\end{thm}

\subsection{An elliptic equation satisfied by $v$}
\label{up45}

In \cite[Prop.~8.1]{Joyc4} we show that if $u,v$ satisfy \eq{up4eq3}
then $v$ satisfies a second-order quasilinear elliptic equation, and
conversely, any solution $v$ of this equation extends to a solution 
$u,v$ of~\eq{up4eq3}.

\begin{prop} Let\/ $S$ be a domain in $\R^2$ and\/ $u,v\in C^2(S)$
satisfy \eq{up4eq3} for $a\ne 0$. Then
\e
Q(v)=\frac{\pd}{\pd x}\Bigl[\bigl(v^2+y^2+a^2\bigr)^{-1/2}
\frac{\pd v}{\pd x}\Bigr]+2\,\frac{\pd^2v}{\pd y^2}=0.
\label{up4eq9}
\e
Conversely, if\/ $v\in C^2(S)$ satisfies \eq{up4eq9} then 
there exists $u\in C^2(S)$, unique up to addition of a 
constant\/ $u\mapsto u+c$, such that $u,v$ satisfy~\eq{up4eq3}.
\label{up4prop4}
\end{prop}

In \cite[Prop 8.7]{Joyc4} we show that solutions of \eq{up4eq9}
satisfying strict inequalities on $\pd S$ satisfy the same
inequality on~$S$.

\begin{prop} Let\/ $S$ be a domain in $\R^2$, let\/ $a\ne 0$,
and suppose $v,v'\in C^2(S)$ satisfy \eq{up4eq9} on $S$. If\/
$v<v'$ on $\pd S$ then $v<v'$ on~$S$.
\label{up4prop5}
\end{prop}

In \cite[Th.~8.8]{Joyc4} we prove existence and uniqueness of 
solutions to the Dirichlet problem for \eq{up4eq9} when $a\ne 0$
in arbitrary domains in~$\R^2$.

\begin{thm} Let\/ $S$ be a domain in $\R^2$. Then whenever $a\ne 0$, 
$k\ge 0$, $\al\in(0,1)$ and\/ $\phi\in C^{k+2,\al}(\pd S)$ there 
exists a unique solution $v\in C^{k+2,\al}(S)$ of\/ \eq{up4eq9} 
with\/ $v\vert_{\pd S}=\phi$. Fix a basepoint\/ $(x_0,y_0)\in S$.
Then there exists a unique $u\in C^{k+2,\al}(S)$ with\/ $u(x_0,y_0)=0$
such that $u,v$ satisfy \eq{up4eq3}. Furthermore, $u,v$ are real
analytic in~$S^\circ$.
\label{up4thm7}
\end{thm}

This was extended to $a=0$ in \cite[Th.~6.1]{Joyc5}, for
a restricted class of domains $S$. Much of the technical work
in \cite{Joyc5} went into proving $u,v$ are continuous.

\begin{thm} Let\/ $S$ be a strictly convex domain in $\R^2$ invariant
under the involution $(x,y)\mapsto(x,-y)$, let\/ $k\ge 0$ and\/ $\al\in
(0,1)$. Suppose $\phi\in C^{k+2,\al}(\pd S)$ with\/ $\phi(x,0)\ne 0$
for points $(x,0)$ in $\pd S$. Then there exists a unique weak solution
$v$ of\/ \eq{up4eq9} in $C^0(S)$ with\/ $a=0$ and\/~$v\vert_{\pd S}=\phi$.

Fix a basepoint\/ $(x_0,y_0)\in S$. Then there exists a unique $u\in
C^0(S)$ with\/ $u(x_0,y_0)=0$ such that\/ $u,v$ are weakly differentiable
in $S$ and satisfy \eq{up4eq2} with weak derivatives. The weak derivatives
$\frac{\pd u}{\pd x},\frac{\pd u}{\pd y},\frac{\pd v}{\pd x},
\frac{\pd v}{\pd y}$ satisfy $\frac{\pd u}{\pd x}=\frac{\pd v}{\pd y}
\in L^p(S)$ for $p\in[1,\frac{5}{2})$, and\/ $\frac{\pd u}{\pd y}\in
L^q(S)$ for $q\in[1,2)$, and\/ $\frac{\pd v}{\pd x}$ is bounded on $S$.
Also $u,v$ are $C^{k+2,\al}$ in $S$ and real analytic in $S^\circ$
except at singular points $(x,0)$ with\/~$v(x,0)=0$.
\label{up4thm8}
\end{thm}

Combined with Proposition \ref{up4prop1} these two theorems give an
existence and uniqueness for nonsingular and singular $\U(1)$-invariant
SL 3-folds in $\C^3$ satisfying certain boundary conditions. In
\cite[Th.~6.2]{Joyc5} we show that in the last two theorems
$u,v$ depend continuously on~$\phi,a$.

\begin{thm} Let\/ $S$ be a strictly convex domain in $\R^2$ invariant
under the involution $(x,y)\mapsto(x,-y)$, let\/ $k\ge 0$, $\al\in(0,1)$,
and\/ $(x_0,y_0)\in S$. Define $X$ to be the set of\/ $\phi\in C^{k+2,\al}
(\pd S)$ with\/ $\phi(x,0)=0$ for some $(x,0)\in\pd S$. Then the map
$C^{k+2,\al}(\pd S)\t\R\sm X\t\{0\}\ra C^0(S)^2$ taking $(\phi,a)
\mapsto(u,v)$ is continuous, where $(u,v)$ is the unique solution of\/
\eq{up4eq3} (with weak derivatives when $a=0$) with\/ $v\vert_{\pd S}=\phi$
and\/ $u(x_0,y_0)=0$, constructed in Theorem \ref{up4thm7} when $a\ne 0$,
and in Theorem \ref{up4thm8} when $a=0$. This map is also continuous in
stronger topologies on $(u,v)$ than the $C^0$ topology.
\label{up4thm9}
\end{thm}

\subsection{A class of solutions of \eq{up4eq2} with singularities}
\label{up46}

In most of the rest of the paper we shall be studying solutions
of \eq{up4eq2} with singularities, and the corresponding SL
3-folds $N$. So we need to know just what we mean by a {\it
singular solution} of \eq{up4eq2}. We give a definition here,
to avoid repeating technicalities about weak derivatives, and
so on, many times.

\begin{dfn} Let $S$ be a domain in $\R^2$ and $u,v\in C^0(S)$. We
say that $u,v$ are a {\it singular solution} of \eq{up4eq2} if
\begin{itemize}
\setlength{\itemsep}{0pt}
\setlength{\parsep}{0pt}
\item[{\rm(i)}] $u,v$ are weakly differentiable, and their weak
derivatives $\frac{\pd u}{\pd x},\frac{\pd u}{\pd y},
\frac{\pd v}{\pd x},\frac{\pd v}{\pd y}$ in $L^1(S)$
satisfy~\eq{up4eq2}.
\item[{\rm(ii)}] $v$ is a {\it weak solution} of \eq{up4eq9} with
$a=0$, as in~\S\ref{up32}.
\item[{\rm(iii)}] Define the {\it singular points} of $u,v$ to be
the $(x,0)\in S$ with $v(x,0)=0$. Then except at singular points,
$u,v$ are $C^2$ in $S$ and real analytic in $S^\circ$, and
satisfy \eq{up4eq2} in the classical sense.
\item[{\rm(iv)}] For $a\in(0,1]$ there exist $u_a,v_a\in C^2(S)$
satisfying \eq{up4eq3} such that as $a\ra 0_+$ we have
\end{itemize}
\e
\begin{gathered}
\ts\text{$u_a\ra u$, $v_a\ra v$ in $C^0(S)$,
$\frac{\pd u_a}{\pd x}\ra\frac{\pd u}{\pd x}$,
$\frac{\pd v_a}{\pd y}\ra\frac{\pd v}{\pd y}$ in $L^2(S)$,}\\
\ts\text{$\frac{\pd u_a}{\pd y}\ra\frac{\pd u}{\pd y}$ in $L^q(S)$,
$q\in[1,2)$, and $\frac{\pd v_a}{\pd x}\ra\frac{\pd v}{\pd x}$ in
$L^r(X)$, $r\in[1,\iy)$.}
\end{gathered}
\label{up4eq10}
\e
\label{up4def3}
\end{dfn}

The reason we choose this definition is that Theorems \ref{up4thm4}
and \ref{up4thm8} prove existence and uniqueness of singular solutions
$u,v$ of \eq{up4eq2} with boundary conditions on certain domains $S$.
Note that part (iv) holds automatically in these theorems because
$u,v$ were constructed as limits of $u_a,v_a$ satisfying \eq{up4eq10}
in the proofs in \cite{Joyc5}. Note also that if $u,v$ are singular
solutions in $S$, then the restrictions to any subdomain $T$ in $S$
are singular solutions in $T$. It is not difficult to show that
$u_0,v_0$ in Theorem \ref{up4thm1} and $u,v$ in Example \ref{up4ex2}
are both examples of singular solutions of \eq{up4eq2} in~$\R^2$.

We show that the singular SL 3-folds corresponding to singular
solutions $u,v$ are {\it special Lagrangian rectifiable currents},
in the sense of~\S\ref{up33}.

\begin{prop} Let\/ $S$ be a domain in $\R^2$ and\/ $u,v\in C^0(S)$
a singular solution of\/ \eq{up4eq2}, in the sense of Definition
\ref{up4def3}. Define $N$ by \eq{up4eq1} with\/ $a=0$. Then $N$ is
the support of a special Lagrangian rectifiable current, in the
sense of Geometric Measure Theory, which we identify with\/ $N$.
The boundary current\/ $\pd N$ is supported on
\e
\begin{split}
\bigl\{(z_1,z_2,z_3)\in\C^3:\,& z_1z_2=v(x,y)+iy,\quad z_3=x+iu(x,y),\\
&\md{z_1}=\md{z_2},\quad (x,y)\in \pd S\bigr\}.
\end{split}
\label{up4eq11}
\e
\label{up4prop6}
\end{prop}

\begin{proof} For $a\in(0,1]$ let $u_a,v_a$ be as in Definition
\ref{up4def3}(iv), and define $N_a$ by \eq{up4eq1} using $u_a,v_a$.
Then $N_a$ is a compact, nonsingular SL 3-fold with boundary, and
so defines a {\it special Lagrangian integral current}. We shall
show that $N_a\ra N$ as rectifiable currents as $a\ra 0_+$, so
that $N$ is an SL rectifiable current.

Calculation shows that the volume of $N_a$ is given by
\e
\vol(N_a)=\int_{N_a}\Re\Om=
\int_S\Bigl(1+\frac{\pd u_a}{\pd x}\,\frac{\pd v_a}{\pd y}-
\frac{\pd u_a}{\pd y}\,\frac{\pd v_a}{\pd x}\Bigr)\,\d x\,\d y.
\label{up4eq12}
\e
Fix $q\in[1,2)$ and $r\in[1,\iy)$ with $\frac{1}{q}+\frac{1}{r}=1$.
Then by \eq{up4eq10}, $\frac{\pd u_a}{\pd x},\frac{\pd v_a}{\pd y},
\frac{\pd u_a}{\pd y},\frac{\pd v_a}{\pd x}$ converge in
$L^2,L^2,L^q,L^r(S)$ respectively as $a\ra 0_+$. Hence $\vol(N_a)$
converges to the well-defined integral $\int_S\bigl(1+\frac{\pd
u}{\pd x}\frac{\pd v}{\pd y}-\frac{\pd u}{\pd y}\frac{\pd v}{\pd
x}\bigr)\,\d x\,\d y$ as $a\ra 0_+$, by H\"older's inequality.
In particular, $\vol(N_a)$ is uniformly bounded for small~$a$.

Using \eq{up4eq10} one can show that the family of rectifiable
currents $N_a$ is {\it Cauchy} in the flat norm on currents as
$a\ra 0_+$. But the family of {\it integral flat chains} supported
in a given compact subset of $\R^n$ is {\it complete} in the flat
norm. Thus $N_a\ra N$ in the flat norm as $a\ra 0_+$ for a unique
integral flat chain~$N$.

As $\vol(N_a)$ is uniformly bounded for small $a$ we have
$\vol(N)<\iy$. So by the Closure Theorem \cite[5.4]{Morg}, $N$
is a {\it rectifiable current}. Since $N_a\ra N$ as currents,
and $N_a$ is special Lagrangian for all $a$, we see that $N$ is
{\it special Lagrangian}. And as $u_a\ra u$, $v_a\ra v$ in $C^0(S)$
the supports of $\pd N_a$ converge to \eq{up4eq11} as $a\ra 0_+$,
so $\pd N$ is supported in~\eq{up4eq11}.
\end{proof}

Note that we do not claim that $N$ is an {\it integral current},
nor that $\pd N$ is {\it rectifiable}. This is because our
assumptions are not strong enough to ensure that $\pd N$ has
finite area (Hausdorff 2-measure) near points $\bigl(0,0,x+
iu(x,0)\bigr)$ for $(x,0)\in\pd S$ with $v(x,0)=0$. For
$u,v$ in Theorem \ref{up4thm4} the estimates of \cite{Joyc5}
do not imply $\pd N$ has finite area, which is why we have
not supposed it. For $u,v$ in Theorem \ref{up4thm8} $\pd N$
does have finite area, and $N$ is an integral current.

\section{$\U(1)$-invariant special Lagrangian cones}
\label{up5}

Combining Proposition \ref{up4prop1} and Theorems \ref{up4thm4} and
\ref{up4thm8} we get powerful existence results for {\it singular}
SL 3-folds of the form \eq{up4eq1} with $a=0$. By Proposition
\ref{up4prop6} we may regard these as {\it minimal rectifiable
currents}, as in \S\ref{up33}, and so by Theorem \ref{up3thm3}
there exists a {\it tangent cone} at each singular point, which
will be a $\U(1)$-invariant SL cone in~$\C^3$.

In this section we will study the possible tangent cones of singular
SL 3-folds of the form \eq{up4eq1} with $a=0$, and find there are
only a few possibilities, which can be written down very explicitly.
We begin by quoting work of the author \cite{Joyc1} and Haskins
\cite{Hask} on $\U(1)$-invariant SL cones in $\C^3$. Our first
result comes from \cite[Th.~8.4]{Joyc1} with $a_1=a_2=-1$ and
$a_3=2$, with some changes in notation.

\begin{thm} Let\/ $N_0$ be a closed special Lagrangian cone in
$\C^3$ invariant under the $\U(1)$-action \eq{up1eq1}, with\/
$N_0\sm\{0\}$ connected. Then there exist\/ $A\in[-1,1]$ and
functions $w:\R\ra(-\ha,1)$ and\/ $\al,\be:\R\ra\R$ satisfying
\e
\begin{gathered}
\Bigl(\frac{\d w}{\d t}\Bigr)^2=4\bigl((1-w)^2(1+2w)-A^2\bigr),\quad
\frac{\d\al}{\d t}=\frac{A}{1-w},\\
\frac{\d\be}{\d t}=\frac{-2A}{1+2w}\quad\text{and\/}\quad
(1-w)(1+2w)^{1/2}\cos(2\al+\be)\equiv A,
\end{gathered}
\label{up5eq1}
\e
such that away from points $(z_1,z_2,z_3)\in\C^3$ with\/ $z_j=0$
for some $j$, we may locally write $N_0$ in the form
$\bigl\{\Phi(r,s,t):r>0$, $s,t\in\R\bigr\}$, where
\e
\begin{split}
\Phi:(r,s,t)\mapsto\bigl(r{\rm e}^{i(\al(t)+s)}\sqrt{1-w(t)},\,
&r{\rm e}^{i(\al(t)-s)}\sqrt{1-w(t)},\\
&r{\rm e}^{i\be(t)}\sqrt{1+2w(t)}\,\,\bigr).
\end{split}
\label{up5eq2}
\e
\label{up5thm1}
\end{thm}

We can say more about $N_0$ by dividing into cases, depending on~$A$.

\begin{thm} In the situation of Theorem \ref{up5thm1}, we have
\begin{itemize}
\setlength{\itemsep}{0pt}
\setlength{\parsep}{0pt}
\item[{\rm(a)}] If\/ $A=1$ then $N_0$ is the $\U(1)^2$-invariant
SL $T^2$-cone
\e
\bigl\{(r{\rm e}^{i\th_1},r{\rm e}^{i\th_2},r{\rm e}^{i\th_3}):
r\ge 0,\quad \th_1,\th_2,\th_3\in\R,\quad \th_1+\th_2+\th_3=0\bigr\}.
\label{up5eq3}
\e
\item[{\rm(b)}] If\/ $A=-1$ then $N_0$ is the $\U(1)^2$-invariant
SL $T^2$-cone
\e
\bigl\{(r{\rm e}^{i\th_1},r{\rm e}^{i\th_2},r{\rm e}^{i\th_3}):
r\ge 0,\quad \th_1,\th_2,\th_3\in\R,\quad \th_1+\th_2+\th_3=\pi\bigr\}.
\label{up5eq4}
\e
\item[{\rm(c)}] If\/ $A=0$ then for some $\phi\in(-\frac{\pi}{2},
\frac{\pi}{2}]$ either $N_0=\Pi^\phi_+$ or $N_0=\Pi^\phi_-$ or $N_0$
is the singular union $\Pi^\phi_+\cup\Pi^\phi_-$, where $\Pi^\phi_\pm$
are the SL\/ $3$-planes
\e
\begin{split}
\Pi^\phi_+&=\bigl\{(z,i{\rm e}^{-i\phi}\bar z,r{\rm e}^{i\phi}):z\in\C,
\quad r\in\R\bigr\}\\
\text{and\/}\quad
\Pi^\phi_-&=\bigl\{(z,-i{\rm e}^{-i\phi}\bar z,r{\rm e}^{i\phi}):z\in\C,
\quad r\in\R\bigr\}.
\end{split}
\label{up5eq5}
\e
\item[{\rm(d)}] If\/ $0<\md{A}<1$ then the function $w(t)$ of Theorem
\ref{up5thm1} may be written explicitly in terms of the Jacobi elliptic
functions. It is non-constant and periodic in $t$ with period\/ $T$
depending only on $A$, and\/ $2\al+\be$ is also non-constant and
periodic in $t$ with period\/~$T$.

Define $\Phi(A)=\be(T)-\be(0)$. Then $\Phi$ depends only on $A$, and\/
$\Phi:(0,1)\ra\R$ is real analytic and strictly monotone decreasing
with\/ $\Phi(A)\ra-\pi$ as $A\ra 0_+$ and\/ $\Phi(A)\ra-2\pi/\sqrt{3}$
as $A\ra 1_-$, and\/ $\Phi:(-1,0)\ra\R$ is real analytic and strictly
monotone decreasing with\/ $\Phi(A)\ra\pi$ as $A\ra 0_-$ and\/
$\Phi(A)\ra 2\pi/\sqrt{3}$ as~$A\ra -1_+$.
\end{itemize}
\label{up5thm2}
\end{thm}

\begin{proof} The division into cases (a)--(d) corresponds to the
cases considered in \cite[\S 7.2--\S 7.5]{Joyc1}. Case (a) comes
immediately from \cite[\S 7.3]{Joyc1}, and (b) follows from (a)
by replacing $z_3$ by $-z_3$. When $A=0$, equation \eq{up5eq1}
shows that $\al,\be$ are constant, and (c) then follows from
Theorem \ref{up5thm1} with a certain amount of work. The point
about the SL 3-planes $\Pi^\phi_\pm$ is that they are invariant
under the $\U(1)$-action \eq{up1eq1} and intersect in the line
$\bigl\{(0,0,r{\rm e}^{i\phi}):r\in\R\bigr\}$. Thus if
$N_0=\Pi^\phi_+\cup\Pi^\phi_-$ then $N_0\sm\{0\}$ is connected,
which is why these two SL 3-planes can be combined under the
hypotheses of Theorem~\ref{up5thm1}.

For case (d), when $A\in(0,1)$ explicit formulae for $w$ in terms
of the Jacobi elliptic functions $\sn(t,k)$ are given in
\cite[Prop.~8.6]{Joyc1} and \cite[Prop.~4.2]{Hask}, and
\cite[Prop.~7.11]{Joyc1} shows that $w$ and $2\al+\be$ are
non-constant and periodic with period $T$. The limiting values
of $\Phi$ as $A\ra 0_+$ and $A\ra 1_-$ follow from
\cite[Prop.~7.13]{Joyc1}, noting that $\Phi$ above is
$\frac{1}{3}\Psi$ in the notation of~\cite{Joyc1}.

By Haskins \cite[p.~20]{Hask}, $\Phi(A)$ is strictly monotone
on $(0,1)$, noting that the $\U(1)$-action \eq{up1eq1} corresponds
to the case $\al=0$ in \cite{Hask}, and $\Phi$ and $A$ in our
notation correspond to $\Th_2$ and $3\sqrt{3}J$ in his notation.
This proves (d) when $A\in(0,1)$. The claims for $A\in(-1,0)$
follow by replacing $z_3$ by~$-z_3$.
\end{proof}

We are interested in the tangent cones not of arbitrary
$\U(1)$-invariant SL 3-folds, but only those which can be written
{\it globally} in the form \eq{up4eq1}. Therefore, in the next
three propositions, we work out which of the SL cones above can
be written globally in the form \eq{up4eq1}. The first follows
from Theorem \ref{up4thm1}, as \eq{up5eq3} agrees with the SL
3-fold $N_0$ of~\eq{up4eq4}.

\begin{prop} The SL\/ $T^2$-cone of\/ \eq{up5eq3} may be written
in the form \eq{up4eq1} with\/ $u=u_0$ and\/ $v=v_0$, where $u_0,v_0$
are as in Theorem \ref{up4thm1}. Similarly, the SL\/ $T^2$-cone of\/
\eq{up5eq4} may be written in the form \eq{up4eq1} with\/ $u=-u_0$
and\/~$v=-v_0$.
\label{up5prop1}
\end{prop}

The second is elementary.

\begin{prop} Let\/ $\phi\in(-\frac{\pi}{2},\frac{\pi}{2})$. Then the
SL\/ $3$-planes $\Pi^\phi_+,\Pi^\phi_-$ of\/ \eq{up5eq5} may be written
\begin{align}
\begin{split}
\Pi^\phi_+=\bigl\{(z_1,z_2,z_3)\in\C^3:\,&
z_1z_2=v+iy,\;\> z_3=x+iu,\;\> \ms{z_1}-\ms{z_2}=0,\\
&u=x\,\tan\phi,\;\> v=y\,\tan\phi, \;\> x\in\R,\;\> y\ge 0\bigr\}.
\end{split}
\label{up5eq6}\\
\begin{split}
\Pi^\phi_-=\bigl\{(z_1,z_2,z_3)\in\C^3:\,&
z_1z_2=v+iy,\;\> z_3=x+iu,\;\> \ms{z_1}-\ms{z_2}=0,\\
&u=x\,\tan\phi,\;\> v=y\,\tan\phi, \;\> x\in\R,\;\> y\le 0\bigr\}.
\end{split}
\label{up5eq7}
\end{align}
Thus the union $\Pi^\phi_+\cup\Pi^\phi_-$ may be written in the form
\eq{up4eq1} with\/ $u=x\,\tan\phi$, $v=y\,\tan\phi$ and\/~$a=0$.
\label{up5prop2}
\end{prop}

However, the cones $N_0$ in case (d) above cannot be written this way.

\begin{prop} None of the SL cones $N_0$ of Theorem \ref{up5thm1}
with\/ $0<\md{A}<1$ may be written globally in the form \eq{up4eq1}
with single-valued functions~$u,v$.
\label{up5prop3}
\end{prop}

\begin{proof} For the cone $N_0$ of Theorem \ref{up5thm1} to be
closed without boundary, we need the function $\Phi$ of \eq{up5eq2}
to be periodic in $t$. Since $w$ and $2\al+\be$ are periodic with
period $T$ by part (d) of Theorem \ref{up5thm2}, the possible periods
of $\Phi$ in $t$ are $qT$ for integers $q\ge 1$, and this will happen
if ${\rm e}^{i\be(t)}$ has period $qT$. Now as $w(t)$ has period $T$,
we see from \eq{up5eq1} and the definition of $\Phi(A)$ that
$\be(t+T)=\be(t)+\Phi(A)$ for all $t$, and hence~$\be(t+qT)=
\be(t)+q\Phi(A)$.

Therefore ${\rm e}^{i\be(t+qT)}={\rm e}^{i\be(t)}$ if and only
if $q\Phi(A)=2\pi p$ for some $p\in\Z$. Hence $\Phi$ has period
$qT$ if $\Phi(A)=2\pi\frac{p}{q}$ for integers $p,q$, where
$q\ge 1$ is as small as possible. But $\Phi(A)$ lies in
$(\pi,2\pi/\sqrt{3})$ or $(-2\pi/\sqrt{3},-\pi)$ by part (d)
of Theorem \ref{up5thm2}, and thus $\frac{p}{q}$ lies in
$(\ha,\frac{1}{\sqrt{3}})$ or $(-\frac{1}{\sqrt{3}},-\ha)$.
It easily follows that $q\ge 3$ and~$\md{p}\ge 2$.

Now by \eq{up5eq2} we have $x+iu=z_3=r\sqrt{1+2w}\,{\rm e}^{i\be}$
and $v+iy=z_1z_2=r^2(1-w){\rm e}^{i(2\al+\be)}{\rm e}^{-i\be}$.
Since $2\al+\be$ is periodic, considering the phases of $x+iu$
and $v+iy$ we see that in one period $qT$ of $t$ the phase of
$x+iu$ rotates through an angle $2\pi p$, and the phase of
$v+iy$ rotates through an angle~$-2\pi p$.

It is not difficult to use this to show that for generic
$(x,y)\in\R^2$ we expect $\md{p}$ points $(u,v)$ in $\R^2$
(possibly counting with multiplicity) for which $x+iu$ and
$v+iy$ can be written in the form above for some $(z_1,z_2,z_3)
\in N_0$. For instance, if $x=0$, $y>0$ and $\cos(2\al(t)+\be(t))
>0$ for all $t$, which is reasonable as $2\al+\be$ is periodic,
then there is one point $(u,v)$ over $(x,y)$ for each $t\in[0,T)$
with ${\rm e}^{i\be(t)}=i$. Since $\be$ increases by $2\pi p$ on
$[0,T]$, there will be $\md{p}$ such points.

Thus, as $\md{p}\ge 2$ there will be at least 2 points $(u,v)$
over each generic $(x,y)$, so $N_0$ cannot be written in the
form \eq{up4eq1} for single-valued functions $u,v$, but only
for multi-valued `functions' $(u,v)$ with at least 2 values
at generic points.
\end{proof}

We can now classify the possible tangent cones in our problem.

\begin{thm} Let\/ $S$ be a domain in $\R^2$ and\/ $u,v\in C^0(S)$
a singular solution of\/ \eq{up4eq2}, as in \S\ref{up46}. Define
$N$ by \eq{up4eq1} with\/ $a=0$. Then by Proposition \ref{up4prop6}
we may regard $N$ as a minimal rectifiable current, as in~\S\ref{up33}.

Let\/ ${\bf z}\in N^\circ$ be a singular point, so that\/ ${\bf z}=
(0,0,z_3)$ by Proposition \ref{up4prop1}(a), and\/ $C$ be a tangent
cone to $N$ at\/ $\bf z$. Then the only possibilities for $C$ are
\begin{itemize}
\setlength{\itemsep}{0pt}
\setlength{\parsep}{0pt}
\item[{\rm(i)}] $C$ is given in \eq{up5eq3}, with multiplicity~$1$.
\item[{\rm(ii)}] $C$ is given in \eq{up5eq4}, with multiplicity~$1$.
\item[{\rm(iii)}] $C$ is $\Pi^\phi_+$ or $\Pi^\phi_-$ or
$\Pi^\phi_+\cup\Pi^\phi_-$ for some $\phi\in(-\frac{\pi}{2},
\frac{\pi}{2})$, where $\Pi^\phi_\pm$ are defined in
\eq{up5eq5} and have multiplicity~$1$.
\item[{\rm(iv)}] $C$ is the sum of\/ $\Pi^{\pi/2}_+$ with
multiplicity $k$ and\/ $\Pi^{\pi/2}_-$ with multiplicity $l$,
for nonnegative integers $k,l$ with\/~$k+l\ge 1$.
\end{itemize}
\label{up5thm3}
\end{thm}

\begin{proof} Let $C$ be a tangent cone to $N$ at $\bf z$. Then
$C$ is a minimal locally rectifiable cone without boundary, and is
special Lagrangian by Theorem \ref{up3thm3}. Since $N$ is invariant
under the $\U(1)$-action \eq{up1eq1} which fixes $\bf z$ it easily
follows that $C$ is invariant under \eq{up1eq1}. Also, $C$ is an
embedded minimal 3-submanifold with positive integer multiplicity
outside a closed singular set $S$ of Hausdorff dimension at most 1,
by Theorem~\ref{up3thm1}.

Since $N\sm\{{\bf z}\}$ is locally connected, by considering the
limiting process defining $C$ one can show that $C\sm\{0\}$ is also
connected. As $C$ is a $\U(1)$-invariant cone, $S$ is also a
$\U(1)$-invariant cone. Therefore $S$ is a union of $(0,0,0)$ and
some collection of rays $\bigl\{(0,0,r{\rm e}^{i\phi}):r>0\bigr\}$,
as otherwise $S$ would contain a Hausdorff dimension 2 cone on a
$\U(1)$-orbit. So~$S\subset\bigl\{(0,0,z):z\in\C\bigr\}$.

Thus $C$ is a $\U(1)$-invariant SL cone in $\C^3$, with multiplicity,
outside a singular set $S$. Such cones are locally classified in
Theorem \ref{up5thm1}, and so $C\sm S$ is locally parametrized by
$\Phi$ as in \eq{up5eq2}. Hence, each connected component of $C\sm S$
fits into the framework of Theorems \ref{up5thm1} and \ref{up5thm2},
with some value of $A$, functions $u,\al$ and $\be$, and so on,
and some positive integer multiplicity.

In cases (a)--(d) of Theorem \ref{up5thm2}, only the cones $N_0$ of
case (c) intersect $\bigl\{(0,0,z):z\in\C\bigr\}$ other than at 0.
Since $S$ lies in $\bigl\{(0,0,z):z\in\C\bigr\}$, and $C$ is closed
and nonsingular except at $S$, we see that if a component of $C\sm\{0\}$
locally agrees with a cone $N_0$ in cases (a), (b), (d) of Theorem
\ref{up5thm2} then $C$ must contain all of $N_0$, since otherwise
the boundary of the included subset of $N_0$ would lie in~$S$.

By an elementary calculation we can prove:

\begin{lem} Let\/ $N$ be an SL\/ $3$-fold of the form \eq{up4eq1}
with\/ $u,v:S\!\ra\!\R$ locally $C^1$ almost everywhere, and let\/
$\chi:\R^2\!\ra\!\R$ be smooth and compactly-supported. Then
\e
\int_N\chi\bigl(\Re z_3,\Im(z_1z_2)\bigr)\cdot
\Re(\d z_1\w\d z_2)\w\Re(\d z_3)=
2\pi\int_S\chi(x,y)\,\d x\,\d y.
\label{up5eq8}
\e
\label{up5lem}
\end{lem}

The important point here is that the right hand side of
\eq{up5eq8} is {\it independent of\/} $u,v$. Now $C$ is the
limit, as a current, of a sequence $r_j^{-1}(N-{\bf z})$ as
$j\ra\iy$. Each $r_j^{-1}(N-{\bf z})$ is of the form
\eq{up4eq1} for single-valued $u,v$, which are locally $C^1$
almost everywhere by Definition \ref{up4def3}(iii), so that
\eq{up5eq8} holds for~$r_j^{-1}(N-{\bf z})$.

We wish to take the limit of \eq{up5eq8} for $r_j^{-1}(N-{\bf z})$
as $j\ra\iy$. Currents and their limits are defined by the integral
of smooth, compactly-supported forms, and the 3-form on the l.h.s.\
of \eq{up5eq8} is {\it not\/} compactly-supported. However, if we
take $\chi$ to be {\it nonnegative} then this 3-form has nonnegative
restriction to each $r_j^{-1}(N-{\bf z})$. Taking the limit $j\ra\iy$
then shows that
\e
\int_C\chi\bigl(\Re z_3,\Im(z_1z_2)\bigr)\cdot
\Re(\d z_1\w\d z_2)\w\Re(\d z_3)\le
2\pi\int_{\R^2}\chi(x,y)\,\d x\,\d y,
\label{up5eq9}
\e
since the effect of a portion of $r_j^{-1}(N-{\bf z})$ going to
infinity in the support of $\chi\bigl(\Re z_3,\Im(z_1z_2)\bigr)$
is that a positive contribution to the l.h.s.\ of \eq{up5eq8} for
$r_j^{-1}(N-{\bf z})$ does not appear in the l.h.s.\ of~\eq{up5eq9}.

Over a small open neighbourhood of a generic point $(x,y)\in\R^2$,
$C$ is nonsingular, and splits into $k$ components, each of the
form \eq{up4eq1} for single-valued $u,v$. But then Lemma \ref{up5lem}
for $\chi$ supported near $x,y$ shows that
\begin{equation*}
\int_C\chi\bigl(\Re z_3,\Im(z_1z_2)\bigr)\cdot
\Re(\d z_1\w\d z_2)\w\Re(\d z_3)=2\pi k
\int_{\R^2}\chi(x,y)\,\d x\,\d y.
\end{equation*}
Comparing this with \eq{up5eq9} shows that the only possibilities
are $k=0,1$. Therefore for {\it generic} $(x,y)\in\R^2$ there are
either 0 or 1 points $(u,v)$, counted with multiplicity, such that
there exists $(z_1,z_2,z_3)\in C$ with $z_1z_2=v+iy$ and $z_3=x+iu$.
So by Proposition \ref{up5prop3}, the cones $N_0$ of part (d) of
Theorem \ref{up5thm2} do not occur, even locally, as tangent cones
$C$ to~$N$.

Hence, $C\sm\bigl\{(0,0,z):z\in\C\bigr\}$ is a union of connected
components, with multiplicity, each of the form $N_0\sm\bigl\{(0,0,z):
z\in\C\bigr\}$, where $N_0$ is one of the SL cones in parts (a)--(c) of
Theorem \ref{up5thm2}. If there is more than one component, their closures
must intersect in $S\sm\{0\}$ to ensure that $C\sm\{0\}$ is connected.
The only way for this to happen is if $C=\Pi^\phi_+\cup\Pi^\phi_-$ as
in part (c) of Theorem \ref{up5thm2}, with $S=\Pi^\phi_+\cap\Pi^\phi_-$,
so that $C\sm S$ has two connected components $\Pi^\phi_+\sm S$ and
$\Pi^\phi_-\sm S$. Thus, $C$ cannot combine more than one possibility
from parts (a)--(c) of Theorem~\ref{up5thm2}.

It remains only to pin down the multiplicities of each component of
$C\sm S$. First note that they are all positive, as SL 3-folds
cannot converge in the sense of currents to SL 3-folds with the
opposite orientation. In cases (i)--(iii), Propositions \ref{up5prop1}
and \ref{up5prop2} show that $C$ can be written in the form \eq{up4eq1}
with $S$ a half-plane $y\ge 0$ or $y\le 0$, or $S=\R^2$. But as above,
for generic $(x,y)\in\R^2$ there corresponds no more than one point
$(u,v)$, counted with multiplicity. Thus the multiplicities in
(i)--(iii) are all 1. In part (iv) all we know is that the
multiplicities are nonnegative integers, and not both zero.
\end{proof}

We finish by calculating the densities of the cones at their
vertices, as in~\S\ref{up33}.

\begin{prop} In cases {\rm(i)} and\/ {\rm(ii)} of Theorem \ref{up5thm3}
we have $\Th(C,0)=\pi/\sqrt{3}\approx 1\cdot 81$. In case {\rm(iii)}
$\Th(C,0)$ is $1$ when $C=\Pi^\phi_+$ or $\Pi^\phi_-$ and\/ $2$
when $C=\Pi^\phi_+\cup\Pi^\phi_-$, and in case {\rm(iv)}
$\Th(C,0)=k+l$.
\label{up5prop4}
\end{prop}

\begin{proof} For case (i), let $\Si=N_0\cap{\cal S}^5$, where
${\cal S}^5$ is the unit sphere in $\C^3$. The metric on $\Si\cong T^2$
is isometric to the quotient of $\R^2$ with its flat Euclidean metric
by the lattice $\Z^2$ with basis $2\pi\bigl(\frac{\sqrt{2}}{\sqrt{3}},
0\bigl)$, $2\pi\bigl(\frac{1}{\sqrt{6}},\frac{1}{\sqrt{2}}\bigr)$.
Taking a $2\t 2$ determinant gives $\area(\Si)=4\pi^2/\sqrt{3}$.
The density follows by dividing by $\area({\cal S}^2)=4\pi$. Cases
(i) and (ii) are isomorphic under $z_3\mapsto -z_3$, and so have
the same density. Cases (iii) and (iv) are immediate, as
$\Pi^\phi_\pm$ are $\R^3$ vector subspaces.
\end{proof}

\section{Multiplicities of zeroes and counting formulae}
\label{up6}

We shall now generalize the material of \S\ref{up43} to the
singular case $a=0$. We begin by defining the {\it multiplicity}
of an isolated zero $(u,v)-(\hat u,\hat v)$, where $u,v$ and
$\hat u,\hat v$ are singular solutions of~\eq{up4eq2}.

\begin{dfn} Suppose $S$ is a domain in $\R^2$, and $(u,v)$,
$(\hat u,\hat v)$ are singular solutions of \eq{up4eq2} in $S$,
as in \S\ref{up46}. We call a point $(b,c)\in S$ a {\it zero} of
$(u,v)-(\hat u,\hat v)$ in $S$ if $(u,v)=(\hat u,\hat v)$ at $(b,c)$.
A zero $(b,c)$ is called {\it singular} if $c=0$ and $v(b,0)=
\hat v(b,0)=0$, so that $(b,c)$ is a {\it singular point} of $u,v$
and $\hat u,\hat v$. Otherwise we say $(b,c)$ is a {\it nonsingular zero}.
We call a zero $(b,c)$ {\it isolated} if for some $\ep>0$ there
exist no other zeroes $(x,y)$ of $(u,v)-(\hat u,\hat v)$ in $S$
with~$0<(x-b)^2+(y-c)^2\le\ep^2$.

Let $(b,c)\in S^\circ$ be an isolated zero of $(u,v)-(\hat u,\hat v)$.
Define the {\it multiplicity} of $(b,c)$ to be the winding number of
$(u,v)-(\hat u,\hat v)$ about 0 along the positively oriented circle
$\ga_\ep(b,c)$ of radius $\ep$ about $(b,c)$, where $\ep>0$ is chosen
small enough that $\ga_\ep(b,c)$ lies in $S^\circ$ and $(b,c)$ is the
only zero of $(u,v)-(\hat u,\hat v)$ inside or on~$\ga_\ep(b,c)$.

As $(u,v)-(\hat u,\hat v)$ is {\it continuous}, though not necessarily
differentiable, its winding number about 0 along $\ga_\ep(b,c)$ is
defined as in Definition \ref{up4def1}. Since winding numbers are
invariant under continuous deformation of the path, this is
independent of $\ep$, so the multiplicity of $(b,c)$ is well-defined.
\label{up6def1}
\end{dfn}

We shall eventually prove that if $(u,v)\not\equiv(\hat u,\hat v)$
then zeroes $(b,c)$ of $(u,v)-(\hat u,\hat v)$ in $S^\circ$ are {\it
isolated}, with {\it positive} multiplicity. It is simple to show
that multiplicities are {\it nonnegative}, and we do this in
\S\ref{up61}. But if $(b,0)$ is a {\it singular} zero of $(u,v)-
(\hat u,\hat v)$, it is more difficult to prove $(b,0)$ is
{\it isolated}, or has {\it positive} multiplicity.

To do this requires a diversion in \S\ref{up62}, to study singular
solutions $u,v$ with $v(x,0)\equiv 0$. We show such $u,v$ have the
symmetry $u(x,-y)=u(x,y)$, $v(x,-y)=-v(x,y)$, and the corresponding
SL 3-folds $N$ are actually the union of two {\it nonsingular}
$\U(1)$-invariant SL 3-folds $N_\pm$ intersecting in a real
analytic real curve~$\ga$.

In \S\ref{up63} we use this to show that if $v<\hat v$ on
$\pd S$, then $v<\hat v$ on $S$. Section \ref{up64} then proves
that multiplicities of zeroes in Definition \ref{up6def1} are
{\it positive}. This is the essential step in showing zeroes
are {\it isolated}, which we do in \S\ref{up65}--\S\ref{up66}.
Sections \ref{up65} and \ref{up67} also prove {\it counting
formulae} for zeroes in $S^\circ$ with multiplicity, in terms
of boundary data.

\subsection{Nonnegativity of winding numbers and multiplicities}
\label{up61}

We begin with four propositions. The first implies that {\it
intersection multiplicities of singular solutions are nonnegative}.
Later in Corollary \ref{up6cor2} we will show that they are in fact
{\it positive}. Much of the intervening discussion is to exclude
the possibility of intersections with multiplicity zero.

\begin{prop} Let\/ $u,v\in C^0(S)$ and\/ $\hat u,\hat v\in C^0(S)$ be
singular solutions of\/ \eq{up4eq2} in a domain $S$ in $\R^2$, as
in \S\ref{up46}, such that\/ $(u,v)\ne(\hat u,\hat v)$ at every point
of\/ $\pd S$. Then the winding number of\/ $(u,v)-(\hat u,\hat v)$
about\/ $0$ along $\pd S$ is a nonnegative integer.
\label{up6prop1}
\end{prop}

\begin{proof} By part (iv) of Definition \ref{up4def3}, we may
write $u,v$ and $\hat u,\hat v$ as the limits in $C^0(S)$ as
$a\ra 0_+$ of solutions $u_a,v_a$ and $\hat u_a,\hat v_a$ in
$C^2(S)$ of \eq{up4eq3} for $a\in(0,1]$. As $(u,v)\ne(\hat u,\hat v)$
on $\pd S$ we see that for small $a\in(0,1]$ we have $(u_a,v_a)\ne
(\hat u_a,\hat v_a)$ on $\pd S$, and the winding numbers of
$(u,v)-(\hat u,\hat v)$ and $(u_a,v_a)-(\hat u_a,\hat v_a)$
about 0 along $\pd S$ are equal, as winding numbers are invariant
under continuous deformation of the path. But by Theorem
\ref{up4thm2} the winding number of $(u_a,v_a)-(\hat u_a,\hat v_a)$
about 0 along $\pd S$ is a nonnegative integer, so the result follows.
\end{proof}

\begin{prop} Let\/ $u,v\in C^0(S)$ and\/ $\hat u,\hat v\in C^0(S)$ be
singular solutions of\/ \eq{up4eq2} in a domain $S$ in $\R^2$, as
in \S\ref{up46}, such that\/ $(u,v)\ne(\hat u,\hat v)$ at every point
of\/ $\pd S$. Suppose $(b_1,c_1),\ldots,(b_n,c_n)$ are isolated
zeroes of\/ $(u,v)-(\hat u,\hat v)$ in $S^\circ$ with multiplicities
$k_1,\ldots,k_n$, but not necessarily the only zeroes. Then the
winding number of\/ $(u,v)-(\hat u,\hat v)$ about\/ $0$ along
$\pd S$ is at least\/~$\sum_{i=1}^nk_i$.
\label{up6prop2}
\end{prop}

\begin{proof} For each $i=1,\ldots,n$, choose $\ep_i>0$ such that
$\ga_{\ep_i}(b_i,c_i)$ lies in $S^\circ$ and $(b_i,c_i)$ is the only
zero of $(u,v)-(\hat u,\hat v)$ inside or on $\ga_{\ep_i}(b_i,c_i)$,
and in addition the circles $\ga_{\ep_i}(b_i,c_i)$ do not intersect.
Clearly this is possible if $\ep_1,\ldots,\ep_n$ are small enough.
Examining the proof of Theorem \ref{up4thm2} in \cite{Joyc4} we see
that it holds not only for domains $S$, which are contractible, but
also for more general compact 2-submanifolds $T$ of $\R^2$ with
finitely many boundary components. Thus Proposition \ref{up6prop1}
also holds for such~$T$.

Let $T=S\sm\bigcup_{i=1}^nB_{\ep_i}(b_i,c_i)$. Then $T$ is a compact
2-submanifold in $\R^2$ whose boundary is the disjoint union of
$\pd S$, and the circles $\ga_{\ep_i}(b_i,c_i)$ for $i=1,\ldots,n$,
with negative orientation. Let the winding number of $(u,v)-
(\hat u,\hat v)$ about 0 along $\pd S$ be $k$. Then the winding number
of $(u,v)-(\hat u,\hat v)$ about 0 along $\pd T$ is $k-\sum_{i=1}^nk_i$.
So $k\ge\sum_{i=1}^nk_i$ by Proposition~\ref{up6prop1}.
\end{proof}

\begin{prop} Let\/ $u,v\in C^0(S)$ and\/ $\hat u,\hat v\in C^0(S)$ be
singular solutions of\/ \eq{up4eq2} in a domain $S$ in $\R^2$, such
that\/ $(u,v)\not\equiv(\hat u,\hat v)$. Suppose $(b,c)\in S^\circ$
is a nonsingular zero of\/ $(u,v)-(\hat u,\hat v)$. Then $(b,c)$ is an
isolated zero of\/ $(u,v)-(\hat u,\hat v)$, and its multiplicity is a
positive integer $k$, with\/ $\pd^ju(b,c)=\pd^j\hat u(b,c)$ and\/
$\pd^jv(b,c)=\pd^j\hat v(b,c)$ for $j=0,\ldots,k-1$, but\/
$\pd^ku(b,c)\ne\pd^k\hat u(b,c)$ and\/~$\pd^kv(b,c)\ne\pd^k\hat v(b,c)$.
\label{up6prop3}
\end{prop}

\begin{proof} As $u,v$ and $\hat u,\hat v$ are nonsingular near $(b,c)$
we can apply the reasoning of \cite[\S 6.1]{Joyc4} for the nonsingular
case $a\ne 0$ to $(u,v)$ and $(\hat u,\hat v)$ near $(b,c)$. By
\cite[Lem.~6.4]{Joyc4} we see that if $(u,v)\not\equiv(\hat u,\hat v)$
then $(b,c)$ is an {\it isolated\/} zero of $(u,v)-(\hat u,\hat v)$, and
it has a unique multiplicity $k$, defined as in Definition \ref{up4def2},
which is a positive integer.

If $\ep>0$ is small enough then $(u,v)$ and $(\hat u,\hat v)$ are
nonsingular on the closed disc $\ovB_\ep(b,c)$, and $(b,c)$ is the
only zero of $(u,v)-(\hat u,\hat v)$ there. The proof of Theorem
\ref{up4thm2} in \cite{Joyc4} when $a\ne 0$ is also valid in this
case on $\ovB_\ep(b,c)$, and shows that the winding number of
$(u,v)-(\hat u,\hat v)$ about 0 along $\ga_\ep(b,c)$ is $k$. Thus
this multiplicity $k$ coincides with that in Definition
\ref{up6def1}, and the proof is complete.
\end{proof}

When $\hat u=\al x+\be$, $\hat v=\al y+\ga$, intersection
multiplicities are positive.

\begin{prop} Let\/ $(u,v)$ be a singular solution of\/ \eq{up4eq2}
in a domain $S$ in $\R^2$, as in \S\ref{up46}, and let\/ $\hat u(x,y)
=\al x+\be$ and $\hat v(x,y)=\al y+\ga$ for $\al,\be,\ga\in\R$, as
in Example \ref{up4ex1}. Suppose $(u,v)-(\hat u,\hat v)$ has a zero
$(b,c)$ in $S^\circ$. Then the winding number of\/ $(u,v)-(\hat u,
\hat v)$ about\/ $0$ along $\pd S$ is a positive integer.
\label{up6prop4}
\end{prop}

\begin{proof} By part (iv) of Definition \ref{up4def3}, we may
write $u,v$ as the limits in $C^0(S)$ as $a\ra 0_+$ of solutions
$u_a,v_a$ in $C^2(S)$ of \eq{up4eq3} for $a\in(0,1]$. Define
$\hat u_a=\al(x-b)+u_a(b,c)$ and $\hat v_a=\al(y-c)+v_a(b,c)$
for $a\in(0,1]$. Then $(\hat u_a,\hat v_a)$ satisfies \eq{up4eq3}
and $\hat u_a,\hat v_a\ra\hat u,\hat v$ in $C^0(S)$ as~$a\ra 0_+$.

Since $(u,v)\ne(\hat u,\hat v)$ on $\pd S$ we see that for small
$a\in(0,1]$ we have $(u_a,v_a)\ne(\hat u_a,\hat v_a)$ on $\pd S$,
and the winding numbers of $(u,v)-(\hat u,\hat v)$ and $(u_a,v_a)
-(\hat u_a,\hat v_a)$ about 0 along $\pd S$ are equal. But as
$(u_a,v_a)=(\hat u_a,\hat v_a)$ at $(b,c)\in S^\circ$, by Theorem
\ref{up4thm2} the winding number of $(u_a,v_a)-(\hat u_a,\hat v_a)$
about 0 along $\pd S$ is a positive integer, so the result follows.
\end{proof}

\subsection{Solutions $u,v$ of \eq{up4eq2} with $v(x,0)\equiv 0$}
\label{up62}

We now study singular solutions $u,v\in C^0(S)$ of \eq{up4eq2}
with $v(x,0)\equiv 0$. Consider the following situation.

\begin{dfn} Let $S$ be a domain in $\R^2$ which intersects the
$x$-axis in $[x_1,x_2]\t\{0\}$ for $x_1<x_2$, with $(x_1,x_2)\t
\{0\}\subset S^\circ$, and $u,v\in C^0(S)$ a singular solution
of \eq{up4eq2} in $S$, with $v(x,0)=0$ for all $x\in[x_1,x_2]$.
That is, $(u,v)$ is singular all along the intersection of $S$
with the $x$-axis. Define subsets $N_\pm$ in $\C^3$ by
\begin{align}
\begin{split}
N_+=\bigl\{(z_1,z_2,z_3)\in\C^3:\,&
z_1z_2=v(x,y)+iy,\quad z_3=x+iu(x,y),\\
&\md{z_1}=\md{z_2},\quad (x,y)\in S,\quad y\ge 0\bigr\},
\end{split}
\label{up6eq1}\\
\begin{split}
N_-=\bigl\{(z_1,z_2,z_3)\in\C^3:\,&
z_1z_2=v(x,y)+iy,\quad z_3=x+iu(x,y),\\
&\md{z_1}=\md{z_2},\quad (x,y)\in S,\quad y\le 0\bigr\}.
\end{split}
\label{up6eq2}
\end{align}
Then $N_+\cap N_-$ is the real curve $\ga=\bigl\{\bigl(0,0,x+iu(x,0)
\bigr):x\in[x_1,x_2]\bigr\}$ in $\C^3$. The {\it end points} of $\ga$
are $\bigl(0,0,x_j+iu(x_j,0)\bigr)$ for $j=1,2$, and the {\it interior}
$\ga^\circ$ of $\ga$ is~$\bigl\{\bigl(0,0,x+iu(x,0)\bigr):x\in(x_1,x_2)
\bigr\}$.
\label{up6def2}
\end{dfn}

We will prove that $N_\pm$ are $\U(1)$-invariant SL 3-folds
with boundaries
\begin{align}
\begin{split}
\pd N_+=\bigl\{(z_1,z_2,z_3)\in\C^3:\,&
z_1z_2=v(x,y)+iy,\quad z_3=x+iu(x,y),\\
&\md{z_1}=\md{z_2},\quad (x,y)\in\pd S,\quad y\ge 0\bigr\},
\end{split}
\label{up6eq3}\\
\begin{split}
\pd N_-=\bigl\{(z_1,z_2,z_3)\in\C^3:\,&
z_1z_2=v(x,y)+iy,\quad z_3=x+iu(x,y),\\
&\md{z_1}=\md{z_2},\quad (x,y)\in\pd S,\quad y\le 0\bigr\},
\end{split}
\label{up6eq4}
\end{align}
which are {\it nonsingular} on their interiors, {\it including}
along~$\ga^\circ$.

Here is a large class of examples $(u,v)$ satisfying Definition
\ref{up6def2}. Note also that Example \ref{up4ex2} has $v(x,0)
\equiv 0$. We shall show in Theorem \ref{up6thm3} that any
$(u,v)$ of this form has the symmetries $u(x,-y)\equiv u(x,y)$
and $v(x,-y)\equiv -v(x,y)$ of these examples.

\begin{ex} Let $S$ be a strictly convex domain in $\R^2$ invariant
under the involution $(x,y)\mapsto(x,-y)$, let $k\ge 0$ and $\al\in
(0,1)$. Suppose $\phi\in C^{k+3,\al}(\pd S)$ with $\phi(x,y)\equiv
-\phi(x,-y)$. Let $f\in C^1(S)$ be the unique weak solution of
\eq{up4eq8} with $a=0$ with $f\vert_{\pd S}=\phi$, which exists
by Theorem \ref{up4thm4}. Then $f'(x,y)=-f(x,-y)$ also satisfies
\eq{up4eq8} with $a=0$, and $f'\vert_{\pd S}=\phi$ as $\phi(x,y)
\equiv-\phi(x,-y)$. Hence $f'=f$ by uniqueness, so
that~$f(x,-y)\equiv -f(x,y)$.

Let $u=\frac{\pd f}{\pd y}$ and $v=\frac{\pd f}{\pd x}$. Then
$u,v\in C^0(S)$ are {\it singular solutions} of \eq{up4eq2}
in the sense of \S\ref{up46}, by Theorem \ref{up4thm4}. Moreover
$f(x,-y)\equiv -f(x,y)$ implies that $u(x,-y)\equiv u(x,y)$ and
$v(x,-y)\equiv -v(x,y)$. This gives $v(x,0)=0$ for all $(x,0)\in S$,
so $(u,v)$ is singular all along the intersection of $S$ with the
$x$-axis.
\label{up6ex}
\end{ex}

We shall show the function $u(x,0)$ used to define $\ga$ is {\it
Lipschitz} on closed subintervals $[b,c]$ of the $x$-axis in~$S^\circ$.

\begin{prop} Let\/ $S,u,v,x_1,x_2$ be as in Definition \ref{up6def2},
and suppose $x_1<b<c<x_2$. Then there exists $K>0$ such that if\/
$b\le x'< x''\le c$ then~$\bmd{u(x'',0)-u(x',0)}\le K\md{x''-x'}$.
\label{up6prop5}
\end{prop}

\begin{proof} As $u,v\in C^0(S)$ we have $(u^2+v^2)^{1/2}\le A$
on $S$ for some $A>0$. Choose $K>0$ large enough that $K\bigl[
(x-x')^2+y^2\bigr]^{1/2}\ge 2A$ on $\pd S$ for all $x'\in[b,c]$,
and $\bmd{u(c,0)-u(b,0)}\le K\md{c-b}$. This is clearly
possible.

Suppose for a contradiction that $b\le x'< x''\le c$
with $\bmd{u(x'',0)-u(x',0)}>K\md{x''-x'}$. Define $\al=(u(x'',0)
-u(x',0))/(x''-x')$, so that $\md{\al}>K$, and $\be=u(x',0)-\al x'$.
Set $\hat u(x,y)=\al(x-x')+\be$ and $\hat v(x,y)=\al y$. Then
$u(x',0)=\al x'+\be$ and~$u(x'',0)=\al x''+\be$.

The conditions $\bmd{u(c,0)-u(b,0)}\le K\md{c-b}$ and
$\md{\al}>K$ imply that $u(b,0)=\al b+\be$ and $u(c,0)=
\al c+\be$ cannot both hold. Thus, $u(x,0)=\al x+\be$ holds
for a proper closed subset of $[b,c]$ containing $x',x''$.
If this subset contains $[x',x'']$, then decrease $x'$ or
increase $x''$ until it no longer lies in this subset, but
$\bmd{u(x'',0)-u(x',0)}>K\md{x''-x'}$ still holds. With these
new $x',x''$ it is not true that $u(x,0)=\al x+\be$ for
all~$x\in[x',x'']$.

Then $\hat u,\hat v$ satisfy \eq{up4eq2}, and $(u,v)=(\hat u,\hat v)$
at $(x',0)$ and $(x'',0)$ by construction. Now on $\pd S$ we have
\begin{align*}
(\hat u^2+\hat v^2)^{1/2}&\ge\md{\al}\bigl((x-x')^2+y^2\bigr)^{1/2}
-\md{u(x',0)}\\
&>K[(x-x')^2+y^2]^{1/2}-A\ge 2A-A=A\ge (u^2+v^2)^{1/2}.
\end{align*}
Hence the winding number of $(u,v)-(\hat u,\hat v)$ about $\pd S$
is the same as that of $-(\hat u,\hat v)$ about $\pd S$, which is 1,
as $(\hat u,\hat v)$ has a single zero on the $x$-axis in~$S^\circ$.

Write $E=\bigl\{(x,y)\in S:(u,v)=(\hat u,\hat v)$ at $(x,y)\bigr\}$.
Then $E$ is compact and contained in $S^\circ$, with $(x',0),(x'',0)
\in E$. Proposition \ref{up6prop3} implies that if $(b,c)\in E$ is
nonsingular then it is isolated, and as $(u,v)-(\hat u,\hat v)$
has winding number 1 about $\pd S$, Proposition \ref{up6prop2}
shows there is at most 1 nonsingular zero $(b,c)$ in $E$. Hence
$E$ is the union of a closed subset of the $x$-axis, and at most
one other point.

Thus, the connected components of $E$ are closed intervals in the
$x$-axis, plus at most one other point. Proposition \ref{up6prop4}
implies that if $T\subset S$ is a subdomain with $\pd T\cap E=
\emptyset$ and at least one connected component of $E$ in $T^\circ$,
then the winding number of $(u,v)-(\hat u,\hat v)$ about 0 along
$\pd T$ is positive. So by the argument of Proposition \ref{up6prop2}
we see that $E$ has at most one connected component. However,
$(x',0)$ and $(x'',0)$ must lie in different connected components
of $E$, as it is not true that $u(x,0)=\al x+\be$ for all
$x\in[x',x'']$, so that $[x',x'']\t\{0\}\not\subseteq E$. So $E$
has at least two connected components, a contradiction.
\end{proof}

Now we can prove the $N_\pm$ are {\it rectifiable currents}, with
boundaries \eq{up6eq3}, \eq{up6eq4}.

\begin{prop} In the situation above, $N_\pm$ are the supports of
special Lagrangian rectifiable currents, which we identify with\/
$N_\pm$, with boundaries $\pd N_\pm$ supported on \eq{up6eq3} and\/
\eq{up6eq4}. In particular, $\ga^\circ$ lies in the interiors
$N_\pm^\circ$, not the boundaries.
\label{up6prop6}
\end{prop}

\begin{proof} Observe that $N_\pm$ are the intersection of $N$
with the subsets $\Im(z_1z_2)\ge 0$, $\Im(z_1z_2)\le 0$ of
$\C^3$, where $N$ is the support of special Lagrangian rectifiable
current with boundary supported on \eq{up4eq11}, by Proposition
\ref{up4prop6}. So by the definition \cite[1.4]{Morg} of
rectifiable $m$-current as a countable union of Lipschitz images
of bounded measurable subsets of $\R^m$, we see that $N_\pm$ are
the supports of {\it special Lagrangian rectifiable currents}, by
making the subsets of $\R^m$ smaller.

It is also easy to see that
\e
\begin{split}
\supp(\pd N_+)\subseteq&\supp(\pd N)\cap
\bigl\{(z_1,z_2,z_3)\in\C^3:\Im(z_1z_2)\ge 0\bigr\}\,\cup\\
&\supp(N)\cap\bigl\{(z_1,z_2,z_3)\in\C^3:\Im(z_1z_2)=0\bigr\},
\end{split}
\label{up6eq5}
\e
and similarly for $\pd N_-$. By Definition \ref{up6def2} and
equations \eq{up4eq1} and \eq{up4eq11} we see that the first
line of the right hand side of \eq{up6eq5} is \eq{up6eq3},
and the second line is~$\ga$.

Now by \cite[Th.~4.7]{Morg}, if $T$ is an $m$-current in $\R^n$
and $\supp(T)$ has zero Hausdorff $m$-measure, then $T=0$. From
Proposition \ref{up6prop5} we see that $\ga^\circ$ has Hausdorff
dimension 1, and zero Hausdorff 2-measure. Hence the portion of
the current $\pd N_+$ supported on $\ga^\circ$ is zero. That is,
$\pd N_+$ is supported on \eq{up6eq3}, with $\ga^\circ$ in the
interior of $N_+$, and similarly for~$N_-$.
\end{proof}

Next we identify the {\it tangent cones} to $N_\pm$ along the
singular curve~$\ga^\circ$.

\begin{prop} In the situation above, at each\/ ${\bf z}\in\ga^\circ$,
any tangent cone $C$ to $N_\pm$ is $\Pi^\phi_\pm$ with multiplicity $1$
for some $\phi\in(-\frac{\pi}{2},\frac{\pi}{2})$, as in part\/
{\rm(iii)} of Theorem~\ref{up5thm3}.
\label{up6prop7}
\end{prop}

\begin{proof} We give the proof for $N_+$. As $N_+$ lies in
the subset $y=\Im(z_1z_2)\ge 0$ in $\C^3$ it follows that
$\Im(z_1z_2)\ge 0$ on $C$. But the cones $N_0$ in parts (i)
and (ii) of Theorem \ref{up5thm3}, and the planes $\Pi^\phi_-$
in part (iii), all contain points with $\Im(z_1z_2)<0$. The
only remaining possibilities for $C$ are $\Pi^\phi_+$ with
multiplicity 1 as in part (iii) of Theorem \ref{up5thm3}, or
part (iv) of Theorem \ref{up5thm3}. So we have to eliminate
possibility~(iv).

As a current, $C$ is the limit as $j\ra\iy$ of a sequence
$r_j^{-1}(N_+-{\bf z})$, where ${\bf z}=(0,0,x'+iu(x',0))$
for some $x'\in(x_1,x_2)$, and $r_j\ra 0_+$. It is easy to
show that $r_j^{-1}(N_+-{\bf z})$ may be written in the form
\eq{up4eq1} with $a=0$ on a neighbourhood of $(0,0)$
in $\{(x,y)\in\R^2:y\ge 0\}$, with $u,v$ replaced by
\begin{equation*}
u_j(x,y)=r_j^{-1}u\bigl(r_j(x-x'),r_j^2y\bigr)
\quad\text{and}\quad
v_j(x,y)=r_j^{-2}v\bigl(r_j(x-x'),r_j^2y\bigr).
\end{equation*}
It then follows from Proposition \ref{up6prop5} that
there exists $K>0$ such that for $\md{x}\le 1$ we have
$\md{u_j(x,0)}\le K\md{x}$ for all sufficiently large $j$.
We see from this that the cone $C$ must contain a ray
$\{(0,0,x+i\al x):x\ge 0\}$ for some $\al\in\R$ with
$\md{\al}\le K$. This holds for $C=\Pi^\phi_+$ provided
$\md{\tan\phi}\le K$, but it does not hold for part (iv)
of Theorem~\ref{up5thm3}.
\end{proof}

We can now show $N_\pm$ are {\it nonsingular\/} SL 3-folds
on their interiors~$N_\pm^\circ$. 

\begin{thm} In the situation of Definition \ref{up6def2}, $N_\pm$
are $\U(1)$-invariant SL\/ $3$-folds with boundaries given in
\eq{up6eq3} and\/ \eq{up6eq4}, and are nonsingular except possibly
at the two boundary points $\bigl(0,0,x_j+iu(x_j,0)\bigr)$ for
$j=1,2$. In particular, the interior $\ga^\circ$ of\/ $\ga$ lies
in the interiors of\/ $N_\pm$, and\/ $N_\pm$ are nonsingular there.
\label{up6thm1}
\end{thm}

\begin{proof} We give the proof for $N_+$. By Definition
\ref{up4def3}, except along the $x$-axis $u,v$ are $C^2$ in
$S$ and real analytic in $S^\circ$, and satisfy \eq{up4eq2}.
Therefore by Proposition \ref{up4prop1}, $N_+$ is a nonsingular
$\U(1)$-invariant SL 3-fold except possibly along $\ga$, where
it may be singular. Also, by Proposition \ref{up6prop6}, $N_+$
is an SL {\it rectifiable current}, and $\ga^\circ$ lies in the
interior $N^\circ_+$, not the boundary.

Let ${\bf z}\in\ga^\circ$. Then by Theorem \ref{up3thm3} there
exists a {\it tangent cone} $C$ to $N_+$ at $\bf z$, which by
Proposition \ref{up6prop7} is of the form $\Pi^\phi_+$ with
multiplicity 1 for some $\phi\in(-\frac{\pi}{2},\frac{\pi}{2})$.
So by Theorem \ref{up3thm4} and Proposition \ref{up5prop4} the
{\it density\/} $\Th(N_+,{\bf z})=1$. Therefore $N_+$ is a
{\it nonsingular\/} SL 3-fold near $\bf z$ by Theorem \ref{up3thm5}.
Thus the only possible singular points of $N_+$ are the endpoints
$\bigl(0,0,x_j+iu(x_j,0)\bigr)$ of~$\ga$.
\end{proof}

Combining Example \ref{up6ex} and Theorem \ref{up6thm1} gives
large families of {\it nonsingular} SL 3-folds invariant under
the $\U(1)$-action \eq{up1eq1}, where the $\U(1)$-action has
{\it nontrivial\/} fixed point set $\ga$. As in Harvey and
Lawson \cite[Th.~III.2.7]{HaLa}, a nonsingular SL 3-fold is
{\it real analytic}, so $N_+$ is real analytic near $\ga^\circ$
by Theorem \ref{up6thm1}. We easily deduce:

\begin{cor} $\ga^\circ$ is a nonsingular real analytic curve in~$\C^3$.
\label{up6cor1}
\end{cor}

Now it has been shown by Bryant \cite[Th.s 2 \& 3]{Brya} that given any
real analytic curve $\de$ in $\bigl\{(0,0,z_3):z_3\in\C\bigr\}$, there
are locally exactly two $\U(1)$-invariant SL 3-folds containing~$\de$.

\begin{thm}[Bryant] Let\/ $\de$ be a nonsingular real analytic curve
lying in the subset\/ $\bigl\{(0,0,z_3):z_3\in\C\bigr\}$ of\/ $\C^3$.
Then there locally exist exactly two SL\/ $3$-folds in $\C^3$
containing $\de$ invariant under the $\U(1)$-action \eq{up1eq1}, which
are exchanged by the involution~$(z_1,z_2,z_3)\mapsto(-z_1,z_2,z_3)$.
\label{up6thm2}
\end{thm}

Bryant's proof involves expanding the $\U(1)$-invariant SL 3-fold
$N$ containing $\de$ as a power series in suitable coordinates.
For the first nontrivial term there are two choices locally,
representing the two possible $\U(1)$-invariant special Lagrangian
choices for $TN$ along $\de$. The higher order terms are then
defined uniquely by a recursive formula derived from the special
Lagrangian condition. Finally, Bryant shows that this formal power
series converges near~$\de$.

We shall use Bryant's result to show that singular solutions $u,v$
of \eq{up4eq2} in which $v$ is zero on an open interval in the
$x$-axis have the symmetries $u(x,-y)=u(x,y)$ and $v(x,-y)=-v(x,y)$,
as in Example~\ref{up6ex}.

\begin{thm} Let\/ $S$ be a domain in $\R^2$ invariant under
$(x,y)\mapsto(x,-y)$. Let\/ $u,v\in C^0(S)$ be singular solutions
of \eq{up4eq2}, as in \S\ref{up46}. Suppose there exist\/ $b<c$ in
$\R$ such that\/ $(x,0)\in S$ and\/ $v(x,0)=0$ for all\/ $x\in[b,c]$.
Then $u(x,-y)=u(x,y)$ and\/ $v(x,-y)=-v(x,y)$ for all\/~$(x,y)\in S$.
\label{up6thm3}
\end{thm}

\begin{proof} Let $T\subset S$ be a subdomain invariant under
$(x,y)\mapsto(x,-y)$ such that the intersection of $T$ with the
$x$-axis is $\bigl\{(x,0):x\in[b,c]\bigr\}$. Define
\begin{align*}
N_+=\bigl\{(z_1,z_2,z_3)\in\C^3:\,&
z_1z_2=v(x,y)+iy,\quad z_3=x+iu(x,y),\\
&\md{z_1}=\md{z_2},\quad (x,y)\in T,\quad y\ge 0\bigr\},
\end{align*}
and define $N_-$ the same way, but with $y\le 0$ rather than
$y\ge 0$. Then Theorem \ref{up6thm1} shows that $N_\pm$ are
compact $\U(1)$-invariant SL 3-folds with boundary, which
contain and are nonsingular along the nonsingular real
analytic curve~$\de=\bigl\{(0,0,x+iu(x,0)):x\in(b,c)\bigr\}$.

By Theorem \ref{up6thm2} there are locally exactly two
$\U(1)$-invariant SL 3-folds containing $\de$, which must
be $N_\pm$, and thus $N_\pm$ are exchanged by $(z_1,z_2,z_3)
\mapsto(-z_1,z_2,z_3)$ near $\de$. This implies that
$u(x,-y)=u(x,y)$ and $v(x,-y)=-v(x,y)$ in $T$ near the
$x$-axis. But $u,v$ are real analytic in $S^\circ$ except
at the $x$-axis, and continuous in $S$, so it easily follows
that $u(x,-y)=u(x,y)$ and $v(x,-y)=-v(x,y)$ in~$S$.
\end{proof}

Later we will use this to prove that $\U(1)$-invariant SL 3-folds
have only {\it isolated\/} singularities, under mild conditions.

\subsection{Inequalities on $v$ for singular solutions of \eq{up4eq2}}
\label{up63}

Next we generalize Proposition \ref{up4prop5} to singular solutions
of~\eq{up4eq2}.

\begin{thm} Let\/ $S$ be a strictly convex domain in $\R^2$ invariant
under the involution $(x,y)\mapsto(x,-y)$, and let\/ $u,v\in C^0(S)$
and\/ $\hat u,\hat v\in C^0(S)$ be singular solutions of \eq{up4eq2}, as in
\S\ref{up46}. If\/ $v<\hat v$ on $\pd S$ then $v<\hat v$ on~$S$.
\label{up6thm4}
\end{thm}

\begin{proof} First we prove the weaker statement that if $v<\hat v$
on $\pd S$ then $v\le \hat v$ on $S$. By part (iv) of Definition
\ref{up4def3}, we may write $v,\hat v$ as the limits in $C^0(S)$ as
$a\ra 0_+$ of solutions $v_a,\hat v_a\in C^2(S)$ of \eq{up4eq9} for
$a\in(0,1]$. Since $v_a\ra v$ and $\hat v_a\ra \hat v$ in $C^0(S)$ as
$a\ra 0_+$, as $v<\hat v$ on $\pd S$ it follows that $v_a<\hat v_a$ on
$\pd S$ for small $a\in(0,1]$. Thus $v_a<\hat v_a$ on $S$ for small
$a\in(0,1]$ by Proposition \ref{up4prop5}. Taking the limit
$a\ra 0_+$ then shows that $v\le \hat v$ in~$S$.

Choose $\phi'\in C^{k+2,\al}(\pd S)$ such that $v<\phi'<\hat v$
on $\pd S$ and $\phi'(x,0)\ne 0$ for points $(x,0)\in\pd S$.
Let $u',v'\in C^0(S)$ be a singular solution of \eq{up4eq2}
with $v'\vert_{\pd S}=\phi'$, which exists by Theorem
\ref{up4thm8}. Since $v<v'<\hat v$ on $\pd S$ we see that $v<v'<\hat v$
near $\pd S$, from above. Thus we can choose a slightly smaller
domain $T\subset S^\circ$ such that $v<v'<\hat v$ on $S\sm T^\circ$,
and in particular on~$\pd T$.

If $t\in\R$ is small then $(x+t,y)\in S$ whenever $(x,y)\in T$.
Define $v_t'\in C^0(T)$ by $v_t'(x,y)=v'(x+t,y)$. Then $v_t'$
is a weak solution of \eq{up4eq9} with $a=0$. Also, since $v_t'$
depends continuously on $t$ we see that there exists $\ep>0$ such
that $v_t'$ is well-defined and $v<v_t'<\hat v$ on $\pd T$ for
$t\in(-\ep,\ep)$. Thus, by the first part of the proof we have
$v\le v_t'\le \hat v$ on $T$ for all~$t\in(-\ep,\ep)$.

Suppose for a contradiction that $v(r,s)=\hat v(r,s)$ for some
$(r,s)\in S$. Then $(r,s)\in T^\circ$, as $v<v'<\hat v$ on
$S\sm T^\circ$. As $v\le v_t'\le \hat v$ on $T$ for all $t\in(-\ep,\ep)$
we have $v_t'(r,s)=v(r,s)$ for all $t\in(-\ep,\ep)$. Hence
$v'(r+t,s)=v(r,s)$ for all~$t\in(-\ep,\ep)$.

Consider the two cases (a) $s=0$ and $v(r,0)=0$ and (b) $s\ne 0$ or
$v(r,s)\ne 0$. In case (a) we have shown that $v'$ is zero on an
open interval $(r-\ep,r+\ep)$ in the $x$-axis. So Theorem \ref{up6thm3}
gives $v'(x,y)=-v'(x,-y)$ in $S$. But then $v'(x,0)\equiv 0$,
contradicting $v'(x,0)=\phi'(x,0)\ne 0$ for~$(x,0)\in\pd S$.

In case (b), as $v'$ is real analytic where it is nonsingular we
see that $v'(x,s)=v(r,s)$ for all $x$ such that $(x,s)\in S$. In
particular, this implies that $\phi'(x_1,s)=\phi'(x_2,s)$ for
the two points $(x_1,s)$, $(x_2,s)$ in $\pd S$ of the form $(x,s)$.
Choosing $\phi'$ so that $\phi'(x_1,s)\ne\phi'(x_2,s)$ we again
derive a contradiction.
\end{proof}

\subsection{Positivity of winding numbers and multiplicities}
\label{up64}

Using the results of \S\ref{up62}--\S\ref{up63}, we shall
complete the arguments of \S\ref{up61} to show that singular
solutions intersect with {\it positive} multiplicity.

\begin{prop} Let\/ $u,v\in C^0(S)$ and\/ $\hat u,\hat v\in C^0(S)$ be
singular solutions of\/ \eq{up4eq2} in a domain $S$ in $\R^2$, as
in \S\ref{up46}, such that\/ $(u,v)\ne(\hat u,\hat v)$ at every point
of\/ $\pd S$. Suppose $(u,v)-(\hat u,\hat v)$ has a zero $(b,c)$ in
$S^\circ$. Then the winding number of\/ $(u,v)-(\hat u,\hat v)$
about\/ $0$ along $\pd S$ is a positive integer.
\label{up6prop8}
\end{prop}

\begin{proof} If $(b,c)$ is a nonsingular zero then by Proposition
\ref{up6prop3} the multiplicity of $(b,c)$ is a positive integer
$k$, and by Proposition \ref{up6prop2} the winding number of $(u,v)-
(\hat u,\hat v)$ about 0 along $\pd S$ is at least $k$, so it is a
positive integer, as we have to prove.

So let $(b,c)$ be a singular zero, giving $c=0$ and $v(b,0)=\hat
v(b,0)=0$. As $(u,v)\ne(\hat u,\hat v)$ along $\pd S$, there exists
a small $\de>0$ such that if $\al\in(-\de,\de)$ then $(u+\al,v)
\ne(\hat u,\hat v)$ along $\pd S$, and the winding numbers of
$(u+\al,v)-(\hat u,\hat v)$ and $(u,v)-(\hat u,\hat v)$ about 0
along $\pd S$ are the same.

We shall show that one can choose $\al\in(-\de,\de)$ such that
$(u+\al,v)-(\hat u,\hat v)$ has a zero $(b',c')$ near $(b,0)$ in
$S^\circ$ with positive multiplicity. Then by Proposition
\ref{up6prop2} the winding number of $(u+\al,v)-(\hat u,\hat v)$
about 0 along $\pd S$ is a positive integer, and thus the winding
number of $(u,v)-(\hat u,\hat v)$ about $0$ along $\pd S$ is a
positive integer, as we want.

Let $\ep>0$ be small enough that $\ovB_\ep(b,0)$ lies in $S^\circ$.
Suppose $v\ne\hat v$ at every point of $\ga_\ep(b,0)$. Then either
$v<\hat v$ on $\ga_\ep(b,c)$ or $v>\hat v$ on $\ga_\ep(b,0)$ by
continuity. Now $\ga_\ep(b,0)$ is the boundary of a strictly
convex domain $\ovB_\ep(b,0)$ invariant under $(x,y)\mapsto(x,-y)$.
Applying Theorem \ref{up6thm4} shows that $v<\hat v$ or $v>\hat v$
on $\ovB_\ep(b,0)$. But this contradicts $v=\hat v=0$ at~$(b,0)$.

Thus, for any small $\ep$ there exists a point $(b',c')$ on
$\ga_\ep(b,0)$ where $v=\hat v$. Set $\al=\hat u(b',c')-u(b',c')$.
Then $(u+\al,v)-(\hat u,\hat v)$ is zero at $(b',c')$. Also, since
$u=\hat u$ at $(b,0)$ and $u,\hat u$ are continuous we see that
$\al\ra 0$ as $\ep\ra 0_+$, and so $\al\in(-\de,\de)$ if $\ep$
is small enough.

For each sufficiently small $\ep>0$ we have constructed
$\al\in(-\de,\de)$ and a zero $(b',c')$ of $(u+\al,v)-(\hat u,
\hat v)$ on the circle $\ga_\ep(b,c)$. Consider the two cases
\begin{itemize}
\setlength{\itemsep}{0pt}
\setlength{\parsep}{0pt}
\item[(a)] $(b',c')$ is a {\it nonsingular} zero for some small
$\ep>0$, and
\item[(b)] $(b',c')$ is a {\it singular} zero for all small~$\ep>0$.
\end{itemize}
In case (a), for this $\ep$ Proposition \ref{up6prop3} shows that
the multiplicity of $(b',c')$ is positive, and we are finished.

In case (b) we must have $c'=0$, so that $b'=b\pm\ep$, and so for
all small $\ep>0$ we have either $v(b+\ep,0)=\hat v(b+\ep,0)=0$ or
$v(b-\ep,0)=\hat v(b-\ep,0)=0$. It is easy to show using continuity
of $v,\hat v$ that $v,\hat v$ are zero on a positive length closed
interval in the $x$-axis, as in \S\ref{up62}. Theorem \ref{up6thm3}
then shows that $u(x,-y)\equiv u(x,y)$, $v(x,-y)\equiv -v(x,y)$,
$\hat u(x,-y)\equiv\hat u(x,y)$ and $\hat v(x,-y)\equiv-\hat v(x,y)$
near $(b,0)$ in~$S$.

By Proposition \ref{up6prop7} the $\U(1)$-invariant SL 3-folds
$N,\hat N$ corresponding to $u+\al,v$ and $\hat u,\hat v$ have
a common singular point at ${\bf z}=\bigl(0,0,b'+iu(b',0)\bigr)$,
with tangent cones $\Pi^\phi_+\cup\Pi^\phi_-$ and $\Pi^{\smash{
\hat\phi}}_+\cup\Pi^{\smash{\hat\phi}}_-$. For generic $\ep$
it is easy to see that $\phi\ne\hat\phi$, and one can then
show using tangent cones that the multiplicity of $(b',0)$
is 1, which is positive. This completes the proof.
\end{proof}

Since the multiplicity of an isolated zero of $(u,v)-(\hat u,\hat v)$
given in Definition \ref{up6def1} is the winding number of
$(u,v)-(\hat u,\hat v)$ about 0 along the boundary of a domain
$\ovB_\ep(b,c)$ containing a zero $(b,c)$ of $(u,v)-(\hat u,\hat v)$,
we deduce:

\begin{cor} Let\/ $u,v\in C^0(S)$ and\/ $\hat u,\hat v\in C^0(S)$
be singular solutions of\/ \eq{up4eq2} in a domain $S$ in $\R^2$.
Then the multiplicity of any isolated zero $(b,c)$ of\/ $(u,v)
-(\hat u,\hat v)$ in $S^\circ$ is a positive integer.
\label{up6cor2}
\end{cor}

\subsection{Counting formulae using winding numbers}
\label{up65}

We can now generalize Theorem \ref{up4thm2} to the singular case.

\begin{thm} Let\/ $S$ be a domain in $\R^2$, and let\/ $u,v$
and\/ $\hat u,\hat v$ be singular solutions of\/ \eq{up4eq2} in $C^0(S)$,
with\/ $(u,v)\neq(\hat u,\hat v)$ at every point of\/ $\pd S$. Then
$(u,v)-(\hat u,\hat v)$ has at most finitely many zeroes in $S^\circ$,
all isolated. Suppose that there are $n$ zeroes with multiplicities
$k_1,\ldots,k_n$. Then the winding number of\/ $(u,v)-(\hat u,\hat v)$
about\/ $0$ along $\pd S$ is~$\sum_{i=1}^nk_i$.
\label{up6thm5}
\end{thm}

\begin{proof} Let the winding number of $(u,v)-(\hat u,\hat v)$
about $0$ along $\pd S$ be $k$. Then $k\ge 0$ by Proposition
\ref{up6prop1}. Suppose there exist $k+1$ isolated zeroes of
$(u,v)-(\hat u,\hat v)$ in $S^\circ$. The multiplicity of each is at
least 1 by Corollary \ref{up6cor2}, so Proposition \ref{up6prop2}
shows that the winding number of $(u,v)-(\hat u,\hat v)$ about 0
along $\pd S$ is at least $k+1$, a contradiction. Hence there
can be at most $k$ isolated zeroes of $(u,v)-(\hat u,\hat v)$ in
$S^\circ$, and in particular there are finitely many.

Next we shall show that there are no {\it nonisolated\/} zeroes of
$(u,v)-(\hat u,\hat v)$. Let $Z$ be the set of nonisolated zeroes of
$(u,v)-(\hat u,\hat v)$. Then $Z$ is a closed subset of $S^\circ$,
as the zero set of $(u,v)-(\hat u,\hat v)$ is closed and $Z$ is
this with finitely many isolated points removed. By Proposition
\ref{up6prop3}, any nonsingular zero is isolated, so $Z$ consists
of singular zeroes, which therefore lie on the $x$-axis, and
$u=\hat u$ and $v=\hat v=0$ on~$Z$.

Each connected component of $Z$ is therefore a closed, connected
subset of the $x$-axis in $S^\circ$. By elementary topology, it
must be either a point or a closed interval of positive length.
Divide into two cases:
\begin{itemize}
\setlength{\itemsep}{0pt}
\setlength{\parsep}{0pt}
\item[(a)] a connected component of $Z$ is an interval
$[\al,\be]\t\{0\}$ for $\al<\be$, and
\item[(b)] all connected components of $Z$ are points.
\end{itemize}

In case (a), as $v_j=0$ on $(\al,\be)\t\{0\}$ Theorem \ref{up6thm3}
shows that $u_j(x,-y)=u_j(x,y)$ and $v_j(x,-y)=-v_j(x,y)$ for $j=1,2$
near $(\al,\be)\t\{0\}$. As $u=\hat u$ and $v=\hat v=0$ on $(\al,\be)\t
\{0\}$, Theorem \ref{up6thm2} implies that $(u,v)\equiv(\hat u,\hat v)$
near $(\al,\be)\t\{0\}$. But the $(u_j,v_j)$ are real analytic where
they are nonsingular in $S^\circ$, and continuous in $S$, so $(u,v)
\equiv(\hat u,\hat v)$ in $S$. This contradicts $(u,v)\ne(\hat u,
\hat v)$ on $\pd S$, and excludes case~(a).

In case (b), suppose $Z$ has at least $k+1$ connected components.
As $Z$ is a closed subset of the $x$-axis whose connected components
are points, it is easy to see that we can find $k+1$ disjoint closed
discs $D_1,\ldots,D_{k+1}$ in $S^\circ$ with centres on the $x$-axis,
such that $(u,v)-(\hat u,\hat v)$ has no zeroes on $\pd D_j$ and
$D_j^\circ$ contains a connected component of $Z$. By Proposition
\ref{up6prop8}, the winding number of $(u,v)-(\hat u,\hat v)$ about
0 along $\pd D_i$ is a positive integer.

Thus the sum of the winding numbers of $(u,v)-(\hat u,\hat v)$ about 0
along $\pd D_i$ for $i=1,\ldots,k+1$ is greater than $k$, the winding
number about 0 along $\pd S$. Reasoning as in the proof of Proposition
\ref{up6prop2} with $T=S\sm\bigcup_{i=1}^{k+1}D_i^\circ$ then gives a
contradiction. Hence $Z$ has at most $k$ connected components, all
single points, and so all isolated zeroes. Thus $Z=\emptyset$ by
definition, and there are finitely many zeroes of $(u,v)-(\hat u,
\hat v)$, all isolated, as we have to prove.

So let the zeroes of $(u,v)-(\hat u,\hat v)$ be $(b_1,c_1),\ldots,
(b_n,c_n)$, with multiplicities $k_1,\ldots,k_n$. Define $\ep_1,
\ldots,\ep_n$ and $T$ as in the proof of Proposition \ref{up6prop2}.
Then the winding number of $(u,v)-(\hat u,\hat v)$ about 0 along $\pd T$
is $k-\sum_{i=1}^nk_i$. But there are no zeroes of $(u,v)-(\hat u,\hat v)$
in $T$, so by the proof of Theorem \ref{up4thm2} in \cite{Joyc4} we see
that this winding number is zero, and hence $k=\sum_{i=1}^nk_i$,
completing the proof.
\end{proof}

\subsection{A criterion for isolated zeroes}
\label{up66}

Without assuming that $(u,v)\ne(\hat u,\hat v)$ at every point of
$\pd S$, we can generalize Theorem \ref{up6thm5} to show that
$(u,v)-(\hat u,\hat v)$ has isolated zeroes in~$S^\circ$.

\begin{thm} Let\/ $S$ be a domain in $\R^2$, and let\/ $u,v$ and\/
$\hat u,\hat v$ be singular solutions of\/ \eq{up4eq2} in $C^0(S)$,
such that\/ $(u,v)\not\equiv(\hat u,\hat v)$. Then there are at most
countably many zeroes of\/ $(u,v)-(\hat u,\hat v)$ in $S^\circ$, all
isolated.
\label{up6thm6}
\end{thm}

\begin{proof} We may surround each isolated zero of $(u,v)-
(\hat u,\hat v)$ in $S^\circ$ by a small disc in $S^\circ$, such that
the collection of such discs is disjoint. As $S^\circ$ can contain
only countably many disjoint discs, there are only countably many
isolated zeroes in $S^\circ$. Since $(u,v)\not\equiv(\hat u,\hat v)$,
Proposition \ref{up6prop3} shows that nonsingular zeroes of
$(u,v)-(\hat u,\hat v)$ in $S^\circ$ are isolated.

Let $Z$ be the set of nonisolated zeroes of $(u,v)-(\hat u,\hat v)$
in $S^\circ$. Then $Z$ is closed in $S^\circ$ (though not necessarily
in $S$) and is a subset of the $x$-axis in $S^\circ$, by the arguments
in the proof of Theorem \ref{up6thm5}. Thus, by elementary topology
each connected component of $Z$ is either a point, or an interval
of positive length. Case (a) in the proof of Theorem \ref{up6thm5}
shows that if a component is an interval of positive length then
$(u,v)\equiv(\hat u,\hat v)$, a contradiction. So, the connected
components of $Z$ are all points.

Given any $(b,c)\in Z$, we can find a subdomain $T\subset S^\circ$
such that $(b,c)\in T^\circ$ and $(u,v)\neq(\hat u,\hat v)$ at every
point of $\pd T$. We must ensure that $\pd T$ avoids the zeroes
of $(u,v)-(\hat u,\hat v)$. There are only countably many isolated
zeroes, so a generic $T$ has none on $\pd T$, and we can also
arrange for $\pd T$ to avoid $Z$, as $Z$ is closed in the $x$-axis in
$S^\circ$ and its connected components are points. Applying Theorem
\ref{up6thm5} on $T$ shows that all zeroes of $(u,v)-(\hat u,\hat v)$
in $T^\circ$ are isolated, a contradiction as $(b,c)\in T^\circ$.
Thus $Z=\emptyset$, and there are no nonisolated zeroes.
\end{proof}

If there are infinitely many zeroes of $(u,v)-(\hat u,\hat v)$ in
$S^\circ$, then they have a limit point in $\pd S$ which is a
nonisolated zero. For all $k\ge 1$ one can write down examples
of holomorphic functions on domains $S$ in $\C$ which are $C^k$
on $\pd S$, but which have countably many zeroes in $S^\circ$
converging to a limit in $\pd S$. Given the strong analogy between
\eq{up4eq2} and the Cauchy--Riemann equations, it is likely
that there exist examples in Theorem \ref{up6thm6} in which
$(u,v)-(\hat u,\hat v)$ has infinitely many zeroes in~$S^\circ$.

From Proposition \ref{up4prop2} and Theorem \ref{up6thm6} we see
that if $u,v$ and $\hat u,\hat v$ satisfy \eq{up4eq2} or \eq{up4eq3}
and $(b,c)\in S^\circ$ is a zero of $(u,v)-(\hat u,\hat v)$, then
{\it either} $(b,c)$ is an isolated zero with a unique multiplicity,
{\it or} $(u,v)\equiv(\hat u,\hat v)$. In effect, $(u,v)\equiv
(\hat u,\hat v)$ means that $(b,c)$ is a zero with `multiplicity $\iy$'.
So we make the following convention, which will simplify the
discussion in~\S\ref{up9}.

\begin{dfn} Let $S$ be a domain in $\R^2$, and let $u,v$ and
$\hat u,\hat v$ be singular solutions of \eq{up4eq2} in $C^0(S)$,
or solutions of \eq{up4eq3} in $C^1(S)$ for some $a\ne 0$. We
say that $(u,v)-(\hat u,\hat v)$ {\it has a zero of multiplicity at
least\/} $k$ at $(b,c)\in S^\circ$ if either $(u,v)-(\hat u,\hat v)$
has an {\it isolated\/} zero of multiplicity at least $k$ at
$(b,c)$, or~$(u,v)\equiv(\hat u,\hat v)$.
\label{up6def3}
\end{dfn}

\subsection{Counting formulae using potentials}
\label{up67}

Here is a generalization of Theorem \ref{up4thm6} to the singular case.

\begin{thm} Suppose $S$ is a strictly convex domain in $\R^2$
invariant under $(x,y)\mapsto(x,-y)$, and\/ $\phi_1,\phi_2\in
C^{3,\al}(\pd S)$ for some $\al\in(0,1)$. Let\/ $u_j,v_j\in
C^0(S)$ be the singular solution of\/ \eq{up4eq2} constructed
in Theorem \ref{up4thm4} from~$\phi_j$.

Suppose $\phi_1-\phi_2$ has exactly $l$ local maxima and\/ $l$
local minima on $\pd S$. Then $(u_1,v_1)-(u_2,v_2)$ has finitely
many zeroes in $S^\circ$, all isolated. Let there be $n$ zeroes
in $S^\circ$ with multiplicities $k_1,\ldots,k_n$.
Then~$\sum_{i=1}^nk_i\le l-1$.
\label{up6thm7}
\end{thm}

\begin{proof} For $j=1,2$ and $a\in(0,1]$, let $f_{j,a}\in
C^{3,\al}(S)$ be the solution of \eq{up4eq7} in $C^{k+2,\al}(S)$ with
$f_{j,a}\vert_{\pd S}=\phi_j$ given in Theorem \ref{up4thm3}, and set
$u_{j,a}=\frac{\pd}{\pd y}f_{j,a}$ and $v_{j,a}=\frac{\pd}{\pd x}
f_{j,a}$. Then Theorem \ref{up4thm6} shows that for all $a\in(0,1]$
there are no more than $l-1$ zeroes of $(u_{1,a},v_{1,a})-(u_{2,a},
v_{2,a})$ in $S^\circ$, counted with multiplicity.

Let $T\subset S$ be a subdomain with no zeroes of $(u_1,v_1)-(u_2,v_2)$
on $\pd T$. Now $(u_{j,a},v_{j,a})\ra(u_j,v_j)$ in $C^0(S)$ as $a\ra 0_+$
by Theorem \ref{up4thm5}. Thus, for small $a\in(0,1]$ there are no
zeroes of $(u_{1,a},v_{1,a})-(u_{2,a},v_{2,a})$ on $\pd T$, and the
winding numbers of $(u_1,v_1)-(u_2,v_2)$ and $(u_{1,a},v_{1,a})-
(u_{2,a},v_{2,a})$ about 0 along $\pd T$ are equal.

But from above there are no more than $l-1$ zeroes of
$(u_{1,a},v_{1,a})-(u_{2,a},v_{2,a})$ in $T^\circ$, counted with
multiplicity. Hence by Theorem \ref{up4thm2} the winding number
of $(u_{1,a},v_{1,a})-(u_{2,a},v_{2,a})$ about 0 along $\pd T$
is no more than $l-1$. As this is the winding number of
$(u_1,v_1)-(u_2,v_2)$, Theorem \ref{up6thm5} shows that there
are no more than $l-1$ zeroes of $(u_1,v_1)-(u_2,v_2)$ in
$T^\circ$, counted with multiplicity.

Since $\phi_1-\phi_2$ is not constant we have $(u_1,v_1)\not\equiv
(u_2,v_2)$, and so by Theorem \ref{up6thm6} there are at most
countably many zeroes of $(u_1,v_1)-(u_2,v_2)$ in $S^\circ$,
all isolated. If there were infinitely many we could choose
$T$ above to contain at least $l$ zeroes, giving a contradiction.
Thus there are only finitely many, we can choose $T$ to contain
them all, and the result follows.
\end{proof}

We can also generalize \cite[Th.~7.10]{Joyc4} to include the
number of nonsingular zeroes on $\pd S$ in the inequality.
However, there may be a problem with including the singular
zeroes of $(u_1,v_1)-(u_2,v_2)$ on~$\pd S$.

\section{Special Lagrangian fibrations}
\label{up7}

In \cite{SYZ}, Strominger, Yau and Zaslow proposed an explanation
of Mirror Symmetry between Calabi--Yau 3-folds $X,\hat X$ in terms
of the existence of dual {\it special Lagrangian fibrations}
$f:X\ra B$, $\hat f:\hat X\ra B$ over the same base space $B$, a
real 3-manifold. This is known as the {\it SYZ Conjecture}. These
fibrations $f,\hat f$ must necessarily contain singular fibres, which
are a source of many of the mathematical difficulties surrounding the
SYZ Conjecture, as the singularities of SL 3-folds are not yet well
understood.

We will now use our results to construct large families of special
Lagrangian fibrations of open subsets of $\C^3$ invariant under
the $\U(1)$-action \eq{up1eq1}, including singular fibres. These
can serve as local models for singularities of SL fibrations of
(almost) Calabi--Yau 3-folds. In \cite{Joyc6} we will discuss
these fibrations at much greater length, and draw some conclusions
on the singular behaviour of SL fibrations of (almost) Calabi--Yau
3-folds, and on how to best formulate the SYZ Conjecture.

\begin{dfn} Let $S$ be a strictly convex domain in $\R^2$ invariant
under $(x,y)\mapsto(x,-y)$, let $U$ be an open set in $\R^3$, and
$\al\in(0,1)$. Suppose $\Phi:U\ra C^{3,\al}(\pd S)$ is a continuous
map such that if $(a,b,c)\ne(a,b',c')$ in $U$ then $\Phi(a,b,c)-
\Phi(a,b',c')$ has exactly one local maximum and one local minimum
in~$\pd S$.

Let ${\bs\al}=(a,b,c)\in U$. If $a\ne 0$, let $f_{\bs\al}\in
C^{3,\al}(S)$ be the unique solution of \eq{up4eq7} with
$f_{\bs\al}\vert_{\pd S}=\Phi({\bs\al})$, which exists by
Theorem \ref{up4thm3}. If $a=0$, let $f_{\bs\al}\in
C^1(S)$ be the unique weak solution of \eq{up4eq8} with
$f_{\bs\al}\vert_{\pd S}=\Phi({\bs\al})$, which exists by
Theorem \ref{up4thm4}. Define $u_{\bs\al}=\frac{\pd f_{\bs\al}}{\pd y}$
and $v_{\bs\al}=\frac{\pd f_{\bs\al}}{\pd x}$. Then
$(u_{\bs\al},v_{\bs\al})$ is a solution of \eq{up4eq3} if
$a\ne 0$, and a singular solution of \eq{up4eq2} if $a=0$.
Also $u_{\bs\al},v_{\bs\al}$ depend continuously on ${\bs\al}\in U$
in $C^0(S)$, by Theorem~\ref{up4thm5}.

For each ${\bs\al}=(a,b,c)$ in $U$, define $N_{\bs\al}$ in $\C^3$ by
\e
\begin{split}
N_{\bs\al}=\bigl\{(z_1,z_2,z_3)\in\C^3:\,&
z_1z_2=v_{\bs\al}(x,y)+iy,\quad z_3=x+iu_{\bs\al}(x,y),\\
&\ms{z_1}-\ms{z_2}=2a,\quad (x,y)\in S^\circ\bigr\}.
\end{split}
\label{up7eq1}
\e
Then $N_{\bs\al}$ is a noncompact SL 3-fold without boundary in $\C^3$,
which is nonsingular if $a\ne 0$, by Proposition~\ref{up4prop1}.
\label{up7def}
\end{dfn}

We shall show that the $N_{\bs\al}$ are the fibres of an SL
fibration. This is one of the main results of the paper,
which will be the central tool in~\cite{Joyc6}.

\begin{thm} In the situation of Definition \ref{up7def}, if\/
${\bs\al}\ne{\bs\al}'$ in $U$ then $N_{\bs\al}\cap N_{{\bs\al}'}
=\emptyset$. There exists an open set\/ $V\subset\C^3$ and a
continuous, surjective map $F:V\ra U$ such that\/ $F^{-1}({\bs\al})
=N_{\bs\al}$ for all\/ ${\bs\al}\in U$. Thus, $F$ is a special
Lagrangian fibration of\/ $V\subset\C^3$, which may include
singular fibres.
\label{up7thm}
\end{thm}

\begin{proof} Let ${\bs\al}=(a,b,c)$ and ${\bs\al}'=(a',b',c')$ be
distinct elements of $U$. As $\ms{z_1}-\ms{z_2}=2a$ on $N_{\bs\al}$
and $\ms{z_1}-\ms{z_2}=2a'$ on $N_{{\bs\al}'}$, clearly $N_{\bs\al}
\cap N_{{\bs\al}'}=\emptyset$ if $a\ne a'$. So suppose $a=a'$. Then
$\Phi(a,b,c)-\Phi(a,b',c')$ has exactly one local maximum and one
local minimum in $\pd S$, by the condition in Definition \ref{up7def}.
Hence by Theorem \ref{up4thm6} when $a\ne 0$ and Theorem \ref{up6thm7}
when $a=0$ we see that $(u_{\bs\al},v_{\bs\al})-(u_{{\bs\al}'},
v_{{\bs\al}'})$ has no zeroes in $S^\circ$, and thus~$N_{\bs\al}\cap
N_{{\bs\al}'}=\emptyset$.

Let $V=\bigcup_{{\bs\al}\in U}N_{\bs\al}$, and define $F:V\ra U$ by
$F({\bf z})={\bs\al}$ if ${\bf z}\in N_{\bs\al}$. As $N_{\bs\al}\cap
N_{{\bs\al}'}=\emptyset$ when ${\bs\al}\ne{\bs\al}'$ this map $F$ is
well-defined, and clearly $F^{-1}({\bs\al})=N_{\bs\al}$. As
$N_{\bs\al}\ne\emptyset$ for all ${\bs\al}\in U$, we see that $F$ is
surjective. It remains only to show that $F$ is continuous, and $V$
is open.

Fix ${\bs\al}'=(a',b',c')\in U$. As $U$ is open in $\R^3=\R\t\R^2$
there exist open $A\subset\R$ and $B\subset\R^2$ with $(a',b',c')
\in A\t B\subset U$. Let $(x,y)\in S^\circ$, and for each $a\in A$
define $\Psi_{a,x,y}:B\ra\R^2$ by $\Psi_{a,x,y}(b,c)=\bigl(
u_{\bs\al}(x,y),v_{\bs\al}(x,y)\bigr)$ for ${\bs\al}=(a,b,c)$.
Then $\Psi_{a,x,y}$ is continuous and depends continuously on
$a,x,y$, as $u_{\bs\al},v_{\bs\al}$ are continuous and depend
continuously on $\bs\al$ from Definition \ref{up7def}. Also
$\Psi_{a,x,y}$ is injective, as $N_{\bs\al}\cap N_{{\bs\al}'}
=\emptyset$ for~${\bs\al}\ne{\bs\al}'\in U$.

Since $\Psi_{a,x,y}$ is continuous and injective, it follows by
elementary topology in $\R^2$ that $W_{a,x,y}=\Psi_{a,x,y}(B)$
is open in $\R^2$, and $\Psi_{a,x,y}^{-1}:W_{a,x,y}\ra B$ is
continuous. Furthermore, $W_{a,x,y}$ and $\Psi_{a,x,y}^{-1}$
depend continuously on~$a,x,y$.

Now if ${\bf z}=(z_1,z_2,z_3)\in V$ and $F({\bf z})=(a,b,c)\in
A\t B$, then $2a=\ms{z_1}-\ms{z_2}$, $x=\Re(z_3)$, $u=\Im(z_3)$,
$v=\Re(z_1z_2)$ and $y=\Im(z_1z_2)$. Therefore
\begin{equation*}
a=\ha(\ms{z_1}-\ms{z_2}),\quad
(b,c)=\Psi_{(\ms{z_1}-\ms{z_2})/2,\Re(z_3),\Im(z_1z_2)}^{-1}
\bigl(\Im(z_3),\Re(z_1z_2)\bigr).
\end{equation*}
As $\Psi_{a,x,y}^{-1}$ is continuous and depends continuously
on $a,x,y$, we see that $(a,b,c)$ depends continuously on
$\bf z$. Hence $F:{\bf z}\mapsto(a,b,c)$ is continuous.

Finally we show $V$ is open. Fix $(z_1',z_2',z_3')\in N_{{\bs\al}'}$,
and let $x'=\Re(z_3')$, $u'=\Im(z_3')$, $v'=\Re(z_1'z_2')$ and
$y'=\Im(z_1'z_2')$. Then $\Psi_{a',x',y'}(b',c')=(u',v')$, so
$(u',v')\in W_{a',x',y'}$. As $W_{a,x,y}$ is open and depends
continuously on $a,x,y$, there exist open neighbourhoods $W$
of $(u',v')$ in $\R^2$ and $X$ of $(a',x',y')$ in $A\t S^\circ$
such that if $(a,x,y)\in X$ then $W\subset W_{a,x,y}$. Define
\begin{align*}
Y=\Bigl\{(z_1,z_2,z_3)\in\C^3:\,&\bigl(\Im(z_3),\Re(z_1z_2)\bigr)\in W,\\
&\bigl(\ha(\ms{z_1}-\ms{z_2}),\Re(z_3),\Im(z_1z_2)\bigr)\in X\Bigr\}.
\end{align*}
Then $Y$ is open in $\C^3$, as $W,X$ are open, and
contains~$(z_1',z_2',z_3')$.

We claim that $Y\subset V$. Let ${\bf z}=(z_1,z_2,z_3)\in Y$,
and set $a=\ha(\ms{z_1}-\ms{z_2})$, $x=\Re(z_3)$, $u=\Im(z_3)$,
$v=\Re(z_1z_2)$ and $y=\Im(z_1z_2)$. Then $(u,v)\in W$ and
$(a,x,y)\in X$. Hence $W\subset W_{a,x,y}$, so $(u,v)\in W_{a,x,y}$.
Thus $\Psi_{a,x,y}^{-1}(u,v)=(b,c)$ is well-defined. From the
definitions we see that ${\bs\al}=(a,b,c)$ lies in $U$ and
$(z_1,z_2,z_3)\in N_{\bs\al}$, so that ${\bf z}\in V$. Thus
$Y\subset V$, so that any $(z_1',z_2',z_3')\in V$ has an open
neighbourhood $Y\subset V$, and $V$ is open.
\end{proof}

Note that in \eq{up7eq1} we have chosen to define the $N_{\bs\al}$
over $S^\circ$, so that they are noncompact SL 3-folds without
boundary. The closures $\,\ov{\!N}_{\bs\al}$ are compact SL 3-folds
with boundary, defined over $S$. The main reason for working
over $S^\circ$ rather than $S$ is to avoid difficulties in
proving that $\,\ov{\!N}_{\bs\al}\cap\,\ov{\!N}_{{\bs\al}'}=\emptyset$
if ${\bs\al}\ne{\bs\al}'$ in $U$. The problem is that
$\,\ov{\!N}_{\bs\al},\,\ov{\!N}_{{\bs\al}'}$ may intersect in a boundary
point, lying over~$\pd S$.

If we strengthen the condition on $\Phi$ in Definition \ref{up7def}
to say that if $(a,b,c)\ne(a,b',c')$ in $U$ then $\Phi(a,b,c)-
\Phi(a,b',c')$ has exactly two stationary points, and the second
derivative of $\Phi(a,b,c)-\Phi(a,b',c')$ is nonzero at each, then
we can use \cite[Th.~7.10]{Joyc4} to show that $\,\ov{\!N}_{\bs\al}\cap
\,\ov{\!N}_{{\bs\al}'}$ contains no nonsingular points. But there
remains the possibility, when $a=0$, that $\,\ov{\!N}_{\bs\al}$,
$\,\ov{\!N}_{{\bs\al}'}$ could intersect in a common singular
point lying over some $(x,0)$ in~$\pd S$.

Here is a simple way to produce families $\Phi$ satisfying
Definition~\ref{up7def}.

\begin{ex} Let $S$ be a strictly convex domain in $\R^2$ invariant
under $(x,y)\mapsto(x,-y)$, let $\al\in(0,1)$ and $\phi\in C^{3,\al}
(\pd S)$. Define $U=\R^3$ and $\Phi:\R^3\ra C^{3,\al}(\pd S)$ by
$\Phi(a,b,c)=\phi+bx+cy$. If $(a,b,c)\ne(a,b',c')$ then
$\Phi(a,b,c)-\Phi(a,b',c')=(b-b')x+(c-c')y\in C^\iy(\pd S)$. As
$b-b',c-c'$ are not both zero and $S$ is strictly convex, it
easily follows that $(b-b')x+(c-c')y$ has exactly one local
maximum and one local minimum in $\pd S$. Hence the conditions
of Definition \ref{up7def} hold for $S,U$ and $\Phi$, and so
Theorem \ref{up7thm} defines an open set $V\subset\C^3$ and
a special Lagrangian fibration~$F:V\ra\C^3$.

Let $a,b,c,c'\in\R$. Then $f_{(a,b,c)}=f_{(a,b,c')}+(c-c')y$,
$u_{(a,b,c)}=u_{(a,b,c')}+(c-c')$ and $v_{(a,b,c)}=v_{(a,b,c')}$.
It follows from \eq{up7eq1} that $N_{(a,b,c)}$ is the translation of
$N_{(a,b,c')}$ by $\bigl(0,0,i(c-c')\bigr)$ in $\C^3$. So, changing
the parameter $c$ in $U=\R^3$ just translates the fibres $N_{\bs\al}$
in~$\C^3$.

One can also show that $v_{(a,b,c)}(x,y)\ra\pm\iy$ as $b\ra\pm\iy$,
for fixed $a,c\in\R$ and $(x,y)\in S^\circ$. Combining these facts
about changing $b,c$ and taking $A=\R$ and $B=\R^2$ we find that
$\Psi_{a,x,y}:\R^2\ra\R^2$ is {\it surjective} for all $a\in\R$
and $(x,y)\in S^\circ$, so that $W_{a,x,y}=\R^2$. From this we
easily prove that
\begin{equation*}
V=\Bigl\{(z_1,z_2,z_3)\in\C^3:\bigl(\Re(z_3),\Im(z_1z_2)\bigr)
\in S^\circ\Bigr\}.
\end{equation*}
\label{up7ex}
\end{ex}

This example, and other families of maps $\Phi$ one can readily
construct, generate many special Lagrangian fibrations of open
subsets of~$\C^3$.

\section{A rough classification of singular points}
\label{up8}

We can now use the work of \S\ref{up6} to study singular points of~$u,v$.

\begin{dfn} Let $S$ be a domain in $\R^2$, and $u,v\in C^0(S)$ a
singular solution of \eq{up4eq2}, as in \S\ref{up46}. Suppose
for simplicity that $S$ is invariant under $(x,y)\mapsto(x,-y)$.
Define $u',v'\in C^0(S)$ by $u'(x,y)=u(x,-y)$ and $v'(x,y)=-v(x,-y)$.
Then $u',v'$ is also a singular solution of~\eq{up4eq2}.

A {\it singular point}, or {\it singularity}, of $(u,v)$ is a point
$(b,0)\in S$ with $v(b,0)=0$. Observe that a singularity of $(u,v)$
is automatically a zero of $(u,v)-(u',v')$. Conversely, a zero of
$(u,v)-(u',v')$ on the $x$-axis is a singularity. A singularity
of $(u,v)$ is called {\it isolated} if it is an isolated zero of
$(u,v)-(u',v')$. Define the {\it multiplicity} of an isolated
singularity $(b,0)$ of $(u,v)$ in $S^\circ$ to be the multiplicity
of $(u,v)-(u',v')$ at $(b,0)$, in the sense of Definition \ref{up6def1}.
By Corollary \ref{up6cor2}, this multiplicity is a positive integer.

If $S$ is not invariant under $(x,y)\mapsto(x,-y)$ then we can
still define what it means for a singular point $(b,0)$ to
be {\it isolated}, and the {\it multiplicity} of an isolated
singular point, by restricting to $T=\bigl\{(x,y)\in S:(x,-y)\in
S\bigr\}$ before defining~$u',v'$.
\label{up8def1}
\end{dfn}

From Theorem \ref{up6thm6} we see that a singular solution $(u,v)$ of
\eq{up4eq2} in $S$ either has the symmetries $u(x,-y)\equiv u(x,y)$
and $v(x,-y)\equiv -v(x,y)$, as in \S\ref{up62}, or else its singular
points in $S^\circ$ are all isolated, and so have a well-defined,
positive multiplicity.

\begin{thm} Let\/ $S$ be a domain in $\R^2$, and\/ $u,v\in C^0(S)$
a singular solution of\/ \eq{up4eq2}, as in \S\ref{up46}. If\/
$u(x,-y)\equiv u(x,y)$ and\/ $v(x,-y)\equiv -v(x,y)$ near the
$x$-axis in $S$ then $(u,v)$ is singular along the intersection
of the $x$-axis with\/ $S$, and the singularities are nonisolated.
Otherwise there are at most countably many singular points of\/
$(u,v)$ in $S^\circ$, all isolated.
\label{up8thm1}
\end{thm}

We may divide isolated singularities $(b,0)$ into four types,
depending on the behaviour of $v(x,0)$ near~$(b,0)$.

\begin{dfn} Let $S$ be a domain in $\R^2$, and $u,v\in C^0(S)$ a
singular solution of \eq{up4eq2}, as in \S\ref{up46}. Suppose
$(b,0)$ is an isolated singular point of $(u,v)$ in $S^\circ$.
Then there exists $\ep>0$ such that $\ovB_\ep(b,0)\subset S^\circ$
and $(b,0)$ is the only singularity of $(u,v)$ in $\ovB_\ep(b,0)$.
Thus, for $0<\md{x-b}\le\ep$ we have $(x,0)\in S^\circ$ and
$v(x,0)\ne 0$. So by continuity $v$ is either positive or
negative on each of $[b-\ep,b)\t\{0\}$ and~$(b,b+\ep]\t\{0\}$.
\begin{itemize}
\setlength{\itemsep}{0pt}
\setlength{\parsep}{0pt}
\item[{\rm(i)}] If $v(x)<0$ for $x\in[b-\ep,b)$ and $v(x)>0$ for
$x\in(b,b+\ep]$ we say the singularity $(b,0)$ is of
{\it increasing type}.
\item[{\rm(ii)}] If $v(x)>0$ for $x\in[b-\ep,b)$ and $v(x)<0$ for
$x\in(b,b+\ep]$ we say the singularity $(b,0)$ is of
{\it decreasing type}.
\item[{\rm(iii)}] If $v(x)<0$ for $x\in[b-\ep,b)$ and $v(x)<0$ for
$x\in(b,b+\ep]$ we say the singularity $(b,0)$ is of
{\it maximum type}.
\item[{\rm(iv)}] If $v(x)>0$ for $x\in[b-\ep,b)$ and $v(x)>0$ for
$x\in(b,b+\ep]$ we say the singularity $(b,0)$ is of
{\it minimum type}.
\end{itemize}
\label{up8def2}
\end{dfn}

The type determines whether the multiplicity of $(b,0)$ is even or odd.

\begin{prop} Let\/ $u,v\in C^0(S)$ be a singular solution of\/
\eq{up4eq2} on a domain $S$ in $\R^2$, and\/ $(b,0)$ be an
isolated singularity of\/ $(u,v)$ in $S^\circ$ with multiplicity
$k$. If\/ $(b,0)$ is of increasing or decreasing type then $k$ is
odd, and if\/ $(b,0)$ is of maximum or minimum type then $k$ is even.
\label{up8prop1}
\end{prop}

\begin{proof} Let $\ep$ be as in Definition \ref{up8def2},
and define $u',v'$ on $\ovB_\ep(b,0)$ by $u'(x,y)=u(x,-y)$
and $v'(x,y)=-v(x,-y)$. Then by definition, $k$ is the winding
number of $(u,v)-(u',v')$ about 0 along $\ga_\ep(b,0)$.
Divide the circle $\ga_\ep(b,0)$ into upper and lower
semicircles, with $y\ge 0$ and~$y\le 0$.

By the involution $(x,y)\mapsto(x,-y)$ we see that $(u,v)-(u',v')$
rotates through the same angle on the upper and lower semicircles.
Now $u(b\pm\ep,0)=u'(b\pm\ep,0)$ and $v(b\pm\ep,0)=-v'(b\pm\ep,0)$,
so $(u,v)-(u',v')$ lies on the positive $y$-axis at $(b\pm\ep,0)$
if $v(b\pm\ep,0)>0$, and on the negative $y$-axis if~$v(b\pm\ep,0)<0$.

Using this we see that in cases (i) and (ii) of Definition
\ref{up8def2}, $(u,v)-(u',v')$ rotates through an angle $(2n+1)\pi$
on the upper semicircle, for $n\in\Z$. Thus $(u,v)-(u',v')$ rotates
through $(2n+1)2\pi$ about $\ga_\ep(b,0)$, and $k=2n+1$ is odd.
Similarly, in cases (iii) and (iv), $(u,v)-(u',v')$ rotates through
an angle $2n\pi$ on the upper semicircle and $4n\pi$ about
$\ga_\ep(b,0)$, and $k=2n$ is even.
\end{proof}

Now an isolated singular point of $u,v$ in $S^\circ$ yields
an isolated singular point $\bf z$ in the interior of the
corresponding $\U(1)$-invariant SL 3-fold $N$. The possible
tangent cones $C$ to $N$ at $\bf z$ were classified in Theorem
\ref{up5thm3}. For cases (i) and (ii) we can use Proposition
\ref{up5prop1} and Theorem \ref{up4thm1} to identify the
multiplicity and type of~$(b,0)$.

\begin{prop} Let\/ $u,v\in C^0(S)$ be a singular solution of\/
\eq{up4eq2} on a domain $S$ in $\R^2$, and\/ $N$ the corresponding
$\U(1)$-invariant SL\/ $3$-fold in $\C^3$. Let\/ $(b,0)$ be an
isolated singularity of\/ $(u,v)$ in $S^\circ$, and\/
${\bf z}=\bigl(0,0,b+iu(b,0)\bigr)$ the corresponding singular
point of\/ $N$. Suppose $C$ is a tangent cone to $N$ at\/ $\bf z$.
If\/ $C$ is as in case {\rm(i)} of Theorem \ref{up5thm3} then
$(b,0)$ has multiplicity $1$ and is of increasing type.
If\/ $C$ is as in case {\rm(ii)} of Theorem \ref{up5thm3} then
$(b,0)$ has multiplicity $1$ and is of decreasing type.
\label{up8prop2}
\end{prop}

The author can also show that if $(b,0)$ has multiplicity
$n\ge 2$ then the unique tangent cone $C$ to $N$ at $\bf z$ is
$\Pi^\phi_+\cup\Pi^\phi_-$ for some $\phi\in(-\frac{\pi}{2},
\frac{\pi}{2}]$, where $\Pi^\phi_\pm$ are defined in
\eq{up5eq5} and have multiplicity~1.

Next we give counting formulae for singularities. Theorem
\ref{up6thm5} yields:

\begin{thm} Let\/ $S$ be a domain in $\R^2$ invariant under
$(x,y)\mapsto(x,-y)$, and\/ $u,v\in C^0(S)$ a singular solution
of\/ \eq{up4eq2}, as in \S\ref{up46}. Define $u',v'\in C^0(S)$
by $u'(x,y)=u(x,-y)$ and\/ $v'(x,y)=-v(x,-y)$. Suppose
$(u,v)\neq(u',v')$ at every point of\/ $\pd S$. Then $(u,v)$ has
finitely many singularities in $S$, all isolated. Let there be
$n$ singularities with multiplicities $k_1,\ldots,k_n$. Then
the winding number of\/ $(u,v)-(u',v')$ about\/ $0$ along
$\pd S$ is at least\/~$\sum_{i=1}^nk_i$.
\label{up8thm2}
\end{thm}

Here in the last line we say at least rather than exactly
$\sum_{i=1}^nk_i$, as $(u,v)-(u',v')$ may have other zeroes not
on the $x$-axis, and so not singular points, which contribute
to the winding number. Similarly, Theorem \ref{up6thm7} gives:

\begin{thm} Suppose $S$ is a strictly convex domain in $\R^2$
invariant under $(x,y)\mapsto(x,-y)$, and\/ $\phi\in C^{3,\al}
(\pd S)$ for some $\al\in(0,1)$. Let\/ $u,v\in C^0(S)$ be the
singular solution of\/ \eq{up4eq2} constructed in Theorem
\ref{up4thm4} from~$\phi$.

Define $\phi'\in C^{3,\al}(\pd S)$ by $\phi'(x,y)=-\phi(x,-y)$.
Suppose $\phi-\phi'$ has exactly $l$ local maxima and\/ $l$ local
minima on $\pd S$. Then $(u,v)$ has finitely many singularities
in $S^\circ$, all isolated. Let there be $n$ singularities in
$S^\circ$ with multiplicities $k_1,\ldots,k_n$.
Then~$\sum_{i=1}^nk_i\le l-1$.
\label{up8thm3}
\end{thm}

\section{Singularities exist with all multiplicities}
\label{up9}

We now prove that there exist singularities with every
multiplicity $n\ge 1$, and every possible type, and that
singularities of multiplicity $n$ occur in codimension $n$
in the family of all $\U(1)$-invariant SL 3-folds, in a
certain sense. For simplicity we work not on a general
domain $S$, but on the unit disc $D$ in~$\R^2$.

\begin{dfn} Let $D$ be the unit disc $\bigl\{(x,y)\in\R^2:x^2+
y^2\le 1\bigr\}$ in $\R^2$, with boundary ${\cal S}^1$, the
unit circle. Define a coordinate $\th:\R/2\pi\Z\ra{\cal S}^1$
by $\th\mapsto(\cos\th,\sin\th)$. Then $\cos(j\th),\sin(j\th)\in
C^\iy({\cal S}^1)$ for~$j\ge 1$.
\label{up9def1}
\end{dfn}

We shall use the functions $\cos(j\th),\sin(j\th)$ for $j=1,\ldots,n$
as a family of perturbations to the boundary data $\phi$ in Theorem
\ref{up4thm4}, and show that for any suitable $\phi$ there exists a
family of perturbations of $\phi$ such that the corresponding
singular solution $u,v$ has a singularity of multiplicity at
least $n$ at~$(0,0)$.

\subsection{Counting stationary points of Fourier sums}
\label{up91}

We begin with two propositions on the stationary points of
linear combinations of $\cos(j\th)$, $\sin(j\th)$. We say that
$\phi$ has a {\it stationary point of multiplicity} $k\ge 1$ at
$\th_0$ if $\d^j\phi/\d\th^j(\th_0)=0$ for $j=1,\ldots,k$,
but~$\d^{k+1}\phi/\d\th^{k+1}(\th_0)\ne 0$.

\begin{prop} Let\/ $\al_1,\be_1,\ldots,\al_n,\be_n\in\R$ be not
all zero for $n\ge 1$, and define $\phi\in C^\iy({\cal S}^1)$ by
$\phi=\sum_{j=1}^n\bigl(\al_j\cos(j\th)+\be_j\sin(j\th)\bigr)$.
Then $\phi$ has at most\/ $2n$ stationary points in ${\cal S}^1$,
counted with multiplicity.
\label{up9prop1}
\end{prop}

\begin{proof} Suppose $\phi$ has more than $2n$ stationary points,
and choose $2n+1$ of them. These divide ${\cal S}^1$ into $2n+1$
intervals. For each interval $\frac{\d\phi}{\d\th}=0$ at the end
points, so $\d^2\phi/\d\th^2=0$ somewhere in the open interval
by Rolle's Theorem. Thus there are at least $2n+1$ points in
${\cal S}^1$ with $\d^2\phi/d\th^2=0$. Similarly, by induction
on $k$ there are at least $2n+1$ points in ${\cal S}^1$
with $\d^k\phi/\d\th^k=0$, for all~$k\ge 1$.

But when $k$ is large $\d^k\phi/\d\th^k$ is dominated by the terms
in $\cos(j\th),\sin(j\th)$ for the biggest $j=1,\ldots,n$ for
which $\al_j\ne 0$. So
\begin{equation*}
\frac{\d^{2k}\phi}{\d\th^{2k}}\approx(-1)^kj^{2k}\bigl(\al_j
\cos(j\th)+\be_j\sin(j\th)\bigr)\quad\text{for large $k$.}
\end{equation*}
Thus $\d^{2k}\phi/\d\th^{2k}$ has at most $2j\le 2n$ zeroes for
large $k$, a contradiction.

This shows that $\phi$ has at most $2n$ stationary points, counted
{\it without\/} multiplicity. To prove the proposition counting
{\it with\/} multiplicity, we have to modify the proof above
slightly. If $\th_0$ is a stationary point of multiplicity $k>1$ then
$\d^j\phi/\d\th^j$ has an extra zero at $\th_0$ for $j=2,\ldots,k$,
which was not taken into account above. Including these, we see that
$\th_0$ contributes $k$ zeroes to $\d^l\phi/\d\th^l$ for all $l\ge k$,
and the result follows.
\end{proof}

\begin{prop} For each\/ $n\ge 1$ there exists $K_n>0$ such that if\/
$\psi\in C^{2n}({\cal S}^1)$ and\/ $\al_1,\be_1,\ldots,\al_n,\be_n\in\R$
with\/ $\sum_{j=1}^n(\al_j^2+\be_j^2)>K_n\cnm{\psi}{2n}^2$, then $\phi=
\psi+\sum_{j=1}^n\bigl(\al_j\cos(j\th)+\be_j\sin(j\th)\bigr)$ has at
most\/ $2n$ stationary points in ${\cal S}^1$, counted with multiplicity.
\label{up9prop2}
\end{prop}

\begin{proof} The idea is that if we perturb $\phi'=\sum_{j=1}^n\bigl(
\al_j\cos(j\th)+\be_j\sin(j\th)\bigr)$ in Proposition \ref{up9prop1}
by a {\it sufficiently small} perturbation $\psi$ then the number of
stationary points, counted with multiplicity, does not increase. The
problem is to decide what `sufficiently small' should mean here.

If $\psi$ is small in $C^1$ then stationary points of $\phi=
\psi+\phi'$ can only appear near stationary points of $\phi'$.
Let $\th_0$ be a stationary point of $\phi'$ with multiplicity
$k$. If $\psi$ is small in $C^{k+1}$ then $\d^k\phi/\d\th^k$ is
nonzero near $\th_0$, and so $\phi$ can have at most $k$
stationary points near $\th_0$, counted with multiplicity.

Now $\phi'$ has at most $2n$ stationary points counted with
multiplicity, by Proposition \ref{up9prop1}, and it has at
least 2. Thus the maximum multiplicity of a stationary point
of $\phi'$ is $2n-1$. Hence if $\psi$ is small in $C^{2n}$
compared to $\phi'$, then the number of stationary points of
$\phi$ does not increase near every stationary point of $\phi'$. 
Using this the proposition easily follows.
\end{proof}

\subsection{Zeroes of multiplicity $n$ in a family of perturbations}
\label{up92}

We now prove a general result on the existence of zeroes of
$(u,v)-(\hat u,\hat v)$ of multiplicity at least $n$ in a family
of perturbations of $(u,v)$. Here and in the rest of the section
we use the convention of Definition~\ref{up6def3}.

\begin{thm} Let\/ $a\ne 0$, $n\ge 1$, $\al_1,\be_1,\ldots,\al_n,
\be_n\in\R$ and\/ $\psi,\hat\phi\in C^\iy({\cal S}^1)$. Define
$\phi=\psi+\sum_{j=1}^n\bigl(\al_j\cos(j\th)+\be_j\sin(j\th)\bigr)$
in $C^\iy({\cal S}^1)$. Let\/ $f,\hat f\in C^\iy(D)$ be the unique
solutions of\/ \eq{up4eq7} on $D$ with\/ $f\vert_{{\cal S}^1}=\phi$,
$\hat f\vert_{{\cal S}^1}=\hat\phi$, which exist by Theorem
\ref{up4thm3}. Let\/ $u=\frac{\pd f}{\pd y}$, $v=\frac{\pd f}{\pd x}$,
$\hat u=\frac{\smash{\pd\hat f}}{\pd y}$ and\/~$\hat v=
\frac{\smash{\pd\hat f}}{\pd x}$.

Then for each\/ $a\ne 0$, $n\ge 1$ and\/ $\psi,\hat\phi\!\in\!
C^\iy({\cal S}^1)$ there exist unique $\al_1,\be_1,\ldots,
\allowbreak
\al_n,\be_n$
with\/ $\sum_{j=1}^n(\al_j^2+\be_j^2)\le K_n\cnm{\psi-\hat\phi}{2n}^2$,
where $K_n$ is as in Proposition \ref{up9prop2}, such that
$(u,v)-(\hat u,\hat v)$ has a zero of multiplicity at least\/ $n$ at\/
$(0,0)$. Furthermore these $\al_j,\be_j$ depend continuously on~$a,\psi$.
\label{up9thm1}
\end{thm}

\begin{proof} We shall prove the theorem by induction on $n$. The
first and inductive steps will be shown together. Let $k\ge 0$ and
suppose the theorem holds for all $n\le k$. We will prove it when
$n=k+1$. Fix $a\ne 0$ and~$\psi,\hat\phi\in C^\iy({\cal S}^1)$.

Let $\ga=\al_{k+1}+i\be_{k+1}\in\C$. When $k=0$, let $f_\ga\in C^\iy(D)$
be the unique solution of \eq{up4eq7} with $f_\ga\vert_{{\cal S}^1}=
\psi+\al_1\cos\th+\be_1\sin\th$. When $k\ge 1$, apply the inductive
hypothesis for $n=k$ with $\psi$ replaced by $\psi+\al_{k+1}\cos(k\!+\!1)
\th+\be_{k+1}\sin(k\!+\!1)\th$. This gives unique $\al_1,\be_1,\ldots,
\al_k,\be_k\in\R$ and a solution $f_\ga\in C^\iy(D)$ of \eq{up4eq7} with
\e
f_\ga\vert_{{\cal S}^1}=\psi+\sum_{j=1}^{k+1}
\bigl(\al_j\cos(j\th)+\be_j\sin(j\th)\bigr).
\label{up9eq1}
\e
Define $u_\ga,v_\ga$ in the obvious way.

When $k\ge 1$ we know by induction that $(u_\ga,v_\ga)-(\hat u,\hat v)$
has a zero of multiplicity at least $k$ at $(0,0)$. Therefore by
Proposition \ref{up4prop2} we may write
\e
\la\,u_\ga(x,y)\!+\!iv_\ga(x,y)=\la\,\hat u(x,y)\!+\!
i\hat v(x,y)\!+\!C\bigl(\la x\!+\!iy\bigr)^k
\!+\!O\bigl(\md{x}^{k+1}\!+\!\md{y}^{k+1}\bigr)
\label{up9eq2}
\e
near $(0,0)$, where $\la>0$, and $C\in\C$ is zero if and only if
$(u_\ga,v_\ga)-(\hat u,\hat v)$ has a zero of multiplicity more
than $k$. Define $F(\ga)=C$, giving a map $F:\C\ra\C$. When $k=0$
this formula is still valid, giving $F(\ga)=\bigl(\la(u_\ga-\hat u)
+i(v_\ga-\hat v)\bigr)\vert_{(0,0)}$. Here are some properties of
the maps~$F$.

\begin{prop} This map $F:\C\ra\C$ is continuous and injective,
and if\/ $\ms{\ga}>K_{k+1}\cnm{\psi-\hat\phi}{2k+2}^2$ then
$F(\ga)\ne 0$, where $K_{k+1}$ is as in Proposition~\ref{up9prop2}.
\label{up9prop3}
\end{prop}

\begin{proof} By the inductive hypothesis $\al_1,\be_1,\ldots,\al_k,
\be_k$ depend continuously on $a,\psi,\hat\phi$ and $\ga$. Hence
$f_\ga,u_\ga,v_\ga$ depend continuously in $C^\iy(D)$ on $a,\psi,\ga$
by \cite[Th.~7.7]{Joyc4}. It follows that $F$ is continuous as a
function of $\ga$, as we have to prove, and also varies continuously
with~$a,\psi,\hat\phi$.

Next we show $F$ is injective. Suppose for a contradiction that
$\ga\ne\ga'\in\C$ and $F(\ga)=F(\ga')$. Then \eq{up9eq2} holds for
$u_\ga,v_\ga$ and $u_{\ga'},v_{\ga'}$ with the same value of $C$.
This implies that $v_\ga(0,0)=v_{\ga'}(0,0)$, so the values of
$\la=\sqrt{2}(v(0,0)^2+a^2)^{1/4}$ for $\ga$ and $\ga'$ are the
same as well. Thus, subtracting \eq{up9eq2} for $\ga,\ga'$ gives
\begin{equation*}
\la\,u_\ga(x,y)+iv_\ga(x,y)=\la\,u_{\ga'}(x,y)+iv_{\ga'}(x,y)
+O\bigl(\md{x}^{k+1}\!+\!\md{y}^{k+1}\bigr).
\end{equation*}
Hence $(u_\ga,v_\ga)-(u_{\ga'},v_{\ga'})$ has a zero of multiplicity
at least $k+1$ at $(0,0)$, by Definition~\ref{up4def2}.

Now by \eq{up9eq1} the potentials $f_\ga,f_{\ga'}$ satisfy
\begin{equation*}
(f_\ga-f_{\ga'})\vert_{{\cal S}^1}=\sum_{j=1}^{k+1}
\bigl((\al_j-\al_j')\cos(j\th)+(\be_j-\be_j')\sin(j\th)\bigr).
\end{equation*}
Thus $(f_\ga-f_{\ga'})\vert_{{\cal S}^1}$ has at most $2k+2$
stationary points by Proposition \ref{up9prop1}, and so at
most $k+1$ local maxima and $k+1$ local minima. Therefore by
Theorem \ref{up4thm6} the number of zeroes of $(u_\ga,v_\ga)-
(u_{\ga'},v_{\ga'})$ in $D^\circ$, counted with multiplicity,
is at most $k$. But this contradicts $(0,0)$ being a zero of
multiplicity at least $k+1$, so $F$ is injective.

Finally, if $\ms{\ga}>K_{k+1}\cnm{\psi-\hat\phi}{2k+2}^2$ then
$(f_\ga-\hat f)\vert_{{\cal S}^1}$ has at most $2k+2$ stationary
points by Proposition \ref{up9prop2}, so there are at most $k$
zeroes of $(u_\ga,v_\ga)-(\hat u,\hat v)$ in $D^\circ$ with
multiplicity by Theorem \ref{up4thm6}. But $C=F(\ga)=0$ if and
only if $(0,0)$ has multiplicity more than $k$, giving
$F(\ga)\ne 0$, as we have to prove.
\end{proof}

We shall show that $F$ has a zero in $\C$. Consider the case
$\psi=\hat\phi=0$. Then when $\ga=0$ we have $u_0=v_0=\hat u=
\hat v=0$, so that $F(0)=0$ by \eq{up9eq2}. As $F$ is injective
and continuous by Proposition \ref{up9prop3}, it follows that
the winding number of $F$ about 0 along the circle $\ga_r(0)$
in $\C$ with radius $r>0$ and centre 0 is~$\pm 1$.

We may deform any $\psi,\hat\phi\in C^\iy({\cal S}^1)$ continuously
to $0,0$ through $t\psi,t\hat\phi$ for $t\in[0,1]$. This gives a
1-parameter family of functions $F_t:\C\ra\C$, with $F_0(0)=0$,
and we seek a zero $\ga$ of $F_1=F$. Now $F_t$ depends continuously
on $t$, as $F$ depends continuously on $\psi,\hat\phi$, and
$F_t(\ga)\ne 0$ if $\ms{\ga}>K_{k+1}\cnm{\psi-\hat\phi}{2k+2}^2$
by Proposition~\ref{up9prop3}.

Thus, if $r^2>K_{k+1}\cnm{\psi-\hat\phi}{2k+2}^2$ then $F_t$ is
nonzero on $\ga_r(0)$ for all $t\in[0,1]$, and so the winding
number $F_1$ about 0 along $\ga_r(0)$ is the same as the winding
number of $F_0$, which is $\pm 1$. Therefore $F=F_1$ must have a
zero inside $\ga_r(0)$, as otherwise the winding number would be~0.

We have shown that the function $F$ has a zero $\ga=\al_{k+1}
+i\be_{k+1}\in\C$. This $\ga$ is unique by injectivity in
Proposition \ref{up9prop3}. The construction above then yields
$\al_1,\be_1,\ldots,\al_{k+1},\be_{k+1}$ such that $(u,v)-
(\hat u,\hat v)=(u_\ga,v_\ga)-(\hat u,\hat v)$ has a zero of
multiplicity at least $k+1$, as we have to prove. The $\al_j,\be_j$
are unique for $j=k+1$ by uniqueness of $\ga$, and for $1\le j\le k$
by the inductive hypothesis.

Now $\ga=\al_{k+1}+i\be_{k+1}$ depends continuously on $a,\psi,
\hat\phi$ as $F$ does, and $\al_1,\be_1,\ldots,\al_k,\be_k$ depend
continuously on $a,\psi,\hat\phi$ and $\ga$ by the inductive
hypothesis. Thus $\al_1,\be_1,\ldots,\al_{k+1},\be_{k+1}$ depend
continuously on~$a,\psi,\hat\phi$.

If $\sum_{j=1}^{k+1}(\al_j^2+\be_j^2)>K_{k+1}\cnm{\psi-
\hat\phi}{2k+2}^2$ then $\phi-\hat\phi$ has at most $2k+2$
stationary points on ${\cal S}^1$ by Proposition \ref{up9prop2},
so there are at most $k$ zeroes of $(u,v)-(\hat u,\hat v)$ in
$D^\circ$ with multiplicity by Theorem \ref{up4thm6}. But $(0,0)$
is a zero with multiplicity at least $k+1$, a contradiction.
Hence $\sum_{j=1}^{k+1}(\al_j^2+\be_j^2)\le K_{k+1}\cnm{\psi-
\hat\phi}{2k+2}^2$. This completes the induction, and the proof
of Theorem~\ref{up9thm1}.
\end{proof}

The following lemma is easily proved using Proposition \ref{up4prop2},
by subtracting \eq{up4eq6} for $u_1,v_1,u_2,v_2$ from \eq{up4eq6}
for~$u_1,v_1,u_3,v_3$.

\begin{lem} Let\/ $S$ be a domain in $\R^2$, let\/ $u_j,v_j\in C^1(S)$
be solutions of\/ \eq{up4eq3} for $j=1,2,3$ and some $a\ne 0$, and
let\/ $(b,c)\in S^\circ$. Suppose $(u_1,v_1)-(u_2,v_2)$ has a zero of
multiplicity at least\/ $k$ and\/ $(u_1,v_1)-(u_3,v_3)$ a zero of
multiplicity at least\/ $l$ at\/ $(b,c)$. Then $(u_2,v_2)-(u_3,v_3)$
has a zero of multiplicity at least\/ $\min(k,l)$ at\/~$(b,c)$.
\label{up9lem}
\end{lem}

Combining the last two results, we prove:

\begin{thm} Let\/ $a\ne 0$, $n\ge 1$, $\ga_1,\ldots,\ga_n\in\R$
and\/ $\psi\in C^\iy({\cal S}^1)$. Then there exists $K>0$ depending
only on $n,\ga_j$ and\/ $\psi$ and unique $\al_1,\be_1,\ldots,\al_n,
\be_n$ with\/ $\sum_{j=1}^n(\al_j^2+\be_j^2)\le K$ such that the
following holds. Furthermore, for fixed\/ $a,n,\psi$ the map
$(\ga_1,\ldots,\ga_n)\mapsto(\al_1,\be_1,\ldots,\al_n,\be_n)$
is injective. Define
\begin{equation*}
\phi=\psi+\sum_{j=1}^n\bigl(\al_j\cos(j\th)+\be_j\sin(j\th)\bigr)
\quad\text{and\/}\quad\hat\phi=\sum_{j=1}^n\ga_j\sin(j\th)
\quad\text{in $C^\iy({\cal S}^1)$.}
\end{equation*}
Let\/ $f,\hat f\in C^\iy(D)$ be the unique solutions of\/ \eq{up4eq7}
on $D$ with\/ $f\vert_{{\cal S}^1}=\phi$, $\hat f\vert_{{\cal S}^1}=
\hat\phi$, which exist by Theorem \ref{up4thm3}. Define $u=\frac{\pd
f}{\pd y}$, $v=\frac{\pd f}{\pd x}$, $\hat u=\frac{\smash{\pd\hat
f}}{\pd y}$, $\hat v=\frac{\smash{\pd\hat f}}{\pd x}$, $u'(x,y)=
u(x,-y)$ and\/ $v'(x,y)=-v(x,-y)$, so that\/ $(u,v)$, $(\hat u,\hat
v)$ and\/ $(u',v')$ satisfy \eq{up4eq3} in $C^\iy(D)$. Then $(u,v)-
(\hat u,\hat v)$ and\/ $(u,v)-(u',v')$ both have a zero of
multiplicity at least\/ $n$ at\/~$(0,0)$.
\label{up9thm2}
\end{thm}

\begin{proof} Most of the theorem is immediate from Theorem
\ref{up9thm1}. We need only prove that $(u,v)-(u',v')$ has
a zero of multiplicity at least $n$ at $(0,0)$, and that
$\ga_j\mapsto\al_j,\be_j$ is injective. As $\hat\phi(x,-y)
\equiv -\hat\phi(x,y)$ we have $\hat f(x,-y)\equiv -\hat
f(x,y)$, $\hat u(x,-y)=\hat u(x,y)$ and $\hat v(x,-y)=
-\hat v(x,y)$. Therefore, as $(u,v)-(\hat u,\hat v)$ has a
zero of multiplicity at least $n$ at $(0,0)$, we see by
applying the symmetry $(x,y)\mapsto(x,-y)$ that $(u',v')-
(\hat u,\hat v)$ has a zero of multiplicity at least $n$
at $(0,0)$. Lemma \ref{up9lem} then shows that $(u,v)-
(u',v')$ has a zero of multiplicity at least $n$ at~$(0,0)$.

Now suppose that for fixed $a,n,\psi$, two distinct $n$-tuples
$\ga_1,\ldots,\ga_n$ and $\ti\ga_1,\ldots,\ti\ga_n$ yield the
same $\al_1,\be_1,\ldots,\al_n,\be_n$. Let $(\hat u,\hat v)$ and
$(\ti v,\ti v)$ be the solutions of \eq{up4eq3} corresponding to
the $\ga_j$ and $\ti\ga_j$, with potentials $\hat f,\ti f$, and
$(u,v)$ the solution corresponding to the~$\al_j,\be_j$.

From above, $(u,v)-(\hat u,\hat v)$ and $(u,v)-(\ti u,\ti v)$
both have a zero of multiplicity at least $n$ at $(0,0)$. So
$(\hat u,\hat v)-(\ti u,\ti v)$ has a zero of multiplicity at
least $n$ at $(0,0)$ by Lemma \ref{up9lem}. But the corresponding
potentials $\hat f,\ti f$ satisfy
\begin{equation*}
(\hat f-\ti f)\vert_{{\cal S}^1}=\ts\sum_{j=1}^n(\ga_j-\ti\ga_j)
\sin(j\th).
\end{equation*}
So $(\hat f-\ti f)\vert_{{\cal S}^1}$ has at most $2n$ stationary
points in ${\cal S}^1$ by Proposition \ref{up9prop1}, as $\ga_j
\not\equiv\ti\ga_j$. Hence $(\hat u,\hat v)-(\ti u,\ti v)$ has at
most $n-1$ zeroes in $D^\circ$ with multiplicity, a contradiction.
So the map $(\ga_1,\ldots,\ga_n)\mapsto(\al_1,\be_1,\ldots,\al_n,
\be_n)$ is injective.
\end{proof}

\subsection{Constructing singularities with multiplicity $n$}
\label{up93}

By taking the limit $a\ra 0$ in Theorem \ref{up9thm2} we prove:

\begin{thm} Let\/ $n\ge 1$, $\ga_1,\ldots,\ga_n\in\R$ and\/
$\psi\in C^\iy({\cal S}^1)$. Then there exist unique $\al_1,\be_1,
\ldots,\al_n,\be_n\in\R$ such that the following holds. Furthermore,
for fixed\/ $a,\psi$ the map $(\ga_1,\ldots,\ga_n)\mapsto(\al_1,\be_1,
\ldots,\al_n,\be_n)$ is injective. Define
\begin{equation*}
\phi=\psi+\sum_{j=1}^n\bigl(\al_j\cos(j\th)+\be_j\sin(j\th)\bigr)
\quad\text{and\/}\quad\hat\phi=\sum_{j=1}^n\ga_j\sin(j\th)
\quad\text{in $C^\iy({\cal S}^1)$.}
\end{equation*}
Let\/ $f,\hat f\in C^1(D)$ be the unique
weak solutions of\/ \eq{up4eq8} on $D$ with\/ $f\vert_{{\cal S}^1}=\phi$,
$\hat f\vert_{{\cal S}^1}=\hat\phi$ and\/ $a=0$, which exist by Theorem
\ref{up4thm4}. Define $u=\frac{\pd f}{\pd y}$, $v=\frac{\pd f}{\pd x}$,
$\hat u=\frac{\smash{\pd\hat f}}{\pd y}$ and\/ $\hat v=\frac{\smash{\pd
\hat f}}{\pd x}$, so that\/ $(u,v)$ and\/ $(\hat u,\hat v)$ are singular
solutions of\/ \eq{up4eq2} in $C^0(D)$. Then $(u,v)-(\hat u,\hat v)$ has
a zero of multiplicity at least\/ $n$ at\/ $(0,0)$, and either $(u,v)$
has an isolated singularity of multiplicity at least\/ $n$ at\/
$(0,0)$, or $u(x,-y)\equiv u(x,y)$ and\/~$v(x,-y)\equiv -v(x,y)$.
\label{up9thm3}
\end{thm}

\begin{proof} Fix $n,\ga_j$ and $\psi$ as above. For each $a\ne 0$,
Theorem \ref{up9thm2} gives unique $\al_1,\be_1,\ldots,\al_n,\be_n$
with $\sum_{j=1}^n(\al_j^2+\be_j^2)\le K$ satisfying certain
conditions. As the set of $\al_j,\be_j$ with $\sum_{j=1}^n
(\al_j^2+\be_j^2)\le K$ is compact, we can choose a sequence
$(a^i)_{i=1}^\iy$ in $(0,1]$ such that $a^i\ra 0$ as $i\ra\iy$,
such that the corresponding $\al_j^i,\be_j^i$ converge to limits
$\al_j,\be_j$ in~$\R^{2n}$.

Let $f^i,u^i,v^i,\hat f^i,\hat u^i,\hat v^i\in C^\iy(D)$ be the
corresponding solutions $f,u,v,\hat f,\hat u,\hat v$ in Theorem
\ref{up9thm2} for $i=1,2,\ldots$, so that $f^i,a^i$ and $\hat f^i,
a^i$ satisfy \eq{up4eq7}, and $u^i,v^i,a^i$ and $\hat u^i,\hat v^i,
a^i$ satisfy \eq{up4eq3} in $D$. Define $\phi=\psi+\sum_{j=1}^n\bigl
(\al_j\cos(j\th)+\be_j\sin(j\th)\bigr)$ in $C^\iy({\cal S}^1)$, and
let $f,\hat f\in C^1(D)$ be the unique weak solutions of \eq{up4eq8}
on $D$ with $f\vert_{{\cal S}^1}=\phi$, $\hat f\vert_{{\cal S}^1}=
\hat\phi$, which exist by Theorem \ref{up4thm4}. Let $u,v$ and
$\hat u,\hat v$ be the corresponding singular solutions
of~\eq{up4eq2}.

Then $\hat f^i\vert_{{\cal S}^1}=\hat f\vert_{{\cal S}^1}=\hat\phi$
for all $i\ge 1$. Hence by Theorem \ref{up4thm5} we see that
$\hat f^i\ra\hat f$ in $C^1(D)$ as $i\ra\iy$, so that
$\hat u^i,\hat v^i\ra\hat u,\hat v$ in $C^0(D)$ as $i\ra\iy$.
Similarly, as $\al_j^i,\be_j^i\ra\al_j,\be_j$ as $i\ra\iy$,
we see that $f^i\vert_{{\cal S}^1}\ra f\vert_{{\cal S}^1}=\phi$ in
$C^\iy({\cal S}^1)$ as $i\ra\iy$. Hence by Theorem \ref{up4thm5}
we see that $f^i\ra f$ in $C^1(D)$ and $u^i,v^i\ra u,v$ in $C^0(D)$
as~$i\ra\iy$.

Thus $(u^i,v^i)-(\hat u^i,\hat v^i)$ has a zero of multiplicity
at least $n$ at $(0,0)$, and $u^i,v^i,\hat u^i,\hat v^i$ converge
to $u,v,\hat u,\hat v$ in $C^0(S)$ as $i\ra\iy$. So $(u,v)-(\hat
u,\hat v)$ has a zero at $(0,0)$. By Theorem \ref{up6thm7}, either
$(u,v)\equiv(\hat u,\hat v)$, or $(0,0)$ is an isolated zero. If
it is isolated, then by $C^0$ convergence we find that $(u^i,v^i)
-(\hat u^i,\hat v^i)$ and $(u,v)-(\hat u,\hat v)$ have the same
winding number about 0 along $\ga_\ep(0,0)$ for small $\ep>0$ and
large $i$. Therefore $(u,v)-(\hat u,\hat v)$ has multiplicity at
least $n$ at~$(0,0)$.

In the same way, if $u'(x,y)=u(x,-y)$ and $v'(x,y)=-v(x,-y)$ we
find that $(u,v)-(u',v')$ has a zero of multiplicity at least
$n$ at $(0,0)$, and so by definition either $(u,v)\equiv(u',v')$
or $(u,v)$ has an isolated singularity of multiplicity at least
$n$ at $(0,0)$. It remains only to show that $\al_1,\be_1,\ldots,
\al_n,\be_n$ are {\it unique}, and the map $\ga_j\ra\al_j,\be_j$
{\it injective}.

Suppose $\al_j,\be_j,\phi,f,u,v$ and $\ti\al_j,\ti\be_j,\ti\phi,
\ti f,\ti u,\ti v$ are two distinct solutions. Then
\begin{equation*}
\phi-\ti\phi=\sum_{j=1}^n\bigl((\al_j-\ti\al_j)\cos(j\th)+
(\be_j-\ti\be_j)\sin(j\th)\bigr),
\end{equation*}
so as $\al_j-\ti\al_j,\be_j-\ti\be_j$ are not all zero
$\phi-\ti\phi$ has at most $2n$ stationary points by
Proposition \ref{up9prop1}. Theorem \ref{up6thm7} then shows
that there are at most $n-1$ zeroes of $(u,v)-(\ti u,\ti v)$
in $D^\circ$, counted with multiplicity.

But as $(u,v)-(\hat u,\hat v)$ and $(\ti u,\ti v)-(\hat u,\hat v)$
both have a zero of multiplicity at least $n$ at $(0,0)$, by a
version of Lemma \ref{up9lem} for $a=0$ we see that $(u,v)-
(\ti u,\ti v)$ has a zero of multiplicity at least $n$ at $(0,0)$, a
contradiction. Thus the $\al_j,\be_j$ are unique. A similar proof,
extending that in Theorem \ref{up9thm2}, shows that the map
$(\ga_1,\ldots,\ga_n)\mapsto(\al_1,\be_1,\ldots,\al_n,\be_n)$
is injective.
\end{proof}

For particular $\psi$ we can pin down the multiplicity exactly.

\begin{thm} Take $\psi=\al_{n+1}\cos(n\!+\!1)\th+\be_{n+1}
\sin(n\!+\!1)\th\in C^\iy({\cal S}^1)$, for $\al_{n+1}\ne 0$
and\/ $\be_{n+1}\in\R$. Then for all\/ $\ga_1,\ldots,\ga_n\in\R$
the singular solution $(u,v)$ of\/ \eq{up4eq2} constructed in
Theorem \ref{up9thm3} from $\psi$ and\/ $\ga_1,\ldots,\ga_n$
has an isolated singularity at\/ $(0,0)$ of multiplicity $n$,
and no other singularities in~$D^\circ$.
\label{up9thm4}
\end{thm}

\begin{proof} Let $n,\al_j,\be_j,\phi,f,u,v$ be as in Theorem
\ref{up9thm3}, and set $\phi'(x,y)\!=\!-\phi(x,-y)$, $f'(x,y)=
-f(x,-y)$, $u'(x,y)=u(x,-y)$ and $v'(x,y)=-v(x,-y)$. Then
$\phi-\phi'=2\sum_{j=1}^{n+1}\al_j\cos(j\th)$. As $\al_{n+1}\ne 0$,
Proposition \ref{up9prop1} shows that $\phi-\phi'$ has at most
$2n+2$ stationary points in ${\cal S}^1$, and Theorem \ref{up6thm7}
shows there are at most $n$ zeroes of $(u,v)-(u',v')$ in $D^\circ$,
counted with multiplicity. But $(0,0)$ is a zero of multiplicity at
least $n$. So it has multiplicity exactly $n$, and there are no
other zeroes in $D^\circ$. The result follows.
\end{proof}

Thus there exist singular solutions $(u,v)$ of \eq{up4eq2} with
isolated singularities of all multiplicities $n\ge 1$ at $(0,0)$.
Now by Proposition \ref{up8prop1}, a singularity of multiplicity
$n$ has one of two types. Clearly, if $(u,v)$ has one type, then
$(-u,-v)$ is also singular of multiplicity $n$ at $(0,0)$, but
with the other type. So we prove:

\begin{cor} There exist examples of singular solutions $u,v$
of\/ \eq{up4eq2} with isolated singularities of every possible
multiplicity $n\ge 1$, and with both possible types allowed by
Proposition~\ref{up8prop1}.
\label{up9cor}
\end{cor}

Combining this with Proposition \ref{up4prop1} gives examples of
singular SL 3-folds in $\C^3$ with an infinite number of different
geometrical/topological types. This is one of the main results of
the paper. Also, by combining Example \ref{up7ex} with Corollary
\ref{up9cor} we can construct SL fibrations including singular
fibres with every possible multiplicity and type.

If $\psi(x,-y)\equiv\psi(x,y)$ in Theorem \ref{up9thm3} we can
set $\be_j,\ga_j,\hat\phi,\hat f,\hat u,\hat v$ to zero, and get
a simpler theorem constructing solutions $u,v$ with the symmetries
$u(x,-y)\equiv -u(x,y)$ and~$v(x,-y)\equiv v(x,y)$.

\begin{thm} Let\/ $\psi\in C^\iy({\cal S}^1)$ with\/ $\psi(x,-y)
\equiv\psi(x,y)$ and\/ $n\ge 1$. Then there exist unique $\al_1,
\ldots,\al_n\in\R$ such that the following holds. Define
$\phi=\psi+\sum_{j=1}^n\al_j\cos(j\th)$ in $C^\iy({\cal S}^1)$.
Let\/ $f\in C^1(D)$ be the unique weak solution of\/ \eq{up4eq8}
on $D$ with\/ $f\vert_{{\cal S}^1}=\phi$ and\/ $a=0$, which exists
by Theorem \ref{up4thm4}. Define $u=\frac{\pd f}{\pd y}$ and\/
$v=\frac{\pd f}{\pd x}$, so that\/ $(u,v)$ is a singular solution
of\/ \eq{up4eq2} in $C^0(D)$ with\/ $u(x,-y)\equiv -u(x,y)$ and\/
$v(x,-y)\equiv v(x,y)$. Then either $(u,v)$ has an isolated
singularity of multiplicity at least\/ $n$ at\/ $(0,0)$,
or~$(u,v)\equiv(0,0)$.
\label{up9thm5}
\end{thm}

\subsection{Discussion}
\label{up94}

In \S\ref{up8}--\S\ref{up9} we have put together a detailed
picture of the singularities of singular solutions $u,v$ of
\eq{up4eq2}. Here are some remarks on this, beginning with
how to interpret Theorem~\ref{up9thm3}.

First note that the singular solution $(\hat u,\hat v)$ in Theorem
\ref{up9thm3} has the symmetries $\hat u(x,-y)=\hat u(x,y)$ and
$\hat v(x,-y)=-\hat v(x,y)$, and so was studied in \S\ref{up62}.
In particular, $\hat u(x,0)$ is a real analytic function of $x\in(-1,1)$.
One can show that $\ga_1,\ldots,\ga_n$ are determined uniquely by
$\frac{\pd^m}{\pd x^m}\hat u(0,0)$ for $m=0,\ldots,n-1$, and vice
versa. It seems likely that $\frac{\pd^m}{\pd x^m}u(0,0)$ exists for
$m=0,\ldots,n-1$, and equals~$\frac{\pd^m}{\pd x^m}\hat u(0,0)$.

Thus, singularities with multiplicity $n$ at $(0,0)$ are locally
described to a first approximation by $n$ real parameters, the
$\ga_j$ in Theorem \ref{up9thm3}, and perhaps equivalently by
$\frac{\pd^m}{\pd x^m}u(0,0)$ for~$m=0,\ldots,n-1$.

Theorem \ref{up9thm3} implies that singularities with
multiplicity $n$ at $(0,0)$ and prescribed values of $\ga_j$
occur in {\it codimension} $2n$ in the family of all singular
solutions $u,v$ of \eq{up4eq2}. Hence, singularities with
multiplicity $n$ at $(0,0)$ but without prescribed $\ga_j$
should occur in {\it codimension} $n$ in the family of all
singular solutions $u,v$ of~\eq{up4eq2}.

More generally, if we consider all solutions of \eq{up4eq2}
and \eq{up4eq3}, and allow singularities anywhere, then we
add one codimension for $a\in\R$, and subtract one because
singularities can occur at $(x,0)$ for any $x\in\R$. This gives:
\medskip

\noindent{\bf Principle.} {\it Singular points with multiplicity
$n\ge 1$ should occur in real codimension $n$ in the family of all
SL\/ $3$-folds invariant under the $\U(1)$-action~\eq{up1eq1}.}
\smallskip

We can relate this to the {\it special Lagrangian fibrations\/}
constructed in \S\ref{up7}. If $\Phi$ is chosen {\it generically}
in Definition \ref{up7def}, or if $\phi$ is chosen {\it generically}
in Example \ref{up7ex}, then we expect $N_{\bs\al}$ to have
singular points of multiplicity $n$ for $\bs\al$ in a
(possibly empty) codimension $n$ subset of $U\subseteq\R^3$.
As $\dim U=3$, it follows that $N_{\bs\al}$ has multiplicity
$n$ singularities for $n=1,2,3$ and $\bs\al$ in a
$(3-n)$-dimensional subset of $U$, and no multiplicity $n$
singularities for~$n>3$.

It is also interesting to consider how the singularities of
$\U(1)$-invariant SL 3-folds studied in \S\ref{up8}--\S\ref{up9}
relate to the broader current state of knowledge on singularities
of SL $m$-folds in (almost) Calabi--Yau $m$-folds, which is
surveyed in \cite{Joyc9}. The class of singular SL $m$-folds
which are best understood theoretically are SL $m$-folds with
{\it isolated conical singularities}. They are the subject of
a series of papers by the author. For a survey and further
references, see~\cite{Joyc8}.

Let $N$ be an SL $m$-fold in $M$ with isolated conical
singularities, and $\bf x$ a singular point of $N$. The
definition \cite[Def.~3.7]{Joyc8} implies that $N$ is an
SL rectifiable current with ${\bf x}\in N^\circ$, and the
(unique) tangent cone $C$ to $N$ at $\bf x$ is an SL cone
in $T_{\bf x}M$ with multiplicity 1 and an {\it isolated\/}
singularity at~0.

From \S\ref{up8}, if $N$ is a $\U(1)$-invariant SL 3-fold in
$\C^3$ and $\bf x$ an isolated singular point of multiplicity
$n\ge 2$, in the sense of Definition \ref{up8def1}, then the
(unique) tangent cone $C$ to $N$ at $\bf x$ is $\Pi^\phi_+\cup
\Pi^\phi_-$ for $\Pi^\phi_\pm$ defined in \eq{up5eq5}, with
multiplicity 1. That is, $C$ is the union of two copies of
$\R^3$ intersecting in $\R$. So $C$ does {\it not\/} have
an isolated singularity at 0, but is singular along~$\R$.

Therefore the singularities of multiplicity $n\ge 2$ above
are examples of isolated singularities of SL 3-folds which
are {\it not\/} isolated conical singularities in the sense
of \cite{Joyc8}. We have been able to understand them well
in this series of papers by assuming $\U(1)$-invariance.
However, there is no known general theory of special
Lagrangian singularities of this kind, similar to the
isolated conical case, {\it without\/} assuming
$\U(1)$-invariance.

Constructing such a theory appears to the author to
be very difficult, because of the lack of good local
models for the singularities. (The singularities
constructed above are not yet suitable as local models,
as we do not understand their asymptotic behaviour
near the singular points.) One interesting feature
of the results above is that they provide a rich class
of fairly generic special Lagrangian singularities,
which are not adequately covered by any known general
analytic theory. This shows we still have a long way
to go in understanding special Lagrangian singularities,
if this is a feasible goal.

\end{document}